\documentclass[11pt,reqno]{amsart}
\usepackage{amssymb}
\usepackage[T1]{fontenc}
\usepackage[utf8]{inputenc}
\usepackage[russian,english]{babel}
\usepackage{graphicx}
\usepackage{amsmath}
\usepackage{tikz}
\usepackage[hidelinks]{hyperref}
\usepackage{lmodern}
\usepackage{mathtools}
\usepackage{mathabx}

\usepackage[spacing=true,kerning=true,babel=true,tracking=true]{microtype}
\usepackage[shortcuts]{extdash}
\usepackage[foot]{amsaddr}
\usepackage[left=1in,right=1in,top=1in,bottom=1in,bindingoffset=0cm]{geometry}
\usepackage[most]{tcolorbox}
\usepackage[inline]{enumitem}
\usepackage{stmaryrd}
\usepackage{bm}
\usepackage{wasysym}
\usepackage{multicol}

\usepackage{algorithm2e}

\usetikzlibrary{graphs, calc}

\usepackage[
    sortcites,
    backend=biber, style=alphabetic, sorting=nyt, maxnames=100,backref=true]{biblatex}
\usepackage{framed}
 {\endMakeFramed}
\definecolor{shadecolor}{gray}{0.9}

\newtheoremstyle{bfnote}%
{}{}%
{\slshape}{}%
{\bfseries}{\bfseries.}%
{ }%
{\thmname{#1}\thmnumber{ #2}\thmnote{ \ep{\normalfont{}#3}}}

\theoremstyle{bfnote}
\newtheorem{theorem}{Theorem}[section]
\newtheorem{lemma}[theorem]{Lemma}
\newtheorem{prop}[theorem]{Proposition}
\newtheorem{cor}[theorem]{Corollary}
\newtheorem{conj}[theorem]{Conjecture}
\newtheorem{claim}[theorem]{Claim}

\theoremstyle{definition}

\newtheorem{definition}[theorem]{Definition}
\newtheorem{remark}[theorem]{Remark}
\newtheorem{question}[theorem]{Question}
\newtheorem*{assumption*}{Assumption}

\makeatletter
\newcommand{\neutralize}[1]{\expandafter\let\csname c@#1\endcsname\count@}
\makeatother

\newenvironment{theocopy}[1]
{\renewcommand{\thetheorem}{\ref{#1}}%
	\neutralize{theorem}\phantomsection
	\begin{theorem}}
	{\end{theorem}}

\newenvironment{conjcopy}[1]
{\renewcommand{\thetheorem}{\ref{#1}}%
	\neutralize{theorem}\phantomsection
	\begin{conj}}
	{\end{conj}}

\newcommand{\Z}{\mathbb{Z}}

\newcommand{\acts}{\mathrel{\reflectbox{$\righttoleftarrow$}}}

\newcommand{\0}{\varnothing}
\newcommand{\set}[1]{\{#1\}}
\newcommand{\N}{{\mathbb{N}}}
\newcommand{\R}{\mathbb{R}}

\renewcommand{\P}{\mathbb{P}}
\newcommand{\E}{\mathbb{E}}
\renewcommand{\epsilon}{\varepsilon}
\newcommand{\eps}{\epsilon}
\renewcommand{\phi}{\varphi}
\renewcommand{\theta}{\vartheta}
\renewcommand{\leq}{\leqslant}
\renewcommand{\geq}{\geqslant}
\newcommand{\defeq}{\coloneqq}

\newcommand{\bemph}[1]{{\normalfont#1}}
\newcommand{\ep}[1]{\bemph{(}#1\bemph{)}}

\newcommand{\emphd}[1]{{\fontseries{b}\selectfont\textsf{#1}}}

\newcommand{\G}{\Gamma}
\newcommand{\pto}{\dashrightarrow}
\newcommand{\dom}{\mathrm{dom}}

\newcommand{\powerset}[1]{\mathcal{P}(#1)}
\newcommand{\symdif}{\mathbin{\triangle}}
\newcommand{\rest}[2]{{{#1}\vert_{#2}}}

\newcommand{\asi}{\mathsf{asi}}

\newcommand{\dist}{\mathsf{dist}}
\newcommand{\fs}[2]{{[{#1}]^{#2}}}
\newcommand{\finset}[1]{\fs{{#1}}{{<\infty}}}

\newcommand{\LOCAL}{$\mathsf{LOCAL}$\xspace}

\newcommand{\inv}{^{-1}}
\newcommand{\blank}{\mathsf{blank}}
\newcommand{\high}{a}
\newcommand{\low}{b}

\newcommand{\algebra}[1]{\mathbb{#1}}

\newcommand{\Sym}[1]{\mathrm{Sym}(#1)}
\newcommand{\level}{\mathsf{level}}

\usepackage{xspace}

\setlist{topsep=3pt,itemsep=3pt}

\title{\sffamily Measurable matchings in unbalanced graphs}
\date{}

\author{Anton~Bernshteyn}
\address{\normalfont (AB) Department of Mathematics, University of California, Los Angeles, CA, USA}
\email{bernshteyn@math.ucla.edu}

\author{Matt Bowen}
\address{\normalfont (MB) 
Mathematical Institute, University of Oxford, Oxford, UK}
\email{matthew.bowen@maths.ox.ac.uk}

\author{Felix Weilacher}
\address{\normalfont (FW) 
Department of Mathematics, University of California, Berkeley, CA, USA}
\email{weilacher@berkeley.edu}

\thanks{AB's research is partially supported by the NSF CAREER grant DMS-2528522 and the Sloan Research Fellowship (2025).  MB is supported by Ben Green's Simons Investigator Grant number 376201. FW is supported by the NSF under award number DMS-2402064.}

\numberwithin{equation}{section}

\newenvironment{scproof}[1][]{\begin{proof}[\textsc{\upshape{Proof}}#1]}{\end{proof}}

\newcommand*{\myproofname}{Proof}
\newenvironment{claimproof}[1][\myproofname]{\begin{proof}[#1]}{\end{proof}}

\usepackage{titlesec}
\usepackage{titletoc}

\titleformat{\section}[block]{\large\bfseries\sffamily}{\thesection.}{1ex}{}
\titleformat{\subsection}[block]{\bfseries\sffamily}{\thesubsection.}{1ex}{}
\titleformat{\subsubsection}[block]{\itshape}{\bfseries\upshape\sffamily\thesubsubsection.}{1ex}{}
\titlespacing*{\section}{0pt}{*3}{*1}
\titlespacing*{\subsection}{0pt}{*3}{*1}
\titlespacing*{\subsubsection}{0pt}{*2.4}{*1}

\titlecontents{section}
 [1em] %
 {\bigskip}
 {\emphd{\thecontentslabel}\hspace{1.02em}\normalfont}
{\hspace*{2.32em}}
 {\,\,\titlerule*[0.77pc]{.}\contentspage}[\smallskip]

 \titlecontents{subsection}
 [3em] %
 {\smallskip}
 {{\thecontentslabel}\hspace{1.02em}\normalfont}
{\hspace*{2.32em}}
 {\,\,\titlerule*[0.77pc]{.}\contentspage}

\renewbibmacro{in:}{}

\renewbibmacro*{volume+number+eid}{%
	\printfield{volume}%
	\setunit*{\addnbspace}
	\printfield{number}%
	\setunit{\addcomma\space}%
	\printfield{eid}}

\DeclareFieldFormat[article]{volume}{\textbf{#1}\space}
\DeclareFieldFormat[article]{number}{\mkbibparens{#1}}

\DeclareFieldFormat{journaltitle}{#1,}
\DeclareFieldFormat[thesis]{title}{\mkbibemph{#1}\addperiod}
\DeclareFieldFormat[article, unpublished, thesis]{title}{\mkbibemph{#1},}
\DeclareFieldFormat[book]{title}{\mkbibemph{#1}\addperiod}
\DeclareFieldFormat[unpublished]{howpublished}{#1, }

\DeclareFieldFormat{pages}{#1}

\DeclareFieldFormat[article]{series}{Ser.~#1\addcomma}

\begin{document}


\maketitle


\begin{abstract}
    Let $G$ be a locally finite 
    multigraph that is 
    bipartite and ``unbalanced,''
    meaning that it has a nontrivial bipartition $(A,B)$ with $\deg(x) > \deg(y)$ for all $x \in A$ and $y \in B$.  
    We explore matchings in such graphs through the lens of descriptive set theory. In particular, we show that when $G$ is Borel and 
    $\mu$ is a Borel probability measure on its vertex set, 
    there is a Borel matching in $G$ that covers $\mu$-almost every vertex in $A$. 
    This was previously known only under the assumption that $\mu$ is $G$-invariant, which we eliminate using a novel probabilistic approach. We also describe various extra conditions that imply the existence of a Borel matching covering every vertex in $A$. Along the way, we confirm a conjecture of the first and third named authors concerning the existence of Borel independent complete sections in Borel graphs of
    finite asymptotic separation index. 
    
    In addition to their intrinsic interest, our results have applications to various other topics, such as 
    edge-colorings, balanced orientations, and equidecomposition theory for group actions. For example, 
    we show 
    that the measurable edge-chromatic number of every Borel multigraph with finite maximum degree $\Delta$ is at most $\lfloor\frac{3\Delta}{2}\rfloor$, matching Shannon's optimal bound for finite multigraphs.  Another example is that paradoxical Borel group actions with finite asymptotic separation index admit paradoxical decompositions with Borel pieces. This refines a result of Marks and Unger.  
\end{abstract}

\renewcommand{\contentsname}{\emphd{Contents}}
\tableofcontents

\section{Introduction}

    \subsection{Matchings and unbalanced bipartitions}

    In this paper we study matchings and their applications from the perspective of \emph{descriptive combinatorics}, a branch of descriptive set theory that investigates measurably or topologically well-behaved 
    solutions to problems in graph theory.  We refer the reader to \cite{kechris, tserunyan2016introduction} for background on descriptive set theory and to \cite{kechris.marks,pikhurko2020borel} for introductions to descriptive combinatorics.  

    Our objects of interest are \emphd{Borel multigraphs}, i.e., triples $G=(V,E,\pi),$ where the \emphd{vertex set} $V$ and the \emphd{edge set} $E$ are standard Borel spaces---such as $[0,1]$ or $2^\mathbb{N}$---and $\pi \colon E \to [V]^2$ is a Borel map sending edges to their pairs of end-vertices. (Here $[V]^2$ denotes the set of all unordered pairs of elements of $V$.) All multigraphs we consider in this paper are loopless (i.e., every edge has two distinct endpoints) and undirected. Unless explicitly stated otherwise, we shall use the words ``graph'' and ``multigraph'' interchangeably. That being said, most of our results are new 
    already when $G$ is a \emphd{simple} graph, i.e., such that the map $\pi \colon E \to [V]^2$ is injective. If $G$ is simple, it is often convenient to 
    assume that $E \subseteq [V]^2$ is a Borel subset 
    and $\pi$ is the inclusion map $E \hookrightarrow [V]^2$. In that case, we omit the map $\pi$ and simply write 
    $G = (V,E)$. 
    
    We use standard graph-theoretic terminology and notation. 
    In particular, 
    $N_G(x)$ and $\deg_G(x)$ denote the neighborhood and the degree (i.e., the number of incident edges) of a vertex $x$, respectively. 
    Note that  when $G$ is simple, $\deg_G(x) = |N_G(x)|$, while in general we have
    $\deg_G(x) \geq |N_G(x)|$. A graph is 
    \emphd{locally finite} if all its vertex degrees are 
    finite. A graph is \emphd{$d$-regular} for $d \in \N$ if all its vertices have degree exactly $d$.

    The study of 
    matchings in Borel graphs is one of the central topics in descriptive combinatorics. 
    A \emphd{matching} in a multigraph $G=(V,E,\pi)$ is a set $M\subseteq E$ of edges such that their endpoint sets are pairwise disjoint. 
    A matching $M$ is \emphd{Borel} if it is a Borel subset of $E$. We say that a matching $M$ \emphd{covers} a set $A \subseteq V$ if for all $x \in A$, there is a (necessarily unique) edge $e \in M$ with $x \in \pi(e)$.

    %
    Borel matchings have numerous applications, 
    for example, to 
    constructive versions of the Banach--Tarski paradox and related problems 
    \cite{gmp,bowen.kun.sabok,marks2016baire,marks-unger,mathe2022circle,ciesla.sabok,bowen2025uniform}, 
    to F\o{}lner tilings of amenable group actions \cite{folner.tilings}, and to the study of stochastic processes 
    \cite{timar.new}. In recent years, 
    Borel matching theory has experienced rapid growth, and several powerful methods for constructing matchings in Borel graphs have emerged. For example, detailed work has been done on matchings in graphs arising from Abelian group actions, where abundant geometric structure can be exploited \cite{bht,grebik.rozhon,weilacher.abelian}. Another fruitful direction employs concepts from large-scale geometry, such as end spaces, to study matchings in hyperfinite graphs \cite{bowen2024measurable,bowen.kun.sabok,conley.miller.pm}. Here we are interested in more general, purely combinatorial techniques. 

    We will specifically focus on matchings in {bipartite} graphs. 
    Recall that a graph $G = (V,E,\pi)$ is \emphd{bipartite} if there exists a partition $V = A \sqcup B$ 
    such that every edge of $G$ has one endpoint in $A$ and the other in $B$; the pair of sets $(A,B)$ is then called a \emphd{bipartition} of $G$. 
    For a Borel graph $G$, a bipartition $(A,B)$ is \emphd{Borel} if $A$ and $B$ are Borel sets. We remark that a bipartite Borel graph may fail to have a Borel bipartition \cite[\S4.3]{kechris.marks}. The theory of matchings in bipartite graphs is particularly well understood, thanks to the following celebrated theorem of Hall \cite{hall1987representatives}: 

    \begin{definition}[Hall's condition]
        A bipartition $(A,B)$ of a locally finite graph $G$ satisfies \emphd{Hall's condition} if for every finite set $S \subseteq A$, we have
        \[
            |N_G(S)| \,\geq\, |S|,
        \]
        where $N_G(S) \defeq \bigcup_{x \in S} N_G(x)$ is the {neighborhood} of $S$ in $G$.
    \end{definition}

    \begin{theorem}[{Hall \cites{hall1987representatives}, Hall~\cite{MHall}; \cite[Thm.~2.1.2, Ex.~8.27]{Diestel}}]\label{theo:Hall}
        A locally finite graph $G$ with a bipartition $(A,B)$ has a matching that covers $A$ if and only if $(A,B)$ satisfies Hall's condition. 
    \end{theorem}

    Although in most graph theory textbooks (and in Philip Hall's original paper \cite{hall1987representatives}), Theorem~\ref{theo:Hall} is only stated for finite graphs, the general locally finite case follows easily by compactness and was proved by Marshall Hall, Jr.~\cite{MHall}. 
    Notice, however, that this argument explicitly relies on the Axiom of Choice and is therefore ``non-constructive.'' In some sense, a fully constructive proof is impossible, as Theorems~\ref{theo:Marks} and \ref{theo:Kun} below demonstrate. In particular, Hall's condition does not generally imply the existence  of a \emph{Borel} matching that covers all, or even almost all, vertices in $A$. Nonetheless, in certain situations large Borel matchings 
    can be obtained under suitably strengthened variants of Hall's condition. In this paper we study one such strengthening: 

    \begin{definition}[Unbalanced bipartitions]\label{defn:unbalanced}
        We say that a bipartition $(A,B)$ of a locally finite graph $G$ is \emphd{of type $(\high,\low)$} for integers $\high$, $\low \in \N$ if
        \[
            \deg_G(x) \,\geq\, \high \text{ for all } x \in A \qquad \text{and} \qquad \deg_G(y) \,\leq\, \low \text{ for all } y \in B.
        \]
        A bipartition $(A,B)$ is \emphd{unbalanced} if it is of type $(\high,\low)$ for some $\high > \low$.
    \end{definition}

    In other words, in an unbalanced bipartition $(A,B)$, we have $\deg_G(x) > \deg_G(y)$ for all $x \in A$ and $y \in B$. Note that unbalanced bipartitions in locally finite Borel graphs are Borel, because $\deg_G$ is a Borel function \cite[Lem.~5.15]{pikhurko2020borel}. A simple double counting argument shows that unbalanced bipartitions satisfy Hall's condition. 
    Our aim in this paper is to show that being unbalanced is a useful substitute for Hall's condition in the Borel setting: on the one hand, it is a strong enough assumption to guarantee the existence of desirable Borel matchings; on the other hand, 
    Borel graphs arising in applications 
    often have unbalanced bipartitions and hence fall into our framework.

    \subsection{Measure theory}\label{subsec:measure}

    Our first result yields measure-theoretically large Borel matchings in unbalanced bipartitions:

     \begin{tcolorbox}
    \begin{theorem}[Measurable matchings]\label{theo:measurable}
        Suppose $(A,B)$ is an unbalanced bipartition of a 
        locally finite Borel graph $G = (V,E,\pi)$. 
        Then, for any probability measure $\mu$ on $V$, $G$ has a Borel matching that covers $\mu$-almost every vertex in $A$.
    \end{theorem}
    \end{tcolorbox}

    Here and in what follows, `measure' means  `Borel measure'. There are several respects in which the statement of Theorem~\ref{theo:measurable} is optimal. First, the loss of a $\mu$-null set in the conclusion of Theorem~\ref{theo:measurable} is necessary. Indeed, while it is natural to seek a Borel matching that covers every vertex in $A$, the determinacy technique of Marks shows that such a matching may not exist: 

        \begin{theorem}[{ess.~Marks \cite{Marks.determinacy}}]\label{theo:Marks}
        For all integers $\high \geq 2$, there exists a simple 
        Borel graph $G$ with a bipartition $(A,B)$ of type $(\high ,2)$ 
            such that no Borel matching in $G$ covers $A$.
        \end{theorem}    

        Although Theorem~\ref{theo:Marks} is not explicitly stated in Marks's paper \cite{Marks.determinacy}, it follows from Marks's analysis of a related combinatorial problem on regular graphs, which in computer science is referred to as \emph{edge grabbing} \cite{btoast}. See \S\ref{subsec:Marks} for the details.

        Second, for Theorem~\ref{theo:measurable} to hold, the inequality $\high > \low$ 
        in Definition~\ref{defn:unbalanced} must be strict. In other words, the theorem fails for bipartitions of type $(\high,\high)$. Indeed, one of the earliest results in modern descriptive combinatorics is the 1988 theorem of Laczkovich \cite{lacz3}, giving an example of a $2$-regular Borel graph $G$ with a Borel bipartition $(A,B)$ (of type $(2,2)$, since $G$ is $2$-regular) and a probability measure $\mu$ such that no Borel matching in $G$ covers $\mu$-almost every vertex in $A$. Laczkovich's example was generalized to $d$-regular graphs for all even $d \geq 2$ by Conley and Kechris \cite{conley.kechris}, while the odd $d$ case was recently settled by Kun \cite{kun2021measurable}:

        \begin{theorem}[{Kun \cite{kun2021measurable}}]\label{theo:Kun}
        For each integer $d \geq 2$, there exists a simple Borel graph $G = (V,E)$ with a Borel bipartition $(A,B)$ and a probability 
        measure $\mu$ on $V$ such that:
        \begin{itemize}
            \item every connected component of $G$ is an infinite $d$-regular tree,
            \item no Borel matching in $G$ covers $\mu$-almost every vertex in $A$.
        \end{itemize}
        \end{theorem}

        Theorem~\ref{theo:Kun} also shows that the unbalanced property of $(A,B)$ in Theorem~\ref{theo:measurable} cannot be replaced by another, less restrictive amplification of Hall's condition, namely combinatorial expansion:

        \begin{definition}[Combinatorial expansion]\label{defn:expansion}
            A bipartition $(A,B)$ of a locally finite graph $G$ has \emphd{combinatorial expansion} if there is $c > 1$ such that for every finite set $S \subseteq A$, we have
            \[
                |N_G(S)| \,\geq\, c\, |S|.
            \]
        \end{definition}

        Clearly, combinatorial expansion implies Hall's condition. Also, if $(A,B)$ is of type $(\high,\low)$ for $\high > \low$, then it has combinatorial expansion with $c = \high/\low$. 
        Combinatorial expansion has been successfully applied to the study of matchings in Borel graphs, notably in the paper \cite{marks2016baire} by Marks and Unger (we describe their results in \S\S\ref{subsec:Baire} and \ref{subsubsec:Folner_intro}). However, since bipartitions of $d$-regular trees with $d \geq 3$ have combinatorial expansion, Theorem~\ref{theo:Kun} demonstrates that combinatorial expansion alone does not suffice for the conclusion of Theorem~\ref{theo:measurable}.\footnote{This can also be shown with other, simpler constructions. See, for example, the discussion at the end of \S\ref{subsec:Baire}.} 
        That being said, some of our results do apply to bipartitions with combinatorial expansion---see \S\ref{subsec:intro_Borel}.
        
        It is important to note that Theorem~\ref{theo:measurable} was known previously under the additional assumption that 
        $\mu$ is \emphd{$G$-invariant}, or, equivalently, $G$ is \emphd{$\mu$-measure-preserving}, meaning that there exists a measure-preserving action $\G \acts (V,\mu)$ of some countable group $\G$ such that for all 
        $\set{x,y} \in \pi(E)$, 
        there is a group element $\gamma \in \G$ with $y = \gamma \cdot x$. Although this is a rather restrictive condition, some of the most familiar examples in descriptive combinatorics arise from ergodic-theoretic considerations and thus are naturally measure-preserving \cite[\S\S6.5--7]{kechris.marks}. For instance, the graphs constructed by Kun in Theorem~\ref{theo:Kun} are measure-preserving.

        In the measure-preserving setting, the proof of Theorem~\ref{theo:measurable} proceeds in two stages.

        \begin{itemize}[wide]
            \item If $(A,B)$ is of type $(\high,\low)$ with $\high > \low$,  
         the $G$-invariance of $\mu$ can be used to show that 
            $(A,B)$ has \emphd{measure expansion}, i.e., every Borel set $S \subseteq A$ satisfies
            $
                \mu(N_G(S)) \geq c \,\mu(S)
            $,
            where $c \defeq \high/\low > 1$. 

            \item An argument due to Lyons and Nazarov \cite{lyons.nazarov}, which also crucially depends on the $G$-invariance of $\mu$, then uses measure expansion to construct the desired matching.
        \end{itemize}

        See \cite[\S2]{folner.tilings} for a detailed presentation of the proof at this level of generality. Unfortunately, the above strategy completely collapses for non-$G$-invariant measures, as neither stage of the proof goes through for non-measure-preserving graphs. (The second stage actually relies on the 
        invariance of measure in two distinct ways, see \cite[Lems.~2.2, 2.5]{folner.tilings}.) Consequently, we must develop a substantially different approach.
        
        While measure-preserving graphs are much better understood than their non-measure-preserving counterparts, 
        there has been a recent surge of results reaching beyond the measure-preserving setting \cite{MCP1,MCP2,MCP3,bowen.poulin,MCP4,MCP5,MCP_prob_1,MCP_prob_2,MCP_prob_3}. 
        Without extra assumptions on $G$ (such as hyperfiniteness), most prior techniques 
        in this research stream 
        rely on 
        controlling the {Radon--Nikodym cocycle} of the measure $\mu$ to quantify its deviation from $G$-invariance 
        \cite{MCP2,MCP5,MCP_prob_2,MCP1,MCP4,MCP_prob_1,MCP_prob_3}. 
    We especially highlight the paper \cite{MCP2} 
    by Greb\'ik, which provided some initial inspiration for our approach. By contrast, 
    our arguments are purely combinatorial and do not utilize any special properties of the measure $\mu$ at all.

        Specifically, we prove Theorem~\ref{theo:measurable} using a \emph{randomized Borel construction}. That is, given a parameter $\omega$ from a certain probability space $(\Omega, \P)$, we define a corresponding Borel matching $M_\omega$ in $G$. We then show that if $\omega$ is chosen at random from $(\Omega,\P)$, then for \emph{every} point $x \in A$, the probability $x$ is covered by $M_\omega$ is $1$. 
        If $\mu$ is a probability measure on $V$, we may then apply Fubini's theorem in the product space $(\Omega,\P) \times (V,\mu)$ to conclude that with probability $1$, $M_\omega$ covers $\mu$-almost all $x \in A$. In particular, a Borel matching covering $\mu$-almost every point in $A$ exists.

        While randomized Borel constructions of a somewhat similar flavor have already been used in descriptive combinatorics, e.g., in \cites[Lem.~5.12]{BernshteynDistributed}[Thm.~1.3]{bowen.kun.sabok}, we think this is a perspective that is worth highlighting and isolating explicitly, as it has the potential for further applications in a non-measure-preserving environment. Furthermore, our analysis leading to Theorem~\ref{theo:measurable} crucially relies on fully embracing the probabilistic outlook.  
        %
        The proof of Theorem~\ref{theo:measurable} is presented in \S\ref{sec:measurable}.
    


        \subsection{Baire category}\label{subsec:Baire}

        In addition to measure theory, results in descriptive combinatorics often invoke {Baire category}. The version of Theorem~\ref{theo:measurable} in the Baire category setting was already proved by Marks and Unger \cite{marks2016baire}, and we record it here for completeness and ease of reference. Let $X$ be a standard Borel space. A \emphd{compatible topology} on $X$ is a Polish topology $\tau$ that generates the Borel $\sigma$-algebra on $X$. Recall that a set $S \subseteq X$ is \emphd{$\tau$-meager} if it is a countable union of $\tau$-nowhere dense sets.

    \begin{theorem}[{ess.~Marks--Unger \cite[Thm.~1.3]{marks2016baire}}]\label{theo:BM}
        Suppose that a Borel bipartition $(A,B)$ of a locally finite Borel graph $G = (V,E,\pi)$ has combinatorial expansion. Then, for any compatible topology $\tau$ on $V$, $G$ has a Borel matching that covers all but a $\tau$-meager set of vertices in $A$.
    \end{theorem}

        In particular, Theorem~\ref{theo:BM} applies when the bipartition $(A,B)$ is unbalanced. 
        We remark that, technically, the statement of
        Theorem~\ref{theo:BM} 
        differs from that of
        \cite[Thm.~1.3]{marks2016baire} in 
        two minor ways. First, Marks and Unger only work with simple graphs, whereas we consider multigraphs (this is not a significant distinction in the context of Theorem~\ref{theo:BM}, because combinatorial expansion is a property of the underlying simple graph of $G$). Second, \cite[Thm.~1.3]{marks2016baire} is a ``two-sided'' variant of Theorem~\ref{theo:BM}: it assumes that combinatorial expansion holds for both $(A,B)$ and $(B,A)$ and concludes that $G$ has a Borel matching covering all but a $\tau$-meager set of vertices in the entire graph. As a result, it is not necessary for the statement of \cite[Thm.~1.3]{marks2016baire} to assume that the bipartition $(A,B)$ is Borel. However, 
        the proof of \cite[Thm.~1.3]{marks2016baire} can be adapted to establish Theorem~\ref{theo:BM} as stated with minimal modifications. 



    As mentioned in \S\ref{subsec:measure}, Kun's Theorem~\ref{theo:Kun} shows that combinatorial expansion is insufficient to obtain a Borel matching covering $\mu$-almost every vertex in $A$. In other words, the obvious analogue of Theorem~\ref{theo:BM} in the measure-theoretic setting fails. 
    Another way to see this is by noting that Theorem~\ref{theo:BM} was used in \cite{marks2016baire} 
    to give a new, simple proof of the result of Dougherty and  Foreman \cite{DF} that the Banach--Tarski paradoxical decomposition of a ball in $\R^3$ can be achieved with Baire-measurable pieces (we discuss this further in \S\ref{subsubsec:Folner_intro}, see Theorem~\ref{theo:MU_paradox}). With Lebesgue-measurable pieces, such a decomposition is 
    obviously impossible. A concrete counterexample to the measure-theoretic version of Theorem~\ref{theo:BM} along these lines is described in \cite[Prop.~16.1]{kechris.marks}.

    \subsection{When can we cover every vertex?}\label{subsec:intro_Borel}

    


    
    In spite of the negative examples in Theorem~\ref{theo:Marks}, we are able to prove several purely Borel results under extra ``large-scale'' assumptions on the graph $G$. This is part of a broad research stream in descriptive combinatorics that aims to capture reasonable conditions under which constructions that normally require the loss of a null set can be performed in a purely Borel fashion; see, e.g., \cite{subexp1,subexp2,subexp3,bernshteyn2023borel,thornton2022orienting,qian2022descriptive,bw_lll,bowen2023definable,conley2020borel,weilacher.abelian,2ends}. Our results 
    employ three of the most prominent ``large-scale'' conditions studied in the above-cited literature. 

    First, we consider Borel graphs of {subexponential growth}. For a graph $G$, a vertex $x$, and $R > 0$, we let $B_G(x,R)$ be the \emphd{$R$-ball} around $x$ in $G$, i.e., the set of all vertices reachable from $x$ by a path of at most $R$ edges. A graph $G$ is \emphd{of subexponential growth} if for all $\epsilon > 0$, there is $R_\epsilon > 0$ such that for all $R \geq R_\epsilon$ and every vertex $x$, we have $|B_G(x,R)| \leq \mathsf{e}^{\epsilon R}$. 
    The usefulness of the subexponential growth assumption in Borel combinatorics is well established \cite{subexp1,subexp2,subexp3,bernshteyn2023borel,thornton2022orienting}, and we show that it is helpful in our setting as well:

    \begin{tcolorbox}
    \begin{theorem}[Graphs of subexponential growth]\label{theo:growth}
        Suppose that a Borel bipartition $(A,B)$ of a locally finite Borel graph $G$ has combinatorial expansion. If $G$ is of subexponential growth, then $G$ has a Borel matching covering $A$.
    \end{theorem}
    \end{tcolorbox}

    Note that our definition of combinatorial expansion is ``asymmetric'': it bounds the neighborhood sizes for finite subsets of $A$ from below, but says nothing about subsets of $B$. This makes the combinatorial expansion condition compatible with subexponential growth.

    Next we consider everywhere $2$-ended graphs. For a graph $G = (V,E, \pi)$ and a set $S \subseteq V$, we let $G - S$ be the graph obtained from $G$ by deleting the vertices in $S$ and their incident edges. The \emphd{number of ends} of a connected locally finite graph $G$ is the supremum, taken over all finite sets $S$ of vertices, of the number of infinite connected components in $G - S$. Given $k \in \N \cup \set{\infty}$, we say that a connected locally finite graph is \emphd{$k$-ended} if it has exactly $k$ ends. A (not necessarily connected) locally finite graph is \emphd{everywhere $k$-ended} if all its components are $k$-ended.

    In descriptive combinatorics, everywhere $k$-ended Borel graphs are interesting when $k$ is either $1$, $2$, or $\infty$. This is because for all other values of $k$, such graphs are \emph{smooth} \cite[Thms.~A, E]{miller.2end}, and hence ``trivial'' as far as Borel combinatorics is concerned \cite[\S5.3]{pikhurko2020borel}. While the cases $k \in \set{1, \infty}$ often pose a particular difficulty, Borel combinatorics tends to behave well on everywhere $2$-ended graphs. As further support for this general philosophy, we prove the following:

    \begin{tcolorbox}
    \begin{theorem}[Everywhere $2$-ended graphs]\label{theo:ends}
        Suppose that a Borel bipartition $(A,B)$ of a locally finite Borel graph $G$ has combinatorial expansion. If $G$ is everywhere $2$-ended, then $G$ has a Borel matching covering $A$.
    \end{theorem}
    \end{tcolorbox}

    Theorems~\ref{theo:growth} and \ref{theo:ends} are proved in \S\ref{sec:Borel_simple}.

    The class of everywhere $2$-ended Borel graphs was vastly generalized by  Conley, Jackson, Marks, Seward, and Tucker-Drob in their landmark paper \cite{conley2020borel}, where---among many other contributions---they introduced the \emph{asymptotic separation index} for locally finite Borel graphs. The asymptotic separation index of $G$, in symbols $\asi(G)$, is a value in $\N \cup \set{\infty}$ that in some sense measures how easily $G$ can be approximated by subgraphs with finite components. The precise definition of $\asi(G)$ is somewhat technical, and we defer it to \S\ref{sec:asi}. For now, we just note that all everywhere $2$-ended Borel graphs have asymptotic separation index $1$, and there exist many other natural classes of Borel graphs with this property, some of which are surveyed in \S\ref{sec:asi}. 
    
    The study of Borel combinatorics under the finite asymptotic separation index assumption has been exceptionally fruitful \cite{qian2022descriptive,bw_lll,bowen2023definable,conley2020borel,weilacher.abelian,Grids2}. We use asymptotic separation index to find matchings in graphs with ``sufficiently high'' combinatorial expansion:

    \begin{tcolorbox}
    \begin{theorem}[Graphs with $\asi < \infty$]\label{theo:asi}
        There is a constant $C > 1$ with the following property.
        Let $G$ be a locally finite Borel graph with a Borel bipartition $(A,B)$ such that for every finite set $S \subseteq A$, we have $|N_G(S)| \geq C\,|S|$. 
        If $\asi(G) < \infty$, then $G$ has a Borel matching covering $A$.
        
        \medskip

        \noindent Furthermore, for graphs $G$ with $\asi(G) \leq 1$, $C = 2$ suffices.
    \end{theorem}
    \end{tcolorbox}

    We also conjecture that, as in Theorems~\ref{theo:BM}, \ref{theo:growth}, and \ref{theo:ends}, any $C > 1$ should be enough:


        \begin{conj}\label{conj:asi}
            Suppose that a Borel bipartition $(A,B)$ of a locally finite Borel graph $G$ has combinatorial expansion. If $\asi(G) < \infty$, then $G$ has a Borel matching covering $A$.
        \end{conj}

    If true, Conjecture~\ref{conj:asi} would generalize Theorem~\ref{theo:BM}, because every locally finite Borel graph on a Polish space has asymptotic separation index at most $1$ on a comeager subset \cite[Thm.~4.8(b)]{conley2020borel}. It would also contain Theorem~\ref{theo:ends} as a special case, since, as mentioned previously, everywhere $2$-ended Borel graphs have asymptotic separation index $1$ \cite[Thm.~1.5]{bowen2023definable}.
    
    Somewhat surprisingly, it is enough to verify Conjecture~\ref{conj:asi} for unbalanced bipartitions:

    \begin{tcolorbox}
    \begin{theorem}[From combinatorial expansion to unbalanced bipartitions when $\asi < \infty$]\label{theo:asi_two_conjectures}
        If Conjecture~\ref{conj:asi} holds whenever the bipartition $(A,B)$ is unbalanced, then it holds in general.
    \end{theorem}
    \end{tcolorbox}

    We stress that for the proof of Theorem~\ref{theo:asi_two_conjectures}, it is crucial that we work with multigraphs. 


    In addition to Theorem~\ref{theo:asi}, we establish a series of other partial results in support of Conjecture~\ref{conj:asi}. To begin with, we show that in the special case of unbalanced bipartitions of {simple} graphs, Theorem~\ref{theo:asi} can be improved by replacing the constant $C$ with a $1 + o(1)$ factor. More generally, this is true for multigraphs in which every pair of vertices is joined by at most $k$ edges for some fixed $k \in \N^+$; we call such graphs \emphd{$k$-simple}. 

    \begin{tcolorbox}
    \begin{theorem}[Almost simple graphs with $\asi < \infty$]\label{theo:asi_simple}
        For each $\epsilon > 0$ and $k \in \N^+$, there exists $b_0 = b_0(\epsilon,k) > 0$ with the following property. 
        Let $G$ be a $k$-simple locally finite Borel graph with a bipartition $(A,B)$ of type $(\high,\low)$, where $\high > (1+\epsilon)\,\low$ and $\low \geq b_0$. If $\asi(G) < \infty$, then $G$ has a Borel matching covering $A$.
    \end{theorem}
    \end{tcolorbox}


    Next, we use Theorem~\ref{theo:asi_simple} to obtain a relaxation of Conjecture~\ref{conj:asi} where each vertex $x \in A$ gets matched to a vertex $y \in B$ that is ``close'' to $x$ but not necessarily adjacent to it. We use $\dist_G(x,y)$ to denote the distance between $x$ and $y$ in $G$, i.e., the minimum length of a path in $G$ joining $x$ to $y$ (if no such path exists, we write $\dist_G(x,y) \defeq \infty$).

    \begin{tcolorbox}
    \begin{theorem}[Bounded distance matchings in graphs with $\asi < \infty$]\label{theo:asi_bounded_distance}
        Suppose that a Borel bipartition $(A,B)$ of a locally finite Borel graph $G$ has combinatorial expansion. If $\asi(G) < \infty$, then there exist $d \in \N$ and a Borel injection $f \colon A \to B$ such that, for all $x \in A$, \[\dist_G(x,f(x)) \,\leq\, d.\] 
    \end{theorem}
    \end{tcolorbox}

    Although Theorem~\ref{theo:asi_bounded_distance} is clearly weaker than Conjecture~\ref{conj:asi}, its conclusion may still fail without the assumption that $\asi(G) < \infty$, even if the bipartition $(A,B)$ is unbalanced and the graph $G$ is simple (Proposition~\ref{prop:no_dist_matching}). Also, Theorem~\ref{theo:asi_bounded_distance} is already sufficient for some interesting applications; for example, it has consequences in equidecomposition theory that are discussed in \S\ref{subsubsec:Folner_intro}.

    Another piece of evidence supporting Conjecture~\ref{conj:asi} is that it holds if the vertices in $B$ have degree at most $2$. We again note that, by Theorem~\ref{theo:Marks}, already in this very simple situation the assumption $\asi(G) < \infty$ cannot be removed.

    \begin{tcolorbox}
    \begin{theorem}[Low degree graphs with $\asi < \infty$]\label{theo:asi_32}
        Suppose that a Borel bipartition $(A,B)$ of a locally finite Borel graph $G$ has combinatorial expansion. 
        If $\asi(G) < \infty$ and $\deg_G(y) \leq 2$ for all $y \in B$, then $G$ has a Borel matching covering $A$.
    \end{theorem}
    \end{tcolorbox}


    The proofs of the above results rely, among other things, on the machinery developed by the first and third named authors in \cite{bw_lll}, which facilitates the use of probabilistic methods and randomized distributed algorithms to perform Borel constructions in a finite $\asi$ environment. In particular, to establish Theorem~\ref{theo:asi}, we prove a conjecture raised in \cite{bw_lll}, namely \cite[Conj.~1.54]{bw_lll}, which concerns the existence of \emph{Borel independent complete sections} in graphs of finite asymptotic separation index. 
    %
    See \S\ref{sec:asi} for the details.

     Finally, we verify Conjecture~\ref{conj:asi} measure-theoretically. Modulo null sets, finiteness of the asymptotic separation index is equivalent to another well-known property, namely \emph{hyperfiniteness} \cite[Cor.~2.1.14, Prop.~2.1.18]{WeilacherThesis} (see Theorem~\ref{theo:asi_ae}). The following result strengthens Theorem~\ref{theo:measurable} in the hyper\-finite setting (as discussed in \S\ref{subsec:measure}, the hyperfiniteness assumption cannot be dropped):

    \begin{tcolorbox}
    \begin{theorem}[Measurable matchings in hyperfinite graphs]\label{theo:hyperfinite}
        Suppose that a Borel bipartition $(A,B)$ of a locally finite Borel graph $G = (V,E,\pi)$ has combinatorial expansion. If $G$ is hyperfinite, then, for any probability measure $\mu$ on $V$, $G$ has a Borel matching that covers $\mu$-almost every vertex in $A$.
    \end{theorem}
    \end{tcolorbox}

    \subsection{Applications}

    \subsubsection{Edge-coloring: Shannon's bound for Borel multigraphs}\label{subsec:edge_col}

    A concept closely related to matchings is \emph{edge-coloring}, and our results on matchings in unbalanced graphs have interesting consequences in edge-coloring theory. Here we focus on edge-colorings with finitely many colors, so all our definitions will be given in that context.
    
    Let $G = (V,E,\pi)$ be a multigraph and let $q \in \N$. A \emphd{$q$-edge-coloring} of $G$ is a function $f \colon E \to q$, where, as usual, we identify the integer $q$ with the $q$-element set $\set{0,1, \ldots, q-1}$. An edge-coloring $f$ is \emphd{proper at a vertex} $x \in V$ if $f(e) \neq f(e')$ for all edges $e \neq e'$ incident to $x$. The \emphd{chromatic index} of $G$, in symbols $\chi'(G)$, is the smallest $q \in \N$ such that $G$ has a $q$-edge-coloring that is proper at every vertex; if no such $q \in \N$ exists, we set $\chi'(G) \defeq \infty$. Equivalently, $\chi'(G)$ is the minimum $q \in \N$ (if it exists) such that the edge set of $G$ can be partitioned into \ep{or covered by} $q$ matchings.    
    
    It is clear that we always have $\chi'(G) \geq \Delta(G)$, where $\Delta(G) \defeq \sup_{x \in V} \deg_G(x)$ 
    is the maximum degree of $G$. 
    On the other hand, a simple greedy algorithm yields the bound 
    $\chi'(G) \leq 2\Delta(G) - 1$. 
    Shannon \cite{shannon1949theorem} showed that this can be improved, roughly by a factor of $3/4$:

    \begin{theorem}[{Shannon \cite{shannon1949theorem}; \cite[Thm.~7.1.13]{West}}]\label{theo:Shannon}
        If $G$ is a graph with $\Delta(G) < \infty$, then \[ \chi'(G) \,\leq\, \left\lfloor \frac{3}{2}\,\Delta(G) \right\rfloor.\]
    \end{theorem}

    Theorem~\ref{theo:Shannon} gives the best upper bound on $\chi'(G)$ in terms of $\Delta(G)$ that is valid for all multigraphs \cite[Ex.~7.1.5]{West}, and we wish to extend this fundamental result to the Borel setting.
    
    Suppose $G = (V,E,\pi)$ is a Borel graph, and let $\mu$ and $\tau$ be a probability measure and a compatible topology on $V$ respectively. The \emphd{Borel}, \emphd{$\mu$-measurable}, and \emphd{$\tau$-Baire-measurable chromatic indices} of $G$---notation: $\chi'_\mathsf{B}(G)$, $\chi'_\mu(G)$, $\chi'_\tau(G)$---are defined as the smallest $q \in \N$ such that $G$ has a Borel $q$-edge-coloring $f \colon E \to q$ that is proper at every vertex, $\mu$-almost every vertex, or all but a $\tau$-meager set of vertices, respectively. Unfortunately, the best general bound on $\chi'_\mathsf{B}(G)$ in terms of $\Delta(G)$ is the greedy one, even for simple graphs $G$:

    \begin{theorem}[{Kechris--Solecki--Todorcevic \cite[15]{kst}, Marks \cite[Thm.~1.4]{Marks.determinacy}}]\label{theo:greedy}
        If $G$ is a Borel graph with $\Delta(G) < \infty$, then $\chi'_\mathsf{B}(G) \leq 2\Delta(G) - 1$. This bound is sharp, even for simple graphs.
    \end{theorem}

    The upper bound in Theorem~\ref{theo:greedy} was proved by Kechris, Solecki, and Todorcevic in their seminal paper \cite{kst}, where they developed a general technique for implementing certain types of greedy algorithms on Borel graphs. The bound's sharpness was demonstrated by Marks using his celebrated determinacy method \cite{Marks.determinacy} (already mentioned in connection with Theorem~\ref{theo:Marks}).

    What about $\chi'_\mu(G)$ and $\chi'_\tau(G)$? While significantly better upper bounds on these parameters are known for simple graphs \cite{MCP2,gp,qian2022descriptive,bernshteyn2023borel}, there has not been any improvement on the greedy bound for general multigraphs with no constraints on edge multiplicities. We remedy this situation by extending Shannon's Theorem~\ref{theo:Shannon} to $\chi'_\mu(G)$ and $\chi'_\tau(G)$:

     \begin{tcolorbox}
    \begin{theorem}[Measurable/Baire-measurable Shannon's bound]\label{theo:Shannon_meas}
        If $G = (V,E,\pi)$ is a Borel graph with $\Delta(G) < \infty$, then for any probability measure $\mu$ and compatible topology $\tau$ on $V$, 
        \[
            \chi'_\mu(G) \,\leq\, \left\lfloor \frac{3}{2}\,\Delta(G) \right\rfloor \qquad \text{and} \qquad \chi'_\tau(G) \,\leq\, \left\lfloor \frac{3}{2}\,\Delta(G) \right\rfloor.
        \]
    \end{theorem}
    \end{tcolorbox}

    Our proof strategy is derived from the work of Ghaffari, Kuhn, Maus, and Uitto in distributed computing \cite{ghaffari2018deterministic}, which reduces the problem to repeatedly finding large matchings in auxiliary bipartite graphs. Exactly the same argument yields Borel versions of Shannon's Theorem~\ref{theo:Shannon} under suitable extra assumptions:

    \begin{tcolorbox}
    \begin{theorem}[Borel Shannon's bound]\label{theo:Shannon_Borel}
        If $G$ is a Borel graph with $\Delta(G) < \infty$, then
        \[
            \chi'_\mathsf{B}(G) \,\leq\, \left\lfloor \frac{3}{2}\,\Delta(G) \right\rfloor,
        \]
        provided that $G$ is of subexponential growth or everywhere $2$-ended.
    \end{theorem}
    \end{tcolorbox}

    If Conjecture~\ref{conj:asi} is true, it would also give the Borel version of Shannon's theorem for graphs with finite asymptotic separation index. We leave it here as an open problem:

    \begin{conj}\label{conj:asi_Shannon}
        If $G$ is a Borel graph with $\Delta(G) < \infty$ and $\asi(G) < \infty$, then
        $
         \displaystyle   \chi'_\mathsf{B}(G) \leq \left\lfloor \frac{3}{2}\,\Delta(G) \right\rfloor
        $.
    \end{conj}

    As partial progress toward Conjecture~\ref{conj:asi_Shannon}, we establish the following weaker bound: 

    \begin{tcolorbox}
    \begin{theorem}\label{theo:Shannon_Borel_asi}
        If $G$ is a Borel graph with $\Delta(G) < \infty$ and $\asi(G) < \infty$, then
        \[
            \chi'_\mathsf{B}(G) \,\leq\, \left\lfloor \frac{3}{2}\,\Delta(G) \right\rfloor \,+\, \asi(G).
        \]
    \end{theorem}
    \end{tcolorbox}

    Theorems~\ref{theo:Shannon_meas}, \ref{theo:Shannon_Borel}, and \ref{theo:Shannon_Borel_asi} are proved in \S\ref{sec:shannon}.

    \subsubsection{Almost balanced orientations}

    Our next application is to the study of \ep{almost} balanced orientations of regular graphs. Let $G$ be a $d$-regular graph, where $d \in \N$. Given an orientation $\mathcal{O}$ of $G$, let the \emphd{defect} of $\mathcal{O}$ at a vertex $x$ be
    \[
        \mathsf{def}_\mathcal{O}(x) \,\defeq\, \left|\deg^+_\mathcal{O}(x) \,-\, \deg^-_\mathcal{O}(x)\right|.
    \]
    Here $\deg^+_\mathcal{O}(x)$ and $\deg^-_\mathcal{O}(x)$ denote the outdegree and the indegree of $x$ in $\mathcal{O}$, respectively. Note that $\mathsf{def}_\mathcal{O}(x)$ has the same parity as $d$. We are interested in orientations with low defect at every vertex. Such objects arise naturally in graph theory and computer science, where they are related to \emph{degree splitting problems} and are useful in divide-and-conquer algorithms; see, e.g., \cite{split1,split2}.

    An orientation $\mathcal{O}$ is \emphd{balanced at a vertex} $x$ if $\mathsf{def}_\mathcal{O}(x) = 0$. Note that this can only happen when $d$ is even. An orientation is \emphd{balanced} if it is balanced at every vertex.\footnote{Balanced orientations are also often called \emph{Eulerian} in the literature.} A simple exercise is to show that if $d$ is even, then every $d$-regular graph has a balanced orientation. The situation becomes more complicated when we look for {Borel} balanced orientations. Indeed, for any $d \geq 1$, Marks \cite{Marks.determinacy} constructed a Borel $d$-regular graph $G$ such that in every Borel orientation of $G$, there is a vertex with defect $d$, and hence $G$ has no Borel balanced orientation. (See \cites[Prop.~8.1]{clp}[Thm.~3.5]{thornton2022orienting} for explicit derivations of this result from Marks's construction.)

    Sacrificing a set of measure $0$ seems to offer no help: for all even $d \geq 2$, Bencs, Hru\v{s}kov\'a, and T\'oth found a $d$-regular Borel graph $G$ with a probability measure $\mu$ such that no Borel orientation of $G$ is balanced at $\mu$-almost every vertex \cite[Thm.~1.3]{bht}. Moreover, their graph is $\mu$-measure-preserving, everywhere $2$-ended, and has subexponential (in fact, linear) growth, which shows that none of these assumptions are helpful as far as balanced orientations are concerned. Additionally, Kun constructed graphs without $\mu$-almost everywhere balanced Borel orientations whose connected components are $d$-regular trees \cite[Cor.~1.6]{kun2021measurable}. That being said, some positive results are known under certain expansion assumptions \cite{thornton2022orienting,bht1,kastner2023baire} and for everywhere one-ended hyperfinite graphs \cite{bht,bowen.poulin,MCP3}.

    Faced with these obstacles, we slightly relax our requirements. Namely, we say that $\mathcal{O}$ is \emphd{almost balanced at a vertex} $x$ if $\mathsf{def}_\mathcal{O}(x) \leq 2$; an orientation is \emphd{almost balanced} if it is almost balanced at every vertex. When $d$ is odd, every vertex in an almost balanced orientation has defect exactly $1$, while in the even case, the defects could be $0$ or $2$. Note that $\mathcal{O}$ is almost balanced at $x$ if and only if
    \[
        \left\lfloor \frac{d-1}{2} \right\rfloor \,\leq\, \deg_\mathcal{O}^+(x) \,\leq\, \left\lceil \frac{d+1}{2} \right\rceil.
    \]
    The problem of bounding the outdegrees in an orientation from above was systematically studied by Thornton \cite{thornton2022orienting}, who in particular established the following positive results:

    \begin{theorem}[{Thornton \cite[Thms.~1.4, 1.5]{thornton2022orienting}}]\label{theo:Thornton}
        Let $G = (V,E,\pi)$ be a $d$-regular Borel graph, where $d \in \N$. Then $G$ has a Borel orientation $\mathcal{O}$ such that $\deg_\mathcal{O}^+(x) \leq \lceil \frac{d+1}{2} \rceil$\ldots
        
        \begin{enumerate}[label=\ep{\normalfont\arabic*}]
            \item\label{item:orient_subexp} \ldots{}for all $x \in V$, provided $G$ is of subexponential growth.

            \item\label{item:orient_meas} \ldots{}for $\mu$-almost all $x \in V$, provided $\mu$ is a $G$-invariant probability measure on $V$.
        \end{enumerate}
    \end{theorem}

     Using the theory of matchings in unbalanced graphs, we are able to considerably strengthen Theorem~\ref{theo:Thornton}: instead of only bounding outdegrees from above, we build almost balanced orientations; in the measure-theoretic context, we do not need to assume that $\mu$ is $G$-invariant; and we extend the results to the Baire category and everywhere $2$-ended settings.

    \begin{tcolorbox}
    \begin{theorem}[Almost balanced orientations]\label{theo:balance}
        Let $G = (V,E,\pi)$ be a $d$-regular Borel graph for $d \in \N$. Let $\mu$ and $\tau$ be a probability measure and a compatible topology on $V$, respectively.
        \begin{enumerate}[label=\ep{\normalfont\arabic*}]
            \item If $G$ is of subexponential growth, then $G$ has a Borel almost balanced orientation.
            
            \item If $G$ is everywhere $2$-ended, then $G$ has a Borel almost balanced orientation.

            \item\label{item:balance_measure} $G$ has a Borel orientation that is almost balanced at $\mu$-almost every vertex.

            \item $G$ has a Borel orientation that is almost balanced at all but a $\tau$-meager set of vertices.
        \end{enumerate}
    \end{theorem}
    \end{tcolorbox}

    All the results in Theorem~\ref{theo:balance} are obtained via the same simple proof that reduces the problem to finding two matchings in an auxiliary bipartite graph. See \S\ref{subsec:orientations} for the details.

    \subsubsection{Equidecomposition theory}\label{subsubsec:Folner_intro}

    We conclude the introduction by discussing another topic closely linked to matchings in bipartite graphs, namely \emph{equidecomposition theory} for group actions. The book \cite{Wagon} by Tomkowicz and Wagon provides an excellent survey of this fascinating subject. The origins of equidecomposition theory date back to the pioneering work of Mazurkiewicz and Sierpi{\'{n}}ski \cite{MS14} and
    Hausdorff \cite{Hausdorff} that led to the discovery of the famous Banach--Tarski paradox \cite{BanachTarski} and the introduction of amenable groups by von Neumann \cite{vNeumann}. The general problem is to determine when one set can be transformed into another one by partitioning it into finitely many pieces and rearranging them. 
        %
        %
    This is formalized in the following definition: 

    \begin{definition}[Equidecomposable sets]\label{defn:equidecomposable}
        Let $\G \acts X$ be an action of a group $\G$ on a set $X$. For a set $A \subseteq X$, a \emphd{$\G$-translate} of $A$ is a set of the form $\gamma \cdot A \defeq \set{\gamma \cdot x \,:\, x \in A}$ for some $\gamma \in \G$. 
        Two sets $A$, $B \subseteq X$ are \emphd{equidecomposable} by $\G$, 
        in symbols $A \approx_{\G} B$, 
        if for some $k \in \N$, there exist partitions 
        \[
            A \,=\, A_1 \sqcup \ldots \sqcup A_k \qquad \text{and} \qquad B \,=\, B_1 \sqcup \ldots \sqcup B_k,
        \]
        where each set $B_i$ is a $\G$-translate of the corresponding set $A_i$. 
        The tuple $(A_1, \ldots, A_k, B_1, \ldots, B_k)$ is called an \emphd{equidecomposition} of $A$ and $B$, and the sets $A_i$, $B_i$ are the \emphd{pieces} of the equidecomposition.
    \end{definition}

    The Banach--Tarski paradox states that for any $n \geq 3$, a ball in $\R^n$ is equidecomposable with two copies of itself by rigid motions \cite[Cor.~3.11]{Wagon}. The idea of ``duplicating'' a set in this way is formally captured by the notion of a \emph{paradoxical decomposition}:

    \begin{definition}[Paradoxical sets]\label{defn:paradoxical}
        Let $\G \acts X$ be an action of a group $\G$ on a set $X$. A set $A \subseteq X$ is \emphd{$\G$-paradoxical} if it has a \emphd{paradoxical decomposition}, i.e., a partition $A = A_0 \sqcup A_1$ into two sets that are both $\G$-equidecomposable 
        with $A$. 
    \end{definition}

    Not all groups $\G$ admit actions with nonempty $\G$-paradoxical sets; the ones that do not are called \emphd{supramenable} \cite[\S14.1]{Wagon}; this term was introduced by Rosenblatt \cite{Rosenblatt}. Finitely generated groups of subexponential growth are supramenable \cite[Thm.~14.21]{Wagon}, and it is an old open question whether there exist any supramenable groups of exponential growth \cite[Q.~14.22]{Wagon}.

    How complicated do the pieces in a paradoxical decomposition need to be? (By that we mean the pieces in the equidecompositions witnessing that $A \approx_\G A_0 \approx_\G A_1$.) This is a natural question that became a focus of considerable interest early on in the development of equidecomposition theory. Clearly, the Banach--Tarski paradox cannot be realized with Lebesgue-measurable pieces. In 1930, Marczewski asked whether Baire-measurable pieces can be used \cite[\S11.2]{Wagon}. (Recall that a set is \emphd{Baire-measurable} if its symmetric difference with some open set is meager.) The positive answer was given by Dougherty and Foreman \cite{DF} and 
    generalized by Marks and Unger as follows: 

    \begin{theorem}[{Marks--Unger \cite[Thm.~1.1]{marks2016baire}}]\label{theo:MU_paradox}
        Let $\G \acts X$ be a Borel action of a group $\G$ on a Polish space $X$ \ep{meaning that for each $\gamma \in \G$, the map $X \to X \colon x \mapsto \gamma \cdot x$ is Borel}. If $X$ is $\G$-paradoxical, then it has a paradoxical decomposition using Baire-measurable pieces.
    \end{theorem}

    Marks and Unger's proof of Theorem~\ref{theo:MU_paradox} relied on matchings in Borel graphs, specifically on Theorem~\ref{theo:BM}, which they proved in the same paper \cite{marks2016baire}. Using our results, namely Theorem~\ref{theo:asi_bounded_distance},\footnote{Theorem~\ref{theo:asi} could also be used instead, although some of the more general results in this subsection, such as Theorem~\ref{theo:Borel_comparison_asi}, do require Theorem~\ref{theo:asi_bounded_distance}.} we can strengthen Theorem~\ref{theo:MU_paradox} as follows:

    \begin{tcolorbox}
    \begin{theorem}[Borel paradoxes under finite $\asi$ assumption]\label{theo:Borel_paradox}
        Let $\G \acts X$ be a Borel action of a finitely generated group $\G$ on a standard Borel space $X$. Suppose that the asymptotic separation index of the action $\G \acts X$ is finite. If a Borel set $A \subseteq X$ is $\G$-paradoxical, then it has a paradoxical decomposition using Borel pieces.
    \end{theorem}
    \end{tcolorbox} 

    See \S\ref{subsec:asi_defns} for the definition of the asymptotic separation index for a finitely generated group action, as well as for some concrete examples of actions to which Theorem~\ref{theo:Borel_paradox} applies, such as the canonical action of a hyperbolic group on its Gromov boundary \cite{NarHyperbolic}. 
    
    Theorem~\ref{theo:Borel_paradox} indeed implies Theorem~\ref{theo:MU_paradox}, because,   
    by \cite[Thm.~4.8(b)]{conley2020borel}, any Borel action of a finitely generated group $\G$ on a Polish space has asymptotic separation index at most $1$ on a comeager $\G$-invariant Borel subset. \ep{Since $\G$-paradoxicality of a set is witnessed by finitely many group elements, it suffices to prove Theorem~\ref{theo:MU_paradox} for finitely generated groups.} 
    Theorem~\ref{theo:Borel_paradox} supports the general philosophy that many results in descriptive combinatorics can be upgraded from measurable/Baire-measurable to Borel when finite $\asi$ is assumed. 

    We can further extend Theorem~\ref{theo:Borel_paradox} by completely describing equidecomposability with Borel pieces for $\G$-paradoxical sets under a finite $\asi$ condition. For that, we need one more definition:

    \begin{definition}[Subequivalence]
        Let $\G \acts X$ be an action of a group $\G$ on a set $X$. A set $A \subseteq X$ is \emphd{subequivalent} to $B \subseteq X$, 
        in symbols $A \preccurlyeq_{\G} B$, 
        if $A$ is equidecomposable by $\G$ with a subset of $B$. 
    \end{definition}

    A basic fact in equidecomposition theory is that $A \approx_\G B$ if and only if $A \preccurlyeq_\G B$ and $B \preccurlyeq_\G A$ \cite[Thm.~3.6]{Wagon}; furthermore, if the relations $A \preccurlyeq_\G B$ and $B \preccurlyeq_\G A$ can be witnessed by equidecompositions with Borel pieces, then $A$ and $B$ are equidecomposable using Borel pieces as well \cite[\S10.3]{Wagon}.

    \begin{tcolorbox}
    \begin{theorem}[Borel equidecompositions for paradoxical sets under finite $\asi$ assumption]\label{theo:Borel_equidecomp}
        Let $\G \acts X$ be a Borel action of a finitely generated group $\G$ on a standard Borel space $X$. Suppose that the asymptotic separation index of the action $\G \acts X$ is finite. Let $A$, $B \subseteq X$ be Borel sets, at least one of which is $\G$-paradoxical. The following statements are equivalent:
        \begin{enumerate}[label=\ep{\normalfont\roman*}]
            \item\label{item:combi_subeq} $A \preccurlyeq_\G B$,
            \item\label{item:Borel_subeq} $A$ is equidecomposable with a Borel subset of $B$ using Borel pieces,
            \item\label{item:cover} $A$ can be covered by finitely many $\G$-translates of $B$.
        \end{enumerate}
        In particular, if $A \approx_\G B$, then there is an equidecomposition of $A$ and $B$ using Borel pieces.
    \end{theorem}
    \end{tcolorbox}

    Theorems~\ref{theo:Borel_paradox} and \ref{theo:Borel_equidecomp} are consequences of an even more general result, which provides a sufficient condition for subequivalence with Borel pieces that reaches beyond paradoxical sets (in particular, it applies to groups of subexponential growth). To state it, we employ a convenient framework for equidecomposition theory developed by Tarski \cite{TarskiSemigroup}; see \cite[\S10]{Wagon} for an in-depth overview. Roughly, Tarski's idea was to extend the relations $\approx_\G$ and $\preccurlyeq_\G$ from subsets of $X$ to their formal linear combinations with natural number coefficients. For example, using this formalism we can say that a set $A \subseteq X$ is $\G$-paradoxical when $A$ is equidecomposable with $2A$, i.e., literally, ``$A$ is equidecomposable with two copies of itself.'' This is made rigorous by encoding the equidecomposability relation for an action $\G \acts X$ in an algebraic structure called a \emph{type semigroup}. We may additionally restrict to equidecompositions that use pieces belonging to some $\G$-invariant subalgebra $\algebra{A} \subseteq \powerset{X}$.  (Here and in what follows, $\powerset{X}$ denotes the powerset of $X$.)
    

    \begin{definition}[Type semigroups]\label{defn:type}
        Let $\G \acts X$ be an action of a group $\G$ on a set $X$ and let $\algebra{A} \subseteq \powerset{X}$ be a $\G$-invariant subalgebra. The \emphd{type semigroup} $\mathsf{S}(X, \G, \algebra{A})$ is the commutative monoid generated by the elements $\set{[A] \,:\, A \in \algebra{A}}$ subject to the relations
        \begin{align*}
            \hspace{4cm}[A] \,&=\, [\gamma \cdot A]  & &\text{for all $A \in \algebra{A}$ and $\gamma \in \G$},\hspace{4cm}\\
            [B \sqcup C] \,&=\, [B] + [C]  & &\text{for all disjoint $B$, $C \in \algebra{A}$}.
        \end{align*}
        We also let
        $\mathsf{S}(X,\G) \defeq \mathsf{S}(X,\G,\powerset{X})$ and call $\mathsf{S}(X,\G)$ the \emphd{full type semigroup} of the action $\G \acts X$. 
        For $\alpha$, $\beta \in \mathsf{S}(X,\G,\algebra{A})$, we write $\alpha \leq \beta$ if there is $\eta \in \mathsf{S}(X,\G,\algebra{A})$ such that $\beta = \alpha + \eta$.
    \end{definition}

    The order relation $\leq$ makes the structure $\mathsf{S}(X, \G, \algebra{A})$ a preordered commutative monoid. Moreover, if $\algebra{A}$ is a $\sigma$-algebra (which is the main case of interest to us), then the relation $\leq$ is antisymmetric \cite[\S10.3]{Wagon}, and thus 
    $\mathsf{S}(X,\G,\algebra{A})$ is a partially ordered commutative monoid.
    
    Subsets $A$, $B \subseteq X$ satisfy $A \approx_\G B$ or $A \preccurlyeq_\G B$ if and only if $[A] = [B]$ or, respectively, $[A] \leq [B]$ in the full type semigroup $\mathsf{S}(X,\G)$. Similarly, equality and the order relation in the type semigroup $\mathsf{S}(X,\G,\algebra{A})$ characterize equidecomposability and subequivalence using pieces from $\algebra{A}$. 

    Note that the inclusion $\algebra{A} \hookrightarrow \powerset{X}$ induces a homomorphism
    $
        \mathsf{S}(X,\G,\algebra{A}) \to \mathsf{S}(X,\G)
    $.
    When no possibility of confusion arises, we shall use the same symbols for elements of $\mathsf{S}(X,\G,\algebra{A})$ and their images in $\mathsf{S}(X,\G)$. We can now state our general result: 

    \begin{tcolorbox}
    \begin{theorem}[Sufficient condition for Borel subequivalence]\label{theo:Borel_comparison_asi}
        Let $\G \acts X$ be a Borel action of a finitely generated group $\G$ on a standard Borel space $X$. Suppose either that $\G$ is of subexponential growth, or that the asymptotic separation index of the action $\G \acts X$ is finite. Let 
        $\alpha$, $\beta \in \mathsf{S}(X,\G,\algebra{B})$, where $\algebra{B}$ is the Borel $\sigma$-algebra of $X$. If for some 
        natural numbers $n > m$,
        the elements $\alpha$, $\beta$ satisfy $n\alpha \leq m\beta$ in $\mathsf{S}(X,\G)$, then $\alpha \leq \beta$ in $\mathsf{S}(X,\G,\algebra{B})$.
    \end{theorem}
    \end{tcolorbox}

    As an illustration, we can quickly derive Theorems~\ref{theo:Borel_paradox} and \ref{theo:Borel_equidecomp} from Theorem~\ref{theo:Borel_comparison_asi}:



    \begin{scproof}[ of Theorem~\ref{theo:Borel_paradox}]
        If $A \subseteq X$ is $\G$-paradoxical, then $[A] = 2[A] = n[A]$ in $\mathsf{S}(X,\G)$ for all $n \in \N^+$. In particular, $2  (2[A]) = [A]$ in $\mathsf{S}(X,\G)$, and thus $2[A] \leq [A]$ in $\mathsf{S}(X,\G,\algebra{B})$, as desired.
    \end{scproof}

    \begin{scproof}[ of Theorem~\ref{theo:Borel_equidecomp}]
        Implications \ref{item:Borel_subeq} $\Longrightarrow$ \ref{item:combi_subeq} $\Longrightarrow$ \ref{item:cover} are trivial. For \ref{item:cover} $\Longrightarrow$ \ref{item:Borel_subeq}, suppose $A$ can be covered by finitely many $\G$-translates of $B$. This means that $[A] \leq n[B]$ in $\mathsf{S}(X,\G)$ for some $n \in \N^+$. If $A$ is $\G$-paradoxical, then $(n+1)[A] = [A] \leq n[B]$ in $\mathsf{S}(X,\G)$, while if $B$ is $\G$-paradoxical, then $2[A] \leq 2n[B] = [B]$ in $\mathsf{S}(X,\G)$. In both cases, we 
        get $[A] \leq [B]$ in $\mathsf{S}(X,\G,\algebra{B})$, as desired.
    \end{scproof}




   Theorem~\ref{theo:Borel_comparison_asi} is closely related to the algebraic property of \emph{almost unperforation}, 
%
   %
    which was introduced by R\o{}rdam \cite{Rordam}: 

    \begin{definition}[Almost unperforation]\label{defn:almost_unperforated}
        A preordered commutative monoid $\mathsf{S}$ is \emphd{almost unperforated} if for all $\alpha$, $\beta \in \mathsf{S}$ and natural numbers $n > m$, $n\alpha \leq m\beta$ implies $\alpha \leq \beta$. 
    \end{definition}
    
    Almost unperforation is widely used in topological dynamics and $C^*$-algebra theory, where it is a major ingredient in the classification program \cite{Strung}. 
    Other relevant notions in those areas are 
    \emph{dynamical comparison} and \emph{almost finiteness}; some references on these topics are \cite{Kerr, KerrSz,DZ_comparison,KerrNar,NarAM,conley2020borel,folner.tilings,AraExel,MoreTypeSemigroup,RordamS,EvenMoreTypeSemigroup,Rordam,Matui,Ma1,Ma2,Melleray,TypeSemigroupGroupoids} (this list is far from exhaustive). For introductions to this circle of ideas with an emphasis on topological dynamics, see the papers \cite{Kerr} by Kerr, \cite{Ma1,Ma2} by Ma, and \cite{Melleray} by Melleray. We also recommend \cite{TypeSemigroupGroupoids} by Kwa{\'{s}}niewski, Meyer, and Prasad and \cite{WehrungBook} by Wehrung for overviews of the purely algebraic aspects of the theory.

    Our methods show that certain type semigroups are almost unperforated:

    \begin{tcolorbox}
    \begin{theorem}[Almost unperforation]\label{theo:unperf}
        Let $\G \acts X$ be an action of a finitely generated group $\G$ on a standard Borel space $X$ and let $\algebra{A}$ be a $\G$-invariant subalgebra of $\mathcal{P}(X)$. Then  the type semigroup $\mathsf{S}(X,\G,\algebra{A})$ is almost unperforated, provided any of the following holds:
        \begin{enumerate}[label=\ep{\normalfont\arabic*}]
            \item\label{item:unperf_meas} $\algebra{A}$ is the $\sigma$-algebra of $\mu$-measurable sets for some probability measure $\mu$ on $X$,
            \item\label{item:unperf_BM} $\algebra{A}$ is the $\sigma$-algebra of $\tau$-Baire-measurable sets for some compatible topology $\tau$ on $X$,
            \item\label{item:unperf_subexp} $\algebra{A}$ is the Borel $\sigma$-algebra of $X$ and $\G$ is of subexponential growth, or
            \item\label{item:unperf_asi} $\algebra{A}$ is the Borel $\sigma$-algebra of $X$ and $\asi(\G \acts X) < \infty$. 
        \end{enumerate}
    \end{theorem}
    \end{tcolorbox} 

    Parts \ref{item:unperf_BM}, \ref{item:unperf_subexp}, and \ref{item:unperf_asi} of Theorem~\ref{theo:unperf} are consequences of Theorem~\ref{theo:Borel_comparison_asi} (part \ref{item:unperf_BM} follows thanks to \cite[Thm.~4.8(b)]{conley2020borel}---see \S\ref{subsubsec:unperf} for details). The proof of part~\ref{item:unperf_meas} instead relies on Theorem~\ref{theo:measurable}, i.e., the existence of almost everywhere matchings in unbalanced graphs.

    \begin{remark}
        While Theorem~\ref{theo:unperf}\ref{item:unperf_meas} can be applied in the special case when the action $\G \acts (X,\mu)$ is measure-preserving, its full scope is much more general. In particular, Theorem~\ref{theo:unperf}\ref{item:unperf_meas} covers the situation when the action $\G \acts X$ is Borel and the measure $\mu$ is \emphd{$\G$-quasi-invariant}, i.e., the $\sigma$-ideal of $\mu$-null sets is preserved by the action. In practice, one can typically assume that the measure is $\G$-quasi-invariant without any loss of generality; see \S\ref{subsubsec:qi} for a further discussion.
    \end{remark}

    \begin{remark}
    There exist free Borel actions $\G \acts X$ such that the type semigroup $\mathsf{S}(X,\G,\algebra{B})$ corresponding to the Borel $\sigma$-algebra $\algebra{B}$ of $X$ fails to be almost unperforated (Proposition~\ref{prop:no_borel_unperf}). Hence, the additional assumptions in parts \ref{item:unperf_subexp} and \ref{item:unperf_asi} of Theorem~\ref{theo:unperf} cannot be dropped.
    \end{remark}

    Lastly, we mention another subject closely related to our results in equidecomposition theory, namely the study of the geometry of 
    amenable groups and their actions. 
    This research direction was initiated by Ornstein and Weiss, who in their landmark paper \cite{OW} showed that every countable amenable group $\G$ can be ``mostly'' covered by pairwise disjoint translates of finitely many finite sets $F_1$, \ldots, $F_n \subseteq \G$ with prescribed approximate invariance. They also extended this result to probability measure-preserving actions of $\G$. This has become known as the \emph{Ornstein--Weiss quasi-tiling theory} and is an important tool in the study of amenable groups and their dynamics.

    A natural question is whether the Ornstein--Weiss quasi-tiling machinery can be improved by not leaving any uncovered points at all. For tilings of $\G$ itself, this was achieved by Downarowicz, Huczek, and Zhang \cite{DHZ_group}. Their result was followed by a paper of Conley, Jackson, Kerr, Marks, Seward, and Tucker-Drob \cite{folner.tilings} in the probability measure-preserving context, where a measure $0$ set of uncovered points is allowed. The paper \cite{folner.tilings} explicitly relied on matchings in unbalanced graphs, specifically on the version of Theorem~\ref{theo:measurable} for $G$-invariant measures. The measure-preservation requirement was eliminated by Conley, Jackson, Marks, Seward, and Tucker-Drob \cite{conley2020borel}, who built Borel F\o{}lner tilings under a finite $\asi$ assumption (which, for actions of amenable groups, is satisfied almost everywhere \cite[\S3.1.3]{WeilacherThesis}).

    As shown by Downarowicz and Zhang \cite{DZ_comparison} and Kerr and Szab\'o \cite{KerrSz}, the existence of F\o{}lner tilings can be characterized via the \emph{comparison property}, which establishes certain relations in the type semigroup using the notion of \emph{Banach density} (the papers \cite{DZ_comparison,KerrSz} focus on continuous constructions, but their arguments can be applied in the Borel setting as well). 
    It is easy to see that the comparison property in the sense of \cite{DZ_comparison} is implied by the almost unperforation of the corresponding type semigroup. Consequently, Theorem~\ref{theo:unperf} can be seen as a strengthening of some earlier results in this field. For example, it provides a new proof of the theorem from \cite{conley2020borel} mentioned above. We shall briefly discuss this aspect of our work in \S\ref{subsubsec:tiling}.

    We find it intriguing that two essentially ``orthogonal'' phenomena---paradoxical decompositions and F\o{}lner tilings---can arise as manifestations of the same underlying mechanism. This is yet another indication of the power of matching theory in descriptive combinatorics and, more specifically, of the range of situations where matchings in unbalanced graphs appear naturally. 

    \subsection*{Outline of the remainder of the paper}


    After some basic preliminaries in \S\ref{sec:prelim}, we proceed to prove Theorem~\ref{theo:measurable} (the existence of a Borel matching covering almost every vertex) using a randomized Borel construction in \S\ref{sec:measurable}.  Then, in \S\ref{sec:Borel_simple}, we study matchings in graphs of subexponential growth and in everywhere $2$-ended graphs and establish Theorems~\ref{theo:growth} and \ref{theo:ends}. In \S\ref{sec:asi}, we formally introduce asymptotic separation index and prove Theorems~\ref{theo:asi}, \ref{theo:asi_two_conjectures}, \ref{theo:asi_simple}, \ref{theo:asi_bounded_distance}, \ref{theo:asi_32}, and \ref{theo:hyperfinite} using tools such as the Lov\'asz Local Lemma and distributed algorithms. We then discuss applications of our main results in \S\ref{sec:applications}. The paper is concluded with a selection of some remaining open problems in \S\ref{sec:problems}.

    \section{Preliminaries}\label{sec:prelim}

    \subsection{Some graph-theoretic notation and terminology}\label{subsec:graphs}

    As mentioned in the introduction, we use standard graph-theoretic terminology. Here we only review a few concepts that are either less well known or particularly prominent in the rest of the paper.

    A \emphd{subgraph} of a graph $G = (V,E,\pi)$ is a graph $G' = (V', E', \pi')$ such that $V' \subseteq V$, $E' \subseteq E$, and $\pi'$ is the restriction of $\pi$ to $E'$. When $G$ and $G'$ are Borel graphs, we will additionally require the inclusions $V' \subseteq V$, $E' \subseteq E$ to be Borel. To indicate that $G'$ is a subgraph of $G$, we write $G' \subseteq G$. A subgraph of $G$ is \emphd{spanning} if its vertex set is $V$.

    As usual, for a set $X \subseteq V$, $G[X]$ is the \emphd{subgraph of $G$ induced by $X$}, i.e., the graph obtained from $G$ by deleting the vertices in $V \setminus X$ and their incident edges. In particular, $G - X = G[V \setminus X]$. For disjoint sets $X$, $Y \subseteq V$, we let $G[X,Y]$ be the \emphd{bipartite subgraph of $G$ induced by $(X,Y)$}, i.e., the graph obtained from $G[X \cup Y]$ by only keeping the edges that join $X$ to $Y$. For a set $F \subseteq E$ of edges, we let $G - F$ be the spanning subgraph of $G$ with edge set $E \setminus F$.

    Recall that for vertices $x$, $y \in V$, $\dist_G(x,y)$ is the graph distance from $x$ to $y$, i.e., the minimum length of a path in $G$ joining $x$ to $y$; if no such path exists, we let $\dist_G(x,y) \defeq \infty$. For a vertex $x$ and sets $X$, $Y \subseteq V$, we let $\dist_G(x,Y) \defeq \inf_{y \in Y} \dist_G(x,y)$ and $\dist_G(X,Y) \defeq \inf_{x \in X, y \in Y} \dist_G(x,y)$.


    A set $X \subseteq V$ is called \emphd{$G$-invariant} if it is a union of connected components of $G$, i.e., if there are no edges joining $X$ to $V \setminus X$. A set $X \subseteq V$ is \emphd{independent} if there are no edges in $G[X]$.

    The \emphd{underlying simple graph} of a multigraph $G = (V,E,\pi)$ is the simple graph with vertex set $V$ where $\set{x,y}$ is an edge if and only if $\set{x,y} \in \pi(E)$. Note that $G$ has a matching that covers some set $A \subseteq V$ if and only if so does its underlying simple graph.

    The \emphd{$k$-th power} of a simple graph $G = (V,E)$, in symbols $G^k$, is the simple graph with vertex set $V$ in which $\set{x,y}$ is an edge if and only if $0 < \dist_G(x,y) \leq k$.

     \subsection{The Luzin--Novikov theorem and its consequences}\label{subsec:LN}

     The most important general tool in working with locally finite Borel graphs is the \emph{Luzin--Novikov theorem} \cite{Luzin,Novikov}:

     \begin{theorem}[{Luzin--Novikov \cite[Thm.~18.10]{kechris}}]\label{theo:LN}
        Let $X$ and $Y$ be standard Borel spaces and let $A \subseteq X \times Y$ be a Borel subset. If for all $x \in X$, the fiber $A_x \defeq \set{y \in Y \,:\, (x,y) \in A}$ is countable, then there exists a sequence $(f_n)_{n \in \N}$ of Borel partial functions $f_n \colon X \pto Y$, where $\dom(f_n)$ is a Borel subset of $X$, such that for all $x \in X$, $A_x = \set{f_n(x) \,:\, n \in \N,\ x \in \dom(f_n)}$.
    \end{theorem}

    The Luzin--Novikov theorem can be used to show that constructions that only involve quantifiers ranging over countable sets produce Borel results. For example, it implies that for a locally finite Borel graph $G$, the functions  $\deg_G \colon V(G) \to \N$ and $\dist_G \colon V(G)^2 \to \N \cup \set{\infty}$ and the underlying simple graph of $G$ are Borel. Also, if $G$ is simple, the Luzin--Novikov theorem shows that $G^k$ is Borel for all $k \in \N$. For a detailed discussion and further examples of using Theorem~\ref{theo:LN} in Borel graph theory, see \cites[\S3.1]{Grids1}[\S5]{pikhurko2020borel}. Routine arguments involving the Luzin--Novikov theorem will be employed in the sequel without mention.
    

    Let us now record a few other specific consequences of Theorem~\ref{theo:LN} that will be used later.

    \begin{lemma}[Luzin--Novikov Uniformization]\label{lemma:LNU}
        Suppose that in the setting of Theorem~\ref{theo:LN}, $A_x \neq \0$ for all $x \in X$. Then there exists a Borel map $\mathsf{choice} \colon X \to Y$ such that for all $x \in X$, $\mathsf{choice}(x) \in A_x$.
    \end{lemma}
    \begin{scproof}
        Let $(f_n)_{n \in \N}$ be given by Theorem~\ref{theo:LN} and define \[
            \mathsf{choice}(x) \,\defeq\, f_{n}(x), \text{ where $n \in \N$ is minimum such that $x \in \dom(f_n)$}. \qedhere
        \]
    \end{scproof}

    Informally, Lemma~\ref{lemma:LNU} says that if each point $x \in X$ is given a countable set $A_x$ of options to choose from, then all the choices can be made in a Borel way.

    \begin{lemma}\label{lemma:enum}
        Suppose that in the setting of Theorem~\ref{theo:LN}, $|A_x| = \high \in \N^+$ for all $x \in X$. Then there exists a Borel map $\mathsf{list} \colon \set{1,\ldots, \high} \times X \to Y$ such that for all $x \in X$, $A_x = \set{\mathsf{list}(1,x), \ldots, \mathsf{list}(\high,x)}$.
    \end{lemma}
    \begin{scproof}
        Fix a Borel linear ordering $\leq$ on $Y$ (for instance, we may assume that $Y$ is a Borel subset of $\R$ by \cite[Thm.~15.6]{kechris} and use the standard ordering on $\R$) and let $\mathsf{list}(i,x)$ be the $i$-th element of $A_x$ in this ordering. \ep{The Luzin--Novikov theorem ensures that the function $\mathsf{list}$ is Borel.}
    \end{scproof}


    \begin{lemma}\label{lemma:reg_subgraph}
        Let $G=(V,E,\pi)$ be a locally finite Borel graph and let $A\subseteq V$ be a Borel independent set in $G$ such that every vertex in $A$ has degree at least $\high \in \N^+$. Then there exists a spanning Borel subgraph $G' = (V,E',\pi')$ in which 
        every vertex of $A$ has degree exactly $\high$. 
        Furthermore, if $M \subseteq E$ is a Borel matching in $G$, it is possible to choose $G'$ so that $M \subseteq E'$.
    \end{lemma}
    \begin{scproof}
        Fix a Borel linear ordering $\leq$ on $E$. For $x \in A$, define $\high_x \in \set{\high, \high-1}$ to be $\high$ minus the number of edges in $M$ incident to $x$. Let $F_x$ be the set of the first $\high_x$ edges in $E \setminus M$ incident to $x$ in the ordering $\leq$. Setting $E' \defeq M \cup \set{e \in E \,:\, e \in F_x \text{ for some } x \in A}$ finishes the proof.
    \end{scproof}

     \subsection{More on Borel independent sets}\label{subsec:chi}

     Another fundamental result in the theory of locally finite Borel graphs is the following theorem of Kechris, Solecki, and Todorcevic:

     \begin{theorem}[{Kechris--Solecki--Todorcevic \cite[Props.~4.5, 4.6, 4.2]{kst}}]\label{theo:KST}
         Let $G = (V,E,\pi)$ be a locally finite Borel graph. 

         \begin{enumerate}[label=\ep{\normalfont\roman*}]
             \item\label{item:loc_fin} There exists a partition $V = \bigsqcup_{n \in \N} U_n$ into countably many Borel independent sets.

             \item\label{item:max_deg} If the maximum degree $\Delta(G)$ is finite, then $V$ can be partitioned into finitely many Borel independent sets, namely at most $\Delta(G) + 1$.

             \item\label{item:maximal} There exists a Borel inclusion-maximal independent set $U \subseteq V$.
         \end{enumerate}
     \end{theorem}

     The minimum cardinal $k$ such that the vertex set of $G$ can be partitioned into $k$ Borel independent sets is called the \emphd{Borel chromatic number} of $G$ and is denoted by $\chi_\mathsf{B}(G)$.\footnote{Technically, this is called the \emph{weak Borel chromatic number} in the literature, but this distinction disappears when $\chi_\mathsf{B}(G)$ is countable. See \cite{Geschke} for more details.} Hence, parts \ref{item:loc_fin} and \ref{item:max_deg} of Theorem~\ref{theo:KST} can be concisely stated as follows: $\chi_\mathsf{B}(G) \leq \aleph_0$ for all locally finite Borel graphs and furthermore $\chi_\mathsf{B}(G) \leq \Delta(G) + 1$ if $\Delta(G) < \infty$. Both bounds are in general optimal even if every component of $G$ is a tree \cites[Ex.~3.2]{kst}{Marks.determinacy}.



     A set $U \subseteq V$ in a graph $G = (V,E,\pi)$ is \emphd{$R$-independent}\footnote{Other common terms for this notion are \emph{$R$-discrete} and \emph{$R$-sparse}.} if $\dist_G(x,y) > R$ for all distinct $x$, $y \in U$.  The following is a useful consequence of Theorem~\ref{theo:KST}:

     \begin{cor}\label{corl:spaced}
         Let $G = (V,E,\pi)$ be a locally finite Borel graph and let $R \in \N$.

         \begin{enumerate}[label=\ep{\normalfont\roman*}]
             \item There is a partition $V = \bigsqcup_{n \in \N} U_n$ into countably many Borel $R$-independent sets.

             \item\label{item:spaced_deg} If $\Delta(G) < \infty$, then $V$ can be partitioned into finitely many Borel $R$-independent sets.

             \item\label{item:max_spaced} There exists a Borel inclusion-maximal $R$-independent set $U \subseteq V$.
         \end{enumerate}
     \end{cor}
     \begin{scproof}
         Apply Theorem~\ref{theo:KST} to the $R$-th power of the underlying simple graph of $G$.
     \end{scproof}

     \subsection{Smooth graphs}\label{subsec:smooth}

     A locally finite Borel graph $G = (V,E,\pi)$ is \emphd{smooth} if it admits a \emphd{Borel transversal}, i.e., a Borel set $T \subseteq V$ that meets every connected component of $G$ in exactly one point. The following is a useful equivalent characterization of smoothness (which also provides the correct way to define smoothness for arbitrary, i.e., not necessarily locally finite or locally countable, Borel graphs):

     \begin{prop}[{\cite[Prop.~20.5]{tserunyan2016introduction}}]
         A locally finite Borel graph $G = (V,E,\pi)$ is smooth if and only if there exists a Borel map $f \colon V \to X$ to some standard Borel space $X$ such that $f(x) = f(y)$ if and only if $x$ and $y$ are in the same component of $G$.
     \end{prop}


     Subgraphs of locally finite smooth graphs are smooth:

     \begin{prop}[{\cite[Prop.~1.6.2]{smooth}}]\label{prop:smooth_sub}
         If $G$ is a smooth locally finite Borel graph, then all Borel subgraphs of $G$ are smooth.
     \end{prop}

    The following lemma will simplify checking that a bipartite graph is smooth:

     \begin{lemma}\label{lemma:smooth_bipartite}
         Let $G$ be a simple locally finite Borel graph with a Borel bipartition $(A,B)$. Then $G$ is smooth if and only if $G^2[A]$ is smooth.
     \end{lemma}
     \begin{scproof}
         Since $G$ and $G^2$ have the same connected components, if $G$ is smooth, then $G^2$ is smooth as well, and hence $G^2[A]$ is also smooth by Proposition~\ref{prop:smooth_sub}. Conversely, if $T \subseteq A$ is a Borel transversal for $G^2[A]$, then $T \cup \set{y \in B \,:\, \deg_G(y) = 0}$ is a Borel transversal for $G$.
     \end{scproof}

     A graph $G$ is \emphd{component-finite} if all its connected components are finite.

     \begin{prop}[{\cite[Lem.~5.21]{pikhurko2020borel}}]\label{prop:comp_fin_smooth}
         Every component-finite Borel graph is smooth.
     \end{prop}

     From the point of view of Borel combinatorics, smooth graphs are ``trivial,'' in the sense that to find a Borel solution to a combinatorial problem on such a graph, it suffices to show that a solution exists abstractly, i.e., with no regard for it being Borel. This is true in a broad context of so-called \emph{locally checkable labeling \ep{LCL} problems} (sometimes also called \emph{local coloring problems}), as defined, for example, in \cites{Ber_smooth}[\S5.3]{pikhurko2020borel}. (See also Definition~\ref{defn:CSP} in \S\ref{subsec:asi_general} for the essentially equivalent notion of \emph{constraint satisfaction problems}.) Finding a matching that covers a given set $A$ of vertices is an example of such a problem---the task is to label the edges with $0$s and $1$s so that no two edges with a common endpoint are labeled $1$, and each vertex in $A$ is incident to at least one edge labeled $1$---and hence we have the following:

     \begin{theorem}[{\cites[Thm.~1.2]{Ber_smooth}[Thm.~5.23]{pikhurko2020borel}}]\label{theo:smooth}
         Let $G = (V,E, \pi)$ be a locally finite Borel graph and let $A \subseteq V$ be a Borel set. If $G$ is smooth and there is a matching in $G$ covering $A$, then there is a Borel such matching.
     \end{theorem}

     Theorem~\ref{theo:smooth} is a folklore result. Two different proofs of its more general form for arbitrary LCL problems can be found in the papers \cite{Ber_smooth,pikhurko2020borel}. The proof in \cite{pikhurko2020borel} is more combinatorial, while the one in \cite{Ber_smooth} is shorter but relies on more advanced descriptive set theory machinery. 

      \subsection{Proof of Theorem~\ref{theo:Marks}}\label{subsec:Marks}

      To finish this preliminary section, we sketch the derivation of Theorem~\ref{theo:Marks} from Marks's results in \cite{Marks.determinacy}. Let $G = (V,E,\pi)$ be a graph. In the \emphd{edge grabbing problem} on $G$, each vertex $x \in V$ must ``grab'' one of its incident edges so that no two vertices grab the same edge. Formally, a solution to the edge grabbing problem is an injective function $f \colon V \to E$ such that $f(x)$ is incident to $x$ for all $x \in V$. As explained in \cite{btoast}, the proof of \cite[Thm.~1.3]{Marks.determinacy} yields the following:

      \begin{theorem}[{Marks \cite{Marks.determinacy}; see \cite{btoast}}]\label{theo:grab}
          For all $d \in \N$, there exists a simple $d$-regular Borel graph $G$ with no Borel solution to the edge grabbing problem. Furthermore, every component of $G$ is a $d$-regular tree and its Borel chromatic index is $d$. 
      \end{theorem}

      Now fix an integer $\high \geq 2$ and let $G = (V,E)$ be the graph given by Theorem~\ref{theo:grab} with $d = \high$. Let $G^*$ be the simple graph with vertex set $V \cup E$ (where we assume the union is disjoint) 
      and edge set
      \[
        \big\{\set{x,e} \,:\, x \in V,\, e \in E,\, \text{$x$ is incident to $e$ in $G$}\big\}.
      \]
      Note that $(V,E)$ is a bipartition of $G^*$ of type $(\high, 2)$. If $M$ is a matching in $G^*$ that covers $V$, then the map $f \colon V \to E$ such that $\set{x, f(x)} \in M$ for all $x \in V$ is a solution to the edge grabbing problem on $G$. It follows that $G^*$ has no Borel matching covering $V$, as desired. 


    \section{Matchings almost everywhere: Theorem~\ref{theo:measurable}}\label{sec:measurable}

    \subsection{Random Borel matchings}

    The goal of \S\ref{sec:measurable} is to prove Theorem~\ref{theo:measurable}, i.e., the existence of a Borel matching covering $\mu$-almost every vertex in $A$, where $(A,B)$ is an unbalanced bipartition of a locally finite Borel graph $G = (V,E,\pi)$. As mentioned in the introduction, 
        %
    we build the desired matching via a \emph{randomized Borel construction}. 
    
    By a \emphd{random Borel matching} in a Borel graph $G = (V,E,\pi)$ we mean a Borel set $M \subseteq \Omega \times E$, where $(\Omega, \P)$ is a standard probability space, such that for all $\omega \in \Omega$, the fiber
    \[
        M_\omega \,\defeq\, \set{e \in E \,:\, (\omega, e) \in M}
    \]
    is a matching in $G$. We refer to $(\Omega, \P)$ as the \emphd{parameter space} of $M$. 
    Of course, although our 
    focus is on matchings, one can also 
    define random versions of various other Borel structures. 

    To reduce notational clutter, we will make use of the standard convention in probability theory and omit the explicit reference to the point in the parameter space whenever that does not lead to confusion. For example, for a random Borel matching $M \subseteq \Omega \times E$, we will use expressions such as ``the probability $M$ has some property $\mathcal{P}$,'' in symbols $\P[\text{$M$ has $\mathcal{P}$}]$, 
    to mean $\P[\set{\omega \in \Omega \,:\, \text{$M_\omega$ has $\mathcal{P}$}}]$. Also, to emphasize the probabilistic perspective, we shall always use $\P$ to denote the measure on the parameter space, while integration with respect to $\P$ will be denoted by $\E$ (for ``expectation''). Basic probabilistic concepts, such as conditional probability, will be used freely.

    We shall prove the following 
    strengthening of Theorem~\ref{theo:measurable}: 

    \begin{theorem}\label{theo:random}
        Let $(A,B)$ be an unbalanced bipartition of a 
        locally finite Borel graph $G = (V,E,\pi)$. Then there exists a random Borel matching $M$ in $G$ such that for all $x \in A$, $\P \left[\text{$M$ covers $x$}\right] = 1$.
    \end{theorem}

    \begin{scproof}[ of Theorem~\ref{theo:measurable}]
        Let $M \subseteq \Omega \times E$ be given by Theorem~\ref{theo:random}. By Fubini's theorem,
        \[
            \E\left[\mu\left(\set{x \in A \,:\, \text{$M$ covers $x$}}\right)\right] \,=\, \int_A \P[\text{$M$ covers $x$}] \,\mathrm{d}\mu(x) \,=\, \mu(A),
        \]
        and hence $M$ covers $\mu$-almost all $x \in A$ with probability $1$. In particular, there exists some $\omega \in \Omega$ such that $M_\omega$ is a Borel matching that covers $\mu$-almost every vertex in $A$, as desired.
    \end{scproof}

    To prove Theorem~\ref{theo:random}, we iteratively generate random Borel matchings with incrementally growing sets of covered vertices.  Individual steps of the iteration are handled by 
    Lemma~\ref{lemma:induction} below, which contains the technical core of the proof. 
    To state it, we need some notation.
    Let $M$ be a matching in a graph $G$. 
    For a vertex $x$, we let $M(x)$ be the unique edge in $M$ incident to $x$ if $x$ is covered by $M$; otherwise, we let $M(x) \defeq \blank$. The set of all vertices not covered by $M$ is denoted by $U_M$.  

    \begin{lemma}\label{lemma:induction}
        For all $\high > \low \in \N$, there are  $p$, $\lambda$, $C > 0$ such that the following holds. Let $G = (V,E,\pi)$ be a locally finite Borel graph with a bipartition $(A,B)$ of type $(\high,\low)$. Then, for every Borel matching $M$ in $G$ and for all $\epsilon > 0$, there exists a random Borel matching $M'$ in $G$ such that for all $x \in A$:
        \begin{enumerate}[label=\ep{\normalfont\Alph*}]
            \item\label{item:grow} If $x$ is not covered by $M$, then $\P[\text{$x$ is covered by $M'$}] \geq p$.
            \item\label{item:not_shrink} If $x$ is covered by $M$, then $\P[\text{$x$ is not covered by $M'$}] \leq \epsilon$.
            \item\label{item:converge} $\P[M'(x) \neq M(x)] \leq C\exp(-\lambda \,\dist_G(x,A \cap U_M))$.
        \end{enumerate}
    \end{lemma}

    Condition~\ref{item:grow} says that we are making substantial progress: each vertex that is not yet covered by the matching becomes so with some guaranteed positive probability. Condition~\ref{item:not_shrink} means that we are not losing the progress that has already been made: since $\epsilon$ can be made arbitrarily small, covered vertices only extremely rarely become uncovered. 
    (With some extra care, it is even possible to allow $\epsilon = 0$, but we prefer to stick with the somewhat easier construction for $\epsilon > 0$.) Finally, condition~\ref{item:converge} constitutes the main innovation of our argument: it provides a way to control the difference between $M'$ and $M$, which will enable us to apply the Borel--Cantelli lemma to pass to a limit after infinitely many iterative applications of the lemma.  

    It should be noted that the measure-preserving case of Theorem~\ref{theo:measurable} is also proved by growing the matching iteratively and invoking the Borel--Cantelli lemma in the end. However, the invariance of measure is used there to control both the amount of progress made on each step \cite[Lem.~2.2]{folner.tilings} and the difference between the consecutive matchings in the resulting sequence \cite[Lem.~2.5]{folner.tilings}. We expect that our approach---using randomized constructions and controlling the probability of a change at a point $x$ as a function of the distance from $x$ to the set of ``problematic'' elements---will find further applications to combinatorial problems in the non-invariant measure setting.

    To facilitate the iterative construction, we will need to apply Lemma~\ref{lemma:induction} in the case when $M$ is itself a \emph{random} Borel matching in $G$. This is possible thanks to the fact that a random Borel matching $M$ in a Borel graph $G = (V,E,\pi)$ can also be viewed as an ordinary (i.e., \emph{non-random}) Borel matching in an auxiliary graph $G^*$. Namely, letting the parameter space for $M$ be $(\Omega, \P)$, we define $G^* = (V^*, E^*, \pi^*)$ by setting $V^* \defeq \Omega \times V$, $E^* \defeq \Omega \times E$, and
    \[
        \pi^*(\omega, e) \,\defeq\, \set{(\omega, x), \, (\omega, y)}, \text{ where } \pi(e) = \set{x,y}.
    \]
    In other words, $G^*$ is formed by taking a family of disjoint copies of $G$ indexed by the points $\omega \in \Omega$. Since $M_\omega$ is a matching in $G$ for all $\omega \in \Omega$, $M$ is indeed a matching in $G^*$, and hence we may apply Lemma~\ref{lemma:induction} to $M$ using $G^*$ in place of $G$. 
    This yields a 
    new random Borel matching $M'$ in $G^*$, say with parameter space $(\Omega', \P')$. But then $M'$ is also 
    a random Borel matching in the original graph $G$ with parameter space $(\Omega' \times \Omega, \P' \times \P)$. 
    In this manner, we can repeatedly invoke Lemma~\ref{lemma:induction} to generate a sequence of random Borel matchings in $G$, with product spaces serving as their parameter spaces. 
    %
    %
    %
    %
    We derive
    Theorem~\ref{theo:random} from Lemma~\ref{lemma:induction} by analyzing this sequence.

    \begin{scproof}[ of Theorem~\ref{theo:random}]
        Suppose $(A,B)$ is of type $(\high,\low)$, where $\high > \low$.
        Using Lemma~\ref{lemma:reg_subgraph}, we may pass to a spanning subgraph of $G$ and assume that $\deg_G(x) = \high$ for all $x \in A$. In particular, the maximum degree of $G$ is at most $\high$. Let $p$, $\lambda$, $C > 0$ be as in Lemma~\ref{lemma:induction}. 
        By a \emphd{random sequence of Borel matchings} we mean a sequence of random Borel matchings with the same parameter space. Starting with $M_0 \defeq \0$, we iteratively apply Lemma~\ref{lemma:induction} and generate a random sequence of Borel matchings 
        $(M_n)_{n \in \N}$ 
        such that for all $n \in \N$ and $x \in A$:
        \begin{enumerate}[label=\ep{\normalfont\Alph*}]
            \item $\P[\text{$x$ is covered by $M_{n+1}$} \,\big\vert\, \text{$x$ is not covered by $M_n$}] \geq p$,
            \item $\P[\text{$x$ is not covered by $M_{n+1}$} \,\big\vert\, \text{$x$ is covered by $M_n$}] \leq \epsilon_n \defeq (p/2)(1-p/2)^n$,
            \item\label{item:key} $\P[M_{n+1}(x) \neq M_n(x) \,\big\vert\, M_n] \leq C\exp(-\lambda \,\dist_G(x,A \cap U_{M_n}))$.
        \end{enumerate}
        (The parameter space for the sequence $(M_n)_{n \in \N}$ is a suitable infinite product space.) Define
            \[
             M \,\defeq\, \liminf_{n \to \infty} M_n \,=\, \set{e \in E \,:\, \text{$e \in M_n$ for all large enough $n$}}.
            \]
        This is clearly a random Borel matching in $G$. We will show that $M$ is as desired, i.e., it covers each individual vertex $x \in A$ almost surely. 

        \begin{claim}\label{claim:prob_bound}
            For all $x \in A$ and $n \in \N$,  $\P[\text{$x$ is not covered by $M_n$}] \leq \left(1 - p/2\right)^n$.
        \end{claim}
        \begin{claimproof}
            For $n = 0$, the claim is clear, and if it holds for some $n$, then 
        \begin{align*}
            \P\left[\text{$x$ is not covered by $M_{n+1}$}\right] \,&\leq\, (1 - p)\, \P\left[\text{$x$ is not covered by $M_n$}\right] \,+\, \epsilon_n \,\P\left[\text{$x$ is covered by $M_n$}\right] \\
            & \leq\, (1-p)\left(1 - p/2\right)^n  + \epsilon_n \,=\, \left(1 - p/2\right)^{n+1}. \qedhere
        \end{align*}
        \end{claimproof}

        \begin{claim}\label{claim:stable}
            For all $x \in A$ and $n \in \N$, $\P[M_{n+1}(x) \neq M_n(x)]$ is exponentially small in $n$.  
        \end{claim}
        \begin{claimproof}
            Given $R \geq 0$, let $B_G(x,R)$ be the $R$-ball around $x$ in $G$. 
            Note that $|B_G(x,R)| \leq 2\high^R$. Hence, for any $n \in \N$ and $R \geq 0$, we can use Claim~\ref{claim:prob_bound} to write
        \begin{align*}
            \P\left[\dist_G(x, A \cap U_{M_n}) \leq R\right] \,&=\, \P\left[B_G(x,R) \cap A \cap U_{M_n} \neq \0\right]\\
            &\leq\, \sum_{y \,\in\, B_G(x,R) \cap A} \P\left[y \in U_{M_n}\right] \,\leq\, 2\high^R\,(1 - p/2)^n.
        \end{align*}
        Fix some $\gamma > 0$ such that $\delta \defeq \high^\gamma(1 - p/2) < 1$. Then
        \begin{equation}\label{eq:close}
            \P\left[\dist_G(x, A \cap U_{M_n}) \leq \gamma n\right] \,\leq\, 2\high^{\gamma n}\, (1 - p/2)^{n} \,=\,2\delta^n. 
        \end{equation}
        Crucially, condition \ref{item:key} implies that
        \begin{equation}\label{eq:far}
            \P\left[M_{n+1}(x) \neq M_n(x) \,\middle\vert\, \dist_G(x, A \cap U_{M_n}) > \gamma n\right] \,\leq\, Ce^{-\lambda\gamma n}. 
        \end{equation}
        Putting \eqref{eq:close} and \eqref{eq:far} together, we obtain the desired  bound
        \[
            \P\left[M_{n+1}(x) \neq M_n(x)\right] \,\leq\, 2\delta^n + Ce^{-\lambda\gamma n}. \qedhere 
        \]
        \end{claimproof}
         Now consider any $x \in A$. By Claims~\ref{claim:prob_bound} and \ref{claim:stable},
         \[
            \sum_{n = 0}^\infty \P\left[\text{$x$ is not covered by $M_n$ or } M_{n+1}(x) \neq M_n(x)\right] \,<\,\infty.
         \]
         It follows by the Borel--Cantelli lemma that, with probability $1$, the sequence $(M_n(x))_{n \in \N}$ is eventually constant and equal to some $e \in E$. Then $M(x) = e$ as well, and the proof is complete.
    \end{scproof}


    We remark that the above proof of Theorem~\ref{theo:random} would work even if condition \ref{item:converge} in Lemma~\ref{lemma:induction} is replaced by $\P[M'(x) \neq M(x)] \leq f(\dist_G(x,A \cap U_M))$ for any summable function $f \colon \N \to [0,\infty)$.

    \subsection{Proof of Lemma~\ref{lemma:induction}}


    We prove Lemma~\ref{lemma:induction} using a standard and classical tool in matching theory, namely \emph{augmenting paths} \cites[\S2.1]{Diestel}[\S3.1]{West}. The key insight is to generate the requisite augmenting paths using a random walk.


    \begin{definition}[Alternating walks]\label{defn:alt}
        Let $G = (V,E,\pi)$ be a graph and let $M$ be a matching in $G$. An \emphd{$M$-alternating walk} of length $n+1$ is a sequence $w = (x_0, e_0, \ldots, x_n, e_n, x_{n+1})$, where:
        \begin{itemize}
            \item $x_0$, \ldots, $x_{n+1} \in V$ and $e_0$, \ldots, $e_n \in E$,
            \item $\pi(e_i) = \set{x_i, x_{i+1}}$ for all $0 \leq i \leq n$, and
            \item $e_i \in M$ for all odd $1 \leq i \leq n$.
        \end{itemize}
        We write ${V}(w) \defeq \set{x_0, \ldots, x_{n+1}}$ and ${E}(w) \defeq \set{e_0, \ldots, e_n}$. Note that we do not require the vertices $x_0$, \ldots, $x_{n+1}$ to be distinct. 
    \end{definition}

    It is well known that alternating walks can be used to shrink the set of uncovered vertices: 

    \begin{prop}[{Augmenting paths}]\label{prop:augment}
        Let $G = (V,E,\pi)$ be a bipartite graph and let $M$ be a matching in $G$. Suppose $w = (x_0, e_0, \ldots, x_n, e_n, x_{n+1})$ is an $M$-alternating walk of odd length such that $x_0$ and $x_{n+1}$ are not covered by $M$. Then there is a set $P \subseteq {E}(w)$ with the following properties:
        \begin{itemize}
            \item $M' \defeq M \symdif P$ is a matching in $G$, and
            \item $U_{M'} = U_M \setminus \set{x_0, x_{n+1}}$.
        \end{itemize}
    \end{prop}
    \begin{scproof}
        Note that $x_0 \neq x_{n+1}$ because $G$ is bipartite and $n$ is even. Upon truncating $w$ if necessary, we may assume that $x_{n+1}$ does not appear among $x_0$, \ldots, $x_n$. Similarly, by passing to a subsequence of $w$, we may arrange that each $x_i$ does not appear among $x_0$, \ldots, $x_{i-1}$, i.e., 
        $x_0$, \ldots, $x_{n+1}$ are all distinct. 
        We may then take 
        $P = {E}(w)$. 
        Note that the 
        set $P$ produced by this construction 
        is the edge set of a path in $G$, called an \emph{$M$-augmenting path} in the literature.
    \end{scproof}

    We also need the following slight variant of Proposition~\ref{prop:augment}:

    \begin{prop}\label{prop:not_augment}
        Let $G = (V,E,\pi)$ be a bipartite graph and let $M$ be a matching in $G$. Suppose $w = (x_0, e_0, \ldots, x_n, e_n, x_{n+1})$ is an $M$-alternating walk of even length such that $x_0$ is not covered by $M$. Then there is a set $P \subseteq {E}(w)$ with the following properties:
        \begin{itemize}
            \item $M' \defeq M \symdif P$ is a matching in $G$, and
            \item $U_{M'} = U_M \symdif \set{x_0, x_{n+1}} = (U_M \setminus \set{x_0}) \cup \set{x_{n+1}}$.
        \end{itemize}
    \end{prop}
    \begin{scproof}
        Again, we may assume $x_0$, \ldots, $x_{n+1}$ are all distinct and take $P = {E}(w)$.
    \end{scproof}

    We now describe a natural way of randomly generating an alternating walk of bounded length starting at a given uncovered vertex.

    \begin{definition}[Random alternating walks]\label{defn:random_walk}
        Let $G = (V,E,\pi)$ be a locally finite graph and let $M$ be a matching in $G$. Fix an even integer $N \geq 2$. Given a vertex $x \in U_M$, we randomly generate an $M$-alternating walk $w = (x_0, e_0, \ldots, x_n, e_n, x_{n+1})$ of length $n+1 \leq N$ as follows. 

        \smallskip

        \begin{algorithm}[H]\DontPrintSemicolon

             Start with $x_0 \defeq x$.\;

            \smallskip
            
            \For{$i =0$, \ldots, $N-1$}{

            \smallskip

            \uIf{\normalfont$i$ is even}{
            
            \smallskip
            
            Pick an edge $e_i$ incident to $x_i$ uniformly at random and let $x_{i+1}$ be its other endpoint.\;

            \smallskip
  }
  \uElseIf{\normalfont$x_i$ is covered by $M$}{
  
            \smallskip
            
            Let $e_i \in M$ be the unique edge incident to $x_i$ and let $x_{i+1}$ be its other endpoint.\;

            \smallskip
  }
  \Else{
  
            \smallskip
    
    Stop the construction. \;

            \smallskip
        }
            }
    \end{algorithm}

    \smallskip

    \noindent We denote the resulting probability distribution on $M$-alternating walks by $\mathsf{Alt}_N(x, M)$.
    \end{definition}



    \begin{lemma}\label{lemma:augment_random}
        Let $G = (V,E,\pi)$ be a locally finite bipartite graph and let $M$ be a matching in~$G$. Suppose $x \in U_M$ and $w = (x_0, e_0, \ldots, x_n, e_n, x_{n+1})$ is sampled from $\mathsf{Alt}_N(x,M)$ for some even $N \geq 2$. Then there is a set $P \subseteq {E}(w)$ with the following properties:
        \begin{itemize}
            \item $M' \defeq M \symdif P$ is a matching in $G$, and
            \item $M'$ covers $x$ as well as every vertex that is covered by $M$, except $x_N$ if it exists.
        \end{itemize}
    \end{lemma}
    \begin{scproof}
        By construction, if $n+1 < N$, then $w$ satisfies the assumptions of Proposition~\ref{prop:augment}, and otherwise it satisfies the assumptions of Proposition~\ref{prop:not_augment}. 
    \end{scproof}

    We now have the following simple but crucial fact:

    \begin{lemma}\label{lemma:decay}
        Let $G = (V,E,\pi)$ be a locally finite bipartite graph with a bipartition $(A,B)$ of type $(\high,\low)$ 
        and let $M$ be a matching in $G$. Fix an even integer $N \geq 2$. Set $U \defeq A \cap U_M$. For each $x \in U$, 
        let $w_x$ be sampled from $\mathsf{Alt}_N(x,M)$. Then, for each $y \in V$ and $i \in \N$,
        \[
            \sum_{x \in U} \P\left[\text{$y$ is the $i$-th vertex in $w_x$}\right] \,\leq\, \left(\frac{\low}{\high}\right)^{i/2}. 
        \]
    \end{lemma}

    Note that the statement of Lemma~\ref{lemma:decay} does not make any assumptions on the joint distribution of the walks $(w_x)_{x \in U}$.

    \begin{scproof}
        For concreteness, let us assume that $y \in A$; the proof in the $y \in B$ case is the same, \emph{mutatis mutandis}. Since $y \in A$, the probability in question is $0$ when $i$ is odd, so we shall assume that $i$ is even. Let $W_i(y)$ be the set of all $M$-alternating walks of length $i$ ending at $y$, and let $W_i(x,y) \subseteq W_i(y)$ be the set of all such walks starting at $x$. Note that 
            \[
                |W_i(y)| \,\leq\, \low^{i/2}.
            \]
        Indeed, every walk $w \in W_i(y)$ has the form $w = (x_0, e_0, \ldots, x_{i-1}, e_{i-1}, x_i)$, where $x_i = y$ and, given $x_j$, there is a unique choice for $e_{j-1}$, $x_{j-1}$ if $j$ is even (as in that case $e_{j-1}$ is the unique edge in $M$ incident to $x_{j}$)  and at most $\low$ such choices if $j$ is odd (since $x_j \in B$ and so $\deg_G(x_j) \leq \low$). Observe also that, for any $x \in U$ and $w \in W_i(x,y)$, we have
        \[
            \P\left[\text{$w$ is a prefix of $w_x$}\right] \,\leq\, \high^{-i/2},
        \]
        because if $w = (x_0, e_0, \ldots, x_{i-1}, e_{i-1}, x_i)$, 
        then for all even $j < i$, on step $j$ of the construction of $w_x$, the edge $e_{j}$ is chosen with probability $1/\deg_G(x_j) \leq \high^{-1}$. Putting these bounds together, we get
        \begin{align*}
            \sum_{x \in U} \P\left[\text{$y$ is the $i$-th vertex in $w_x$}\right] \,&=\, \sum_{x \in U} \sum_{w \in W_i(x,y)} \P\left[\text{$w$ is a prefix of $w_x$}\right]\\
            &\leq\, \sum_{x \in U} |W_i(x,y)|\,\high^{-i/2}\\
            &\leq\, |W_i(y)|\, \high^{-i/2} \,\leq\, \left(\frac{\low}{\high}\right)^{i/2}. \qedhere
        \end{align*}
     \end{scproof}

    After these preliminaries, we are ready to proceed with the proof of Lemma~\ref{lemma:induction}.

    \begin{scproof}[ of Lemma~\ref{lemma:induction}]
        Let $G = (V,E,\pi)$ be a locally finite Borel graph with a bipartition $(A,B)$ of type $(\high,\low)$, where $\high > \low$, and let $M$ be a Borel matching in $G$. Fix $\epsilon > 0$. Our aim is to construct a random Borel matching $M'$ in $G$ satisfying conditions \ref{item:grow}--\ref{item:converge}, for a suitable choice of $p$, $\lambda$, $C > 0$.

        Let $U \defeq A \cap U_M$. Fix a large even integer $N \geq 2$ (it will become clear later, in Claim~\ref{claim:B}, how large $N$ needs to be). 
        Using Lemma~\ref{lemma:reg_subgraph}, we replace $G$ by a spanning subgraph and arrange that $\deg_G(x) = \high$ for all $x \in A$ (this is not strictly necessary for the proof but helps to simplify the presentation). 
        To make navigating the graph easier, we apply Lemma~\ref{lemma:enum} and fix a Borel function
        \[
            \mathsf{edge} \colon \set{1,\ldots, \high} \times A \to E
        \]
        such that an edge $e \in E$ is incident to a vertex $x \in A$ if and only if $e = \mathsf{edge}(i,x)$ for some $1 \leq i \leq \high$. Using the function $\mathsf{edge}$, the randomized procedure in Definition~\ref{defn:random_walk} can be reorganized as follows: First choose a tuple $\bm{c} = (c_0, c_2, \ldots, c_{N-2}) \in \set{1,\ldots, \high}^{N/2}$ uniformly at random, and then, at every even step $i$, let the next edge $e_i$ in the sequence be
        \[
            e_i \,\defeq\, \mathsf{edge}(c_i, x_i).
        \]
        The resulting walk $w$ will clearly be distributed according to $\mathsf{Alt}_N(x,M)$.

        We now wish to create a random assignment $x \mapsto \bm{c}(x) \in \set{1,\ldots,\high}^{N/2}$ and use it to construct random walks $w_x \sim \mathsf{Alt}_N(x,M)$ starting at every vertex $x \in U$. A slight technical issue is that we cannot just choose the tuples $\bm{c}(x)$ for all $x \in U$ independently, because the mapping $x \mapsto \bm{c}(x)$ must be Borel. We get around this problem via a \emph{randomness conservation trick}. That is, we reuse the same value $\bm{c}(x)$ for multiple vertices $x$. Although the resulting assignments are not mutually independent, they can be made ``locally independent,'' as explained below. This idea has already been successfully employed in descriptive combinatorics, see, e.g., \cites[\S5]{BerShift}{subexp1}.

        We use Corollary~\ref{corl:spaced}\ref{item:spaced_deg} to partition the set $U$ into finitely many Borel sets:
        \[
            U \, =  \, U_1 \sqcup \ldots \sqcup U_k,
        \]
        so that the distances in $G$ between the vertices in each $U_i$ strictly exceed $4N$. We now pick
        \[
            \bm{c}_1,\, \ldots,\, \bm{c}_k \,\in\, \set{1,\ldots, \high}^{N/2}
        \]
        independently and uniformly at random and set $\bm{c}(x) \defeq \bm{c}_i$ for all $x \in U_i$. 
        For each choice of the tuples $\bm{c}_1$, \ldots, $\bm{c}_k$, the corresponding assignment $x \mapsto \bm{c}(x)$ is clearly Borel, and we can use it to define a Borel mapping $x \mapsto w_x$ as described above. 
        By construction, we have $w_x \sim \mathsf{Alt}_N(x,M)$ for all $x \in U$. Furthermore, for each $y \in V$, the random walks
        \[
            \big\{w_x \,:\, x \in U \text{ and } \dist_G(x,y) \leq 2N \big\}
        \]
        are mutually independent, because the corresponding vertices $x$ lie in distinct sets $U_i$ and hence get assigned mutually independent tuples $\bm{c}_i$. This is what we meant by ``local independence'' earlier.


        Next we want to make the walks $w_x$ vertex-disjoint. To this end, we say that a vertex $y \in U$ \emphd{threatens} a vertex $x \in U$ if $x \neq y$ and for some $i \leq j$, the $i$-th vertex in $w_x$ coincides with the $j$-th vertex in $w_y$. Note that if $x \neq y$ and ${V}(w_x) \cap {V}(w_y) \neq \0$, then at least one of $x$, $y$ threatens the other. 
        At this point, we make one more random choice. Namely, we pick elements
        \[
            \sigma_1, \, \ldots, \, \sigma_k \,\in\, \set{0,1}
        \]
        independently and uniformly at random and call a vertex $x$ \emphd{active} if $x \in U_i$ with $\sigma_i = 1$. 
        Let $U'$ be the set of all active vertices that are not threatened by any other active vertex. Since the vertices in $U'$ cannot threaten each other, the walks $\set{w_x \,:\, x \in U'}$ are indeed vertex-disjoint. 

        Finally, we apply Lemma~\ref{lemma:augment_random} to pick, for every $x \in U'$, a subset $P_x \subseteq {E}(w_x)$ such that $M \symdif P_x$ is a matching that covers both $x$ and every vertex covered by $M$, except the $N$-th vertex in $w_x$ if it exists. Since there are only finitely many choices for $P_x$ given $w_x$, the map $x \mapsto P_x$ can be made Borel by Luzin--Novikov Uniformization \ep{i.e., Lemma~\ref{lemma:LNU}}. 
        Thus, we can define a Borel set $P \defeq \bigcup_{x \in U'} P_x$. 
        As the walks $\set{w_x \,:\, x \in U'}$ are vertex-disjoint, the set \[M' \,\defeq\, M \symdif P\] is a matching in $G$. 
        This completes our construction. Note that the parameter space for our random Borel matching $M'$ is a finite set with the uniform probability measure, namely
        \[
            \Omega \,\defeq\, \left(\set{1,\ldots, \high}^{N/2} \times \set{0,1}\right)^k.
        \]

        For a vertex $x \in U$, let $T(x)$ be the (random) set of all $y \in U$ that threaten $x$. We start our analysis of $M'$ by bounding the expected size of $T(x)$:

        \begin{claim}\label{claim:threats}
            For all $x \in U$, $\E[|T(x)|] \leq 4\high^2$.
        \end{claim}
        \begin{claimproof}
            Let $D(x)$ be the set of all vertices $y \in U$ such that $0 < \dist_G(x,y) \leq 2N$ and note that we always have $T(x) \subseteq D(x)$. Consider an arbitrary $M$-alternating walk $w = (x_0, e_0, \ldots, x_n, e_n, x_{n+1})$ of length $n+1 \leq N$ starting at $x_0 = x$. By the linearity of expectation,
            \begin{align*}
                \E\left[|T(x)| \,\middle\vert\, w_x = w\right] \,&\leq\, \sum_{i = 0}^{n+1} \sum_{j = i}^N \sum_{y \in D(x)} \P\left[\text{$x_i$ is the $j$-th vertex in $w_y$} \,\vert\, w_x = w\right] \\
                &=\, \sum_{i = 0}^{n+1} \sum_{j = i}^N \sum_{y \in D(x)} \P\left[\text{$x_i$ is the $j$-th vertex in $w_y$}\right],
            \end{align*}
            where we use that $w_x$ and $w_y$ are probabilistically independent when $0 < \dist_G(x,y) \leq 2N$. By Lemma~\ref{lemma:decay}, for each $0 \leq i \leq n+1$ and $j \in \N$, we have
            \begin{align*}
                \sum_{y \in D(x)} \P\left[\text{$x_i$ is the $j$-th vertex in $w_y$}\right] \,\leq\, \left(\frac{\low}{\high}\right)^{j/2}. 
            \end{align*}
            Therefore, using that $\low \leq \high-1$ and hence $(\low/\high)^{1/2} \leq 1 - 1/(2\high)$, we obtain
            \begin{align*}
                \E\left[|T(x)| \,\middle\vert\, w_x = w\right] \,&\leq\, \sum_{i=0}^{n+1} \sum_{j = i}^N \left(\frac{\low}{\high}\right)^{j/2} \,\leq\, \sum_{i=0}^{\infty} \sum_{j = i}^\infty \left(1 - \frac{1}{2\high}\right)^{j} \,=\, 
                4\high^2.
            \end{align*}
            Since this holds for every $w$, the desired unconditional bound also follows.
        \end{claimproof}
        
        Now we can check that $M'$ has the properties stated in the lemma.

        \begin{claim}
            $M'$ satisfies condition \ref{item:grow} with $p \defeq 2^{-1-4\high^2}$. That is, for all $x \in U$,
            \[
                \P\left[\text{$x$ is covered by $M'$}\right] \,\geq\, p.
            \]
        \end{claim}
        \begin{claimproof}
            By construction, $x \in U$ is covered by $M'$ if and only if $x \in U'$, which happens if $x$ is active and no $y \in T(x)$ is active. Since every vertex becomes active with probability $1/2$ and the random choices for the vertices at distance at most $2N$ from $x$ are mutually independent, we can write
            \[
                \P\left[x \in U' \,\middle\vert\, T(x) \right] \,=\, \frac{1}{2^{1 + |T(x)|}}.
            \]
            Using the convexity of the function $t \mapsto 1/2^t$ and Claim~\ref{claim:threats}, we obtain the desired bound
            \[
                \P\left[x \in U'\right] \,=\, \E\left[\frac{1}{2^{1 + |T(x)|}}\right] \,\geq\, \frac{1}{2^{1 + \E\left[|T(x)|\right]}} \,\geq\, \frac{1}{2^{1 + 4\high^2}}. \qedhere
            \]
        \end{claimproof}

        \begin{claim}\label{claim:B}
            Assuming $N$ is sufficiently large, $M'$ satisfies condition \ref{item:not_shrink}, i.e., for all $x \in A \setminus U$,
            \[
                \P\left[\text{$x$ is not covered by $M'$}\right] \,\leq\, \epsilon.
            \]
        \end{claim}
        \begin{claimproof}
            By construction, $x \in A \setminus U$ is not covered by $M'$ if and only if $x$ is the $N$-th vertex of the walk $w_y$ for some $y \in U'$. By Lemma~\ref{lemma:decay}, the probability of this event is at most
            \[
                \sum_{y \in U}\P\left[\text{$x$ is the $N$-th vertex in $w_y$}\right] \,\leq\, \left(\frac{\low}{\high}\right)^{N/2},
            \]
            which can be made arbitrarily small by choosing $N$ large enough.
        \end{claimproof}

        \begin{claim}
            $M'$ satisfies condition \ref{item:converge} with $C \defeq 2\high$ and $\lambda \defeq 1/(2\high)$. That is, for all $x \in A$,
            \[
                \P\left[M'(x) \neq M(x)\right] \,\leq\, C \exp(-\lambda\, \dist_G(x,U)).
            \]
        \end{claim}
        \begin{claimproof}
            Let $R \defeq \dist_G(x,U)$. If $M'(x) \neq M(x)$, then, in particular, $x \in V(w_y)$ for some $y \in U$. Note that $x$ cannot be the $i$-th vertex in $w_y$ for $i < R$. Hence, Lemma~\ref{lemma:decay} yields
            \begin{align*}
                \P\left[M'(x) \neq M(x)\right] \,&\leq\, \sum_{i = R}^N \sum_{y \in U} \P\left[\text{$x$ is the $i$-th vertex in $w_y$}\right] \\
                &\leq\, \sum_{i = R}^\infty \left(\frac{\low}{\high}\right)^{i/2} \,\leq\, 2\high\,\left(1 - \frac{1}{2\high}\right)^R \,\leq\, 2\high \, e^{-R/(2\high)},
            \end{align*}
            as desired. Here we again use that $(\low/\high)^{1/2} \leq 1 - 1/(2\high)$, together with the bound $1- t \leq e^{-t}$. 
        \end{claimproof}

        We have shown that the random Borel matching $M'$ has all the desired properties. This completes the proof of Lemma~\ref{lemma:induction}, and hence also of Theorems~\ref{theo:random} and \ref{theo:measurable}.
    \end{scproof}

    \section{Graphs of subexponential growth and everywhere 2-ended graphs}\label{sec:Borel_simple}

    We now turn to Borel constructions that avoid relying on measure theory or Baire category. In \S\ref{sec:Borel_simple} we consider two classes of graphs: subexponential growth and everywhere $2$-ended. Our results for them are proved with elementary and relatively standard techniques. In \S\ref{sec:asi} we will move on to graphs of finite asymptotic separation index, which require significantly more advanced machinery.

    \subsection{Graphs of subexponential growth: Theorem~\ref{theo:growth}}

    In the statement of Theorem~\ref{theo:growth}, we are given a Borel graph $G$ of subexponential growth and a bipartition $(A,B)$ that has combinatorial expansion.\footnote{It is not hard to show that under these assumptions, the bipartition $(A,B)$ must be Borel. However, the proof of Theorem~\ref{theo:growth} does not make explicit use of this fact.} In other words, there is a constant $c > 1$ such that for every finite set $S \subseteq A$, we have $|N_G(S)| \geq c|S|$. (We again note that no assumption is made about subsets of $B$, so there is no contradiction with subexponential growth.)
    The goal is to find a Borel matching in $G$ covering every vertex in $A$.
    The proof of this result, like that of Theorem~\ref{theo:measurable}, relies on alternating walks and augmenting paths. See Definition~\ref{defn:alt} and Proposition~\ref{prop:augment} for the relevant concepts.
    
    A common method for solving combinatorial problems on Borel graphs of subexponential growth is to argue that any partial but incomplete solution can always be improved via a ``local'' modification; see, e.g., \cites{bernshteyn2023borel}[\S8]{GV_from}. In line with this philosophy, we show that in the context of Theorem~\ref{theo:growth}, a matching in $G$ that does not cover all of $A$ admits an augmenting path of bounded length.

    \begin{lemma}\label{lem:expand}
        Let $G$ be a graph of subexponential growth and let $(A,B)$ be a bipartition of $G$ with combinatorial expansion. Then there exists an integer $L \in \N$ such that, for every matching $M$ in $G$ and every uncovered vertex $x \in A$, 
        there is an $M$-alternating walk of length at most $L$ starting at $x$ and 
        ending at an uncovered vertex in $B$.
    \end{lemma}

    Note that Lemma~\ref{lem:expand} is a purely graph-theoretic statement that makes no reference to the Borel structure of $G$.
    
    \begin{scproof}
        Fix $c > 1$ such that for every finite set $S \subseteq A$, we have $|N_G(S)| \geq c\,|S|$ \ep{such $c$ exists by combinatorial expansion}. Since $G$ is of subexponential growth, there is $L \in \N$ such that for every vertex $x$, we have $|B_G(x,2L)| < c^{L}$. We claim that this choice of $L$ works.

        Let $M$ be an arbitrary matching in $G$ and suppose $x \in A$ is an uncovered vertex. For $S \subseteq B$, let $N_M(S)$ be the set of all vertices $y \in A$ that are joined to some vertex in $S$ by an edge in $M$. Define
        \[
            A_0 \,\defeq\, \set{x}, \qquad B_n \,\defeq\, N_G(A_n), \qquad A_{n+1} \,\defeq\, N_M(B_n).
        \]
        Then $B_n$ is precisely the set of all vertices that can be reached from $x$ by an $M$-alternating walk of length $n+1$. Therefore, we are done if for some $n < L$, $B_n$ includes an uncovered vertex. Otherwise, we have $|A_{n+1}| = |B_n|$ for all $n < L$. We also have $|B_n| \geq c\,|A_n|$ for all $n$ by combinatorial expansion. It follows that $|A_L| \geq c^L$, which is impossible as $A_L \subseteq B_G(x, 2L)$.
    \end{scproof}

    On the other hand, as observed by Elek and Lippner \cite{elek2010borel}, it is always possible to find a Borel matching without short augmenting paths:

    \begin{prop}[{Elek--Lippner \cite[Prop.~1.1]{elek2010borel}}]\label{prop:no_short_improvement}
        Let $G$ be a locally finite Borel graph. Then, for any $L \in \N$, there exists a Borel matching $M$ in $G$ such that no $M$-alternating walk of length at most $L$ starts  and ends at two distinct uncovered vertices. 
    \end{prop}

    Proposition~\ref{prop:no_short_improvement} is stated in \cite{elek2010borel} for graphs $G$ of finite maximum degree. However, exactly the same proof works in the more general locally finite case as well. \ep{That being said, we will only use it for graphs of subexponential growth, whose maximum degree is by definition finite.}

    \begin{scproof}[ of Theorem~\ref{theo:growth}]
        Let $G$ be a Borel graph of subexponential growth and let $(A,B)$ be a bipartition of $G$ with combinatorial expansion. Fix an integer $L \in \N$ as in Lemma~\ref{lem:expand} and let $M$ be a Borel matching in $G$ given by Proposition~\ref{prop:no_short_improvement}. Then $M$ must cover $A$, since if $x \in A$ is uncovered, then there is an $M$-alternating walk of length at most $L$ starting at $x$ and ending at an uncovered vertex in $B$, contradicting the choice of $M$.
    \end{scproof}

    \subsection{Everywhere 2-ended graphs: Theorem~\ref{theo:ends}}

    Here we consider everywhere $2$-ended graphs. The theory of such graphs in the Borel setting was developed by Miller in \cite{miller.2end}, and our approach draws on the techniques in that paper as well as on their refinements in \cite{2ends}. The particular way in which combinatorial expansion comes into play in our proof is inspired by the arguments in \cite{marks2016baire} due to Marks and Unger.

    Recall the notation $G[X]$ and $G[X, Y]$ for induced subgraphs and induced bipartite subgraphs from \S\ref{subsec:graphs}. A \emphd{separator} in a connected graph $G$ is a finite set $F$ of vertices such that:
    \begin{itemize}
        \item the graph $G[F]$ is connected, and
        \item $G - F$ has at least two infinite components.
    \end{itemize}
    We will need the following simple combinatorial fact regarding separators in $2$-ended graphs:

    \begin{lemma}\label{lemma:separators}
       Let $G$ be a $2$-ended locally finite graph. Suppose $\mathcal{F}$ is a family of pairwise disjoint separators in $G$ and let $X \defeq \bigcup \mathcal{F}$. Then every connected component of the graph $G - X$ is adjacent in $G$ to at most two sets in $\mathcal{F}$. 
    \end{lemma}
    \begin{scproof}
        Suppose 
     there are three distinct sets $F_1$, $F_2$, $F_3 \in \mathcal{F}$ adjacent to the same component $C$ of $G - X$. Pick any vertex $x \in C$ and let $Y_i$ be an infinite component of $G - F_i$ not containing $x$. We claim that $Y_i \cap Y_j = \0$ for $i \neq j$. This implies that the graph $G - F_1 - F_2 - F_3$ has at least $3$ infinite connected components, contradicting the $2$-endedness of $G$. To prove the claim, assume that $y \in Y_i \cap Y_j$ and let $P$ be a shortest path in $G$ from $y$ to $F_i \cup F_j$. Without loss of generality, say $P$ ends at a vertex in $F_i$ and thus avoids $F_j$. Since $F_i$ is adjacent to $C$, there is also a path $P'$ from $x$ to $F_i$ avoiding $F_j$. We can then join $P'$ and $P$ inside $G[F_i]$ to form a path from $x$ to $y$ avoiding $F_j$. Therefore, $x$ and $y$ are in the same component of $G - F_j$, which contradicts $y$ being in $Y_j$.
    \end{scproof}

    In applications \ep{for example, in Theorems~\ref{theo:Shannon_meas} and~\ref{theo:Shannon_Borel}}, it is often necessary to work not just with a given graph itself, but also with its subgraphs. A slight technical issue is that subgraphs of everywhere $2$-ended graphs are not necessarily everywhere $2$-ended themselves. To deal with this, we prove the following minor extension of Theorem~\ref{theo:ends} that applies to subgraphs of everywhere $2$-ended graphs (this generalization has no effect on the proof):

    \begin{theorem}[{Theorem~\ref{theo:ends} for subgraphs}]\label{theo:ends1}
        Suppose that a Borel bipartition $(A,B)$ of a locally finite Borel graph $G$ has combinatorial expansion. If $G$ is a subgraph of some everywhere $2$-ended locally finite Borel graph $H$, then $G$ has a Borel matching covering $A$.
    \end{theorem}

    Throughout the proof, we shall rely on the facts about smooth graphs discussed in \S\ref{subsec:smooth}. 
    
    \begin{scproof} 
        We begin with some general preparations. We may replace $G$ and $H$ by their underlying simple graphs and assume that $G$ and $H$ are simple. Write $H = (V,E)$. Replacing $H$ by $H^2$, we may arrange that $G^2 \subseteq H$. 
        Let $V_k \subseteq V$ be the union of all components of $H$ in which the minimum size of a separator is $k$. Then $V = V_1 \sqcup V_2 \sqcup \ldots $ is a partition of $V$ into Borel $H$-invariant sets, and we may work in each $V_k$ separately. Thus, from now on, we 
        fix $k \in \N^+$ and assume that $V = V_k$.

        Fix $c > 1$ such that for every finite set $S \subseteq A$, we have $|N_G(S)| \geq c\,|S|$. Let $\Phi \subseteq [V]^k$ be the family of all $k$-element sets $F \subseteq V$ such that $F$ is a separator in a connected component of $H$. By the Kechris--Solecki--Todorcevic Theorem~\ref{theo:KST}\ref{item:maximal}, there is a Borel subset $\mathcal{F} \subseteq \Phi$ such that:
        \begin{itemize}
            \item for distinct $F$, $F' \in \mathcal{F}$, we have $\dist_H(F, F') > R \defeq 1 + 4k/(c-1)$, and
            \item the set $\mathcal{F}$ is maximal with this property.
        \end{itemize}
        Let $X \defeq \bigcup \mathcal{F}$ and define 
        \[
            A_0 \,\defeq\, A \cap X, \qquad A' \,\defeq\, \set{x \in A \,:\, \dist_H(x, X) \leq 2k/(c-1)}, \qquad A_1 \,\defeq\, A \setminus A_0.
        \]
        A key observation is that the subgraphs $G[A', B]$ and $G[A_1,B]$ are smooth.
        

        \begin{claim}\label{claim:smooth0}
            The graph $G[A', B]$ is component-finite, hence smooth.
        \end{claim}
        \begin{claimproof}
            Since $G^2 \subseteq H$, it is enough to show that $H[A']$ is component-finite. To this end, suppose $x$, $y \in A'$ are adjacent in $H$ and let $F_x$, $F_y \in \mathcal{F}$ be such that $\dist_H(x,F_x)$, $\dist_H(y,F_y) \leq 2k/(c-1)$. We claim that $F_x = F_y$. Indeed, if $F_x \neq F_y$, then $\dist_H(F_x,F_y) > R$, and thus
            \[
                \dist_H(x,y) \,>\, R - 4k/(c-1) \,=\, 1,
            \]
            which is a contradiction. Hence, if $C$ is a component of $H[A']$, then there is a set $F \in \mathcal{F}$ such that every vertex in $C$ is at distance at most $2k/(c-1)$ from $F$ in $H$. Consequently, $C$ must be finite.
        \end{claimproof}

        \begin{claim}\label{claim:smooth1}
            The graph $G[A_1,B]$ is smooth.
        \end{claim}
        \begin{claimproof}
            It is enough to show that $H - X$ is smooth, as then the graphs $G^2[A_1] \subseteq H[A_1] \subseteq H - X$ would be smooth as well, and we would be done by Lemma~\ref{lemma:smooth_bipartite}. Consider an arbitrary point $x \in V \setminus X$ and let $C$ be the component of $H - X$ containing $x$. Since every component of $H$ contains a separator of size $k$, the maximality of $\mathcal{F}$ implies that $C$ is adjacent in $H$ to at least one set in $\mathcal{F}$. On the other hand, by Lemma~\ref{lemma:separators}, $C$ is adjacent to at most two such sets. Therefore, letting
            \[
                f(x) \,\defeq\, 
                \set{y \in C \,:\, \text{$y$ has a neighbor in $X$}},
            \]
            we see that $f \colon V \setminus X \to \finset{V}$ is a Borel function such that $x$, $y \in V \setminus X$ are in the same component of $H - X$ if and only if $f(x) = f(y)$. Therefore, $H - X$ is smooth.
        \end{claimproof}

        By Hall's Theorem~\ref{theo:Hall}, $G$ has a (not necessarily Borel) matching that covers $A$. It follows by Claim~\ref{claim:smooth0} and Theorem~\ref{theo:smooth} that $G[A',B]$ has a Borel matching $M'$ covering $A'$. Let $M_0 \subseteq M'$ be the  set of all edges in $M'$ incident to a vertex in $A_0$. Note that $M_0$ is a Borel matching that covers $A_0$. Let $B_0 \subseteq B$ be the set of all vertices in $B$ incident to an edge in $M_0$ and define $B_1 \defeq B \setminus B_0$. 

        \begin{claim}\label{claim:Hall}
            The bipartition $(A_1,B_1)$ of $G[A_1,B_1]$ satisfies Hall's condition.
        \end{claim}
        \begin{claimproof}
            Let $S \subseteq A_1$ be a finite set. We need to show that $|N_G(S) \cap B_1| \geq |S|$. Note that
            \[
                |N_G(S) \cap B_0| \,\leq\, |N_H(S) \cap A_0| \,\leq\, |N_H(S) \cap X|,
            \]
            because $G^2 \subseteq H$ and $M_0$ is a matching.  Therefore, we can  use combinatorial expansion to write
            \[
                |N_G(S) \cap B_1| \,=\, |N_G(S)| \,-\, |N_G(S) \cap B_0| \,\geq\, c\,|S| \,-\, |N_H(S) \cap X|.
            \]
            We are done if $N_H(S) \cap X = \0$, so let us assume that $N_H(S) \cap X \neq \0$. 
            We may also assume that the graph $G^2[S]$ is connected, since a minimal counterexample to Hall's condition would have this property. Now we consider two cases depending on the size of $S$.

            \smallskip
            
            \underline{\emph{Case} $1$:} $|S| \leq 2k/(c-1)$. Since $H[S] \supseteq G^2[S]$ is connected and $S$ has a neighbor in $X$, it follows that every vertex in $S$ has distance at most $2k/(c-1)$ to $X$ in $H$, and thus $S \subseteq A'$. In particular, the matching $M'$ covers $S$, and, letting $S'$ be the set of all vertices joined to $S$ by an edge in $M'$, 
            \[
                |N_G(S) \cap B_1| \,\geq\, |S'| \,=\, |S|.
            \]

            \smallskip

            \underline{\emph{Case} $2$:} $|S| \geq 2k/(c-1)$. Note that $S$ is contained in a single connected component $C$ of $H - X$. By Lemma~\ref{lemma:separators}, there exist at most two sets in $\mathcal{F}$ that are adjacent to $C$ in $H$.  Therefore,
            \[
                |N_H(S) \cap X| \,\leq\, |N_H(C) \cap X| \,\leq\, 2k.
            \]
            Using the assumption $|S|  \geq 2k/(c-1)$, we obtain the desired bound
            \[
                |N_G(S) \cap B_1| \,\geq\, c\,|S| \,-\, 2k \,\geq\, |S|. \qedhere
            \]
        \end{claimproof}

        By Claim~\ref{claim:Hall} and Hall's Theorem~\ref{theo:Hall}, the graph $G[A_1,B_1]$ contains a matching $M_1$ covering $A_1$. Since this graph is smooth by Claim~\ref{claim:smooth1}, $M_1$ can be chosen to be Borel. Then $M_0 \cup M_1$ is a Borel matching in $G$ that covers $A$, and the proof is complete.
    \end{scproof}

    \section{Matchings in graphs of finite asymptotic separation index}\label{sec:asi}

    \subsection{Asymptotic separation index: definition and examples}\label{subsec:asi_defns}

    As mentioned in the introduction, asymptotic separation index was introduced by Conley, Jackson, Marks, Seward, and Tucker-Drob in \cite{conley2020borel}. The origins of this notion are twofold. First, a version of it appears implicitly in the paper \cite{conley.miller.toast} by Conley and Miller, although it is not isolated there as a general combinatorial parameter. A second major source of inspiration is large-scale geometry. \emph{Asymptotic dimension} of metric spaces, defined by Gromov \cite[\S1.E]{gromov1993asymptotic}, is particularly relevant here, as the principal aim of \cite{conley2020borel}~was to develop a theory of asymptotic dimension for Borel graphs. Asymptotic separation index is a relaxation of Borel asymptotic dimension that is particularly useful in combinatorial constructions. It is also worth mentioning that some closely related ideas appeared independently in computer science under the names \emph{network decomposition} and \emph{shattering}; see, e.g., \cite{decompose1,decompose2,decompose4,decompose5,decompose6,decompose7,decompose8}.
    

     In the following definition, we say that a set $X$ of vertices in a graph $G$
        is \emphd{$G$-finite} if the induced subgraph $G[X]$ is component-finite. Recall that $G^k$ denotes the $k$-th power of $G$ (see \S\ref{subsec:graphs}).

    \begin{definition}[Asymptotic separation index]\label{defn:asi}
        Let $G$ 
        be a simple locally finite Borel graph.  
             The \emphd{asymptotic separation index} of $G$, in symbols $\asi(G)$, is the smallest $s \in \N$ such that for every $k \in \N$, the vertex set of $G$ can be covered by $s+1$ Borel $G^k$-finite sets. 
             If no such $s \in \N$ exists, we set $\asi(G) \defeq \infty$. 
        If $G$ is a multigraph, we let $\asi(G)$ 
        be equal to that of its 
        underlying simple graph.
    \end{definition}

    It is clear from the definition that $\asi(G) \leq \asi(H)$ if $G$ is a subgraph of $H$, and $\asi(G^k) = \asi(G)$ for simple graphs $G$ and all $k \in \N^+$. 
    The following facts are already implicit in \cite{conley.miller.toast}:

    \begin{theorem}[{\cites[Thm.~4.8(b)]{conley2020borel}[Cor.~2.1.14]{WeilacherThesis}}]\label{theo:asi_ae}
        Let $G = (V,E,\pi)$ be a locally finite Borel graph and let $\mu$ and $\tau$ be a probability measure and a compatible topology on $V$, respectively.

        \begin{itemize}
            \item There exists a $\tau$-meager Borel $G$-invariant set $X \subseteq V$ such that $\asi(G - X) \leq 1$.

            \item If $G$ is {\upshape\emphd{hyperfinite}}, i.e., an increasing union of Borel component-finite subgraphs, then there exists a $\mu$-null Borel $G$-invariant set $X \subseteq V$ such that $\asi(G-X) \leq 1$.
        \end{itemize}
    \end{theorem}

    Everywhere $2$-ended graphs have asymptotic separation index at most $1$. This result is essentially due to Miller \cite{miller.2end}; see \cite[\S2]{2ends} and \cite[Thm.~1.5]{bowen2023definable} for an explicit proof.

    \begin{theorem}[{\cites{miller.2end}[\S2]{2ends}[Thm.~1.5]{bowen2023definable}}]
        Everywhere $2$-ended locally finite Borel graphs have asymptotic separation index at most $1$.  
    \end{theorem}

    Perhaps the most interesting examples of graphs with finite asymptotic separation index arise from Borel actions of certain finitely generated groups.

    \begin{definition}[Schreier graphs]\label{defn:schreier}
        Let $\G$ be a group generated by a finite set $S \subseteq \G$. Suppose $\G \acts X$ is a Borel action of $\G$ on a standard Borel space $X$. The corresponding \emphd{Schreier graph} $\mathsf{Sch}(X, \G,S)$ is the simple graph with vertex set $X$ and edge set
        \[
            \big\{\set{x, \sigma \cdot x} \,:\, x \in X, \, \sigma \in S, \, x \neq \sigma \cdot x\big\}.
        \] 
        The \emphd{asymptotic separation index} of the action $\G \acts X$ is defined as $\asi(\G \acts X) \defeq \asi(\mathsf{Sch}(X, \G, S))$ 
        for any (equivalently, every) finite generating set $S$.
    \end{definition}

    \begin{theorem}[{Conley--Jackson--Marks--Seward--Tucker-Drob \cite{conley2020borel}}]\label{theo:asi1_group}
        Let $\G \acts X$ be a free Borel action of a finitely generated group $\G$ on a standard Borel space $X$. Then $\asi(\G \acts X) \leq 1$, provided that $\G$ has at least one of the following properties: 
        
        \medskip
        
        \noindent
        \begin{minipage}[c]{0.5\linewidth}
        \begin{itemize}
            \item $\G$ is virtually nilpotent,
            \item $\G$ is polycyclic, 
            \item $\G$ is solvable and linear over $\mathbb{Q}$,
        \end{itemize}
        \end{minipage}\hspace*{-30pt} 
        \begin{minipage}[c]{0.6\linewidth}
        \begin{itemize}
            \item $\G = BS(1,2)$, the Baumslag--Solitar group,
            \item $\G = \Z_2 \wr \Z$, the lamplighter group.
        \end{itemize}
        \end{minipage}
    \end{theorem}

    See \cite{conley2020borel} for further examples. Note that most classes of groups in Theorem~\ref{theo:asi1_group} (with the exception of virtually nilpotent ones) can have exponential growth. It is not known whether subexponential growth implies finite asymptotic separation index. However, we have the following:

    \begin{theorem}[{AB--Yu \cite{Grids1}}]\label{theo:asi_poly}
        Borel graphs of polynomial growth have asymptotic separation index at most $1$. 
    \end{theorem}

    If a free probability measure-preserving action of $\G$ has finite asymptotic separation index, then $\G$ is amenable. However, there do exist natural Borel actions of non-amenable groups whose asymptotic separation index is finite. Here are some examples in addition to the ones provided by Theorem~\ref{theo:asi_ae} \ep{see the cited papers for the definitions of the terms used}:

    \begin{theorem}[{Naryshkin--Vaccaro \cite{NarHyperbolic}}]\label{theo:asi_hyperbolic}
        The canonical action of a finitely generated hyperbolic group on its Gromov boundary has asymptotic separation index at most $1$.
    \end{theorem}

    \begin{theorem}[{Iyer--Shinko \cite{IyerShinko}}]\label{theo:asi_generic}
        If $\G$ is a finitely generated group of finite asymptotic dimension, then a generic continuous action of $\G$ on $2^\N$ has asymptotic separation index at most $1$.
    \end{theorem}

    The papers cited above for Theorems~\ref{theo:asi1_group}, \ref{theo:asi_poly}, \ref{theo:asi_hyperbolic}, and \ref{theo:asi_generic} in fact show that the relevant graphs have finite Borel asymptotic dimension (see \cite{conley2020borel}), which implies the upper bound of $1$ on their asymptotic separation index by \cite[Thm.~4.8(a)]{conley2020borel}.  
        %
    We remark that it is not known whether there exist \emph{any} locally finite Borel graphs $G$ with $1 < \asi(G) < \infty$ \cite[15]{conley2020borel}. 

    \subsection{From combinatorial expansion to unbalanced bipartitions}\label{subsec:reduction} 

    In this section we perform the reduction from bipartitions with combinatorial expansion to unbalanced bipartitions in graphs with finite $\asi$. This will allow us to derive Theorems~\ref{theo:asi_two_conjectures} and \ref{theo:hyperfinite}, as well as the part of Theorem~\ref{theo:asi} concerning graphs with asymptotic separation index $1$.
    
    We rely on the machinery developed by the second and third named authors in \cite{bowen2023definable}, which involves combinatorial structures that we call \emph{shattering systems}:

    \begin{definition}[Shattering systems]\label{defn:shattering}
        Let $G = (V,E,\pi)$ be a locally finite Borel graph with the underlying simple graph $\underline{G}$. Given $s$, $k \in \N$, an \emphd{$(s,k)$-shattering system} in $G$ is a tuple $\mathcal{S} = (\mathcal{S}[i])_{i=0}^{s-1}$ of Borel subsets of $V$ such that:
        \begin{itemize}
            \item for each $0 \leq i < s$, the set $\mathcal{S}[i]$, called the \emphd{$i$-th layer} of $\mathcal{S}$, is $\underline{G}^k$-finite,
            \item the \emphd{leftover set} $\mathcal{S}^* \defeq V \setminus \bigcup_{j=0}^{s-1} \mathcal{S}[i]$ is also $\underline{G}^k$-finite.
        \end{itemize}
    \end{definition}

    Comparing Definitions~\ref{defn:asi} and \ref{defn:shattering} makes it clear that a locally finite Borel graph $G$ satisfies $\asi(G) \leq s$ if and only if it admits an $(s,k)$-shattering system for every $k \in \N$. It turns out that, moreover, if $\asi(G) \leq s$, then it is possible to find an arbitrarily large family of $(s,k)$-shattering systems in $G$ that are in some sense ``well-spaced.''

    \begin{definition}[Distance between shattering systems]
        Let $G$ be a locally finite Borel graph. If $\mathcal{S}_0$ and $\mathcal{S}_1$ are $(s,k)$-shattering systems in $G$, then we let \[\dist_G(\mathcal{S}_0, \mathcal{S}_1) \,\defeq\, \min_{0\leq i < s} \dist_G(\mathcal{S}_0[i], \mathcal{S}_1[i]).\] A family $\mathcal{S}_1$, \ldots, $\mathcal{S}_n$ of $(s,k)$-shattering systems is \emphd{$d$-spaced} if $\dist_G(\mathcal{S}_i, \mathcal{S}_j) \geq d$ for all $1 \leq i < j \leq n$.
    \end{definition}

    \begin{lemma}[{MB--FW \cite[Lem.~2.1]{bowen2023definable}}]\label{lemma:shatter}
        Let $G$ be a locally finite Borel graph with $\asi(G) = s < \infty$. Then, for any $k$, $d$, $n \in \N$, $G$ admits a $d$-spaced family $\mathcal{S}_1$, \ldots, $\mathcal{S}_n$ of $(s,k)$-shattering systems. 
    \end{lemma}

    \begin{remark}
        Technically, Lemma~\ref{lemma:shatter} is only proved in \cite{bowen2023definable} for the special case of $k =1$ and $d = 5$. However, the general statement follows immediately by applying \cite[Lem.~2.1]{bowen2023definable} to a suitably high power of the underlying simple graph of $G$.
    \end{remark}

    We shall use the following consequence of Lemma~\ref{lemma:shatter}:

    \begin{lemma}\label{lemma:high_cover}
        Let $G$ be a locally finite Borel graph such that 
        $\asi(G) = s < \infty$. Suppose $(A,B)$ is a Borel bipartition of $G$ that satisfies Hall's condition. Then, for any $n \in \N^+$, there exist Borel matchings $M_1$, \ldots, $M_n$ in $G$ such that every vertex in $A$ is covered by at least $n- s$ of them.
    \end{lemma}
    \begin{scproof}
        By Hall's Theorem~\ref{theo:Hall}, $G$ has a (not necessarily Borel) matching that covers $A$. Let $\mathcal{S}_1$, \ldots, $\mathcal{S}_n$ be a $1$-spaced family of $(s,3)$-shattering systems in $G$. 
        For each $1 \leq j \leq n$, the set \[U_j \,\defeq\, \set{x \in V \,:\, \dist_G(x, \mathcal{S}_j^*) \leq 1}\] is $G$-finite (this follows from the fact that the leftover set $\mathcal{S}_j^*$ itself is $G^3$-finite). Therefore, by applying Theorem~\ref{theo:smooth} to the graph $G[U_j]$, we find a Borel matching $M_j$ in $G$ that covers every vertex in $A \cap \mathcal{S}_j^*$. If a vertex $x \in A$ is not covered by some $M_j$, then $x$ belongs to one of the layers of the shattering system $\mathcal{S}_j$. If $x$ is in both $\mathcal{S}_j[i]$ and $\mathcal{S}_{j'}[i']$ for some $j \neq j'$, then $i \neq i'$, because $\dist_G(\mathcal{S}_{j}, \mathcal{S}_{j'}) \geq 1$. Hence, $x$ can belong to a layer in at most $s$ of the shattering systems $\mathcal{S}_1$, \ldots, $\mathcal{S}_n$, and thus it must be covered by at least $n - s$ of the matchings $M_1$, \ldots, $M_n$, as desired.
    \end{scproof}

    Next we apply Lemma~\ref{lemma:high_cover} to formulate a general result that converts combinatorial expansion into unbalancedness under a finite asymptotic separation index assumption:

    \begin{lemma}[Making a graph unbalanced]\label{lemma:make_unbalanced}
        Let $G = (V,E)$ be a simple locally finite Borel graph such that 
        $\asi(G) = s < \infty$. Let $c > 1$ and suppose $(A,B)$ is a Borel bipartition of $G$ such that every finite set $S \subseteq A$ satisfies $|N_G(S)| \geq c\, |S|$. Then, for any $c' < c$, there exists a locally finite Borel multigraph $G'$ with the following properties:
        \begin{enumerate}[label=\ep{\normalfont\arabic*}]
            \item the vertex set of $G'$ is $V$,
            
            \item the underlying simple graph of $G'$ is a subgraph of $G$, and 
            
            \item\label{item:3_unbalanced} the bipartition $(A,B)$ in $G'$ is of type $(\high, \low)$ for some $\high$, $\low$ such that $\high > c' \,\low$.
        \end{enumerate}
    \end{lemma}
    \begin{scproof}
        Pick $p$, $q \in \N^+$ so that $c' < p/q \leq c$. Define a simple graph $G^*$ as follows. Let
        \[
            A^* \,\defeq\, A \times \set{1, \ldots, p}, \qquad B^* \,\defeq\, B \times \set{1, \ldots, q}, \qquad \text{and} \qquad V^* \,\defeq\, A^* \sqcup B^*.
        \]
        The vertex set of $G^*$ is $V^*$, and its edge set is
        \[
            E^* \,\defeq\, \big\{\,\set{(x,i), (y,j)} \,:\, x\in A, \, y \in B,\, \set{x,y} \in E\,\big\}.
        \]
        In other words, $G^*$ is obtained by ``blowing up'' $G$: each vertex in $A$ (resp.~$B$) is replaced by a set of size $p$ (resp.~$q$), and every edge is replaced by a complete bipartite graph between the corresponding sets. Note that we have a natural $(pq)$-to-one projection map $\pi^* \colon E^* \to E$ given by
        \[
            \pi^*\big(\set{(x,i), (y,j)}\big) \,\defeq\, \set{x,y}.
        \]

        \begin{claim}\label{claim:now_Hall}
            The bipartition $(A^*, B^*)$ of $G^*$ satisfies Hall's condition. 
        \end{claim}
        \begin{claimproof}
           Suppose $S^* \subseteq A^*$ is a finite set. Let \[S \defeq \set{x \in A \,:\, (x,i) \in S^* \text{ for some $i$}}.\] Notice that $S^* \subseteq S \times \set{1,\ldots, p}$ and
            $
                N_{G^*}(S^*) = N_G(S) \times \set{1,\ldots, q}
            $.
            Hence, using the assumption on $G$ and the inequality $p/q \leq c$, we obtain the desired bound
            \[
                |N_{G^*}(S^*)| \,=\, q \,|N_G(S)| \,\geq\, q c \,|S| \, \geq\, \frac{qc}{p} \, |S^*| \,\geq\, |S^*|. \qedhere
            \]
        \end{claimproof}
        
        Fix a large integer $n \in \N^+$; it will become clear later how large $n$ needs to be exactly. 
        Since $\asi(G^*) = \asi(G) = s < \infty$, Lemma~\ref{lemma:high_cover} provides Borel matchings $M_1$, \ldots, $M_n$ in $G^*$ that cover every vertex in $A^*$ at least $n - s$ times. The graph $G' = (V, E', \pi)$ is now defined as follows. Let $E'$ be the disjoint union of the matchings $M_1$, \ldots, $M_n$; explicitly,
        \[
            E' \,\defeq\, (M_1 \times \set{1}) \sqcup \ldots \sqcup (M_n \times \set{n}).
        \]
        The map $\pi \colon E' \to [V]^2$ is given by composing the natural projections
        \[
            E' \,\to\, M_1 \cup \ldots \cup M_n \,\to\, E.
        \]
        That is to say, for all $1 \leq j \leq n$  and $e^* \in M_j$, we let $\pi(e^*,j) \defeq \pi^*(e^*)$. 
        
        By construction, the vertex set of $G'$ is $V$, and its underlying simple graph is a subgraph of $G$. It remains to check condition~\ref{item:3_unbalanced}.
        Since each $M_j$ is a matching in $G^*$, for any $y \in B$, there are at most $q$ edges in $M_j$ incident to a vertex of the form $(y,i) \in B^*$, which shows that $\deg_{G'}(y) \leq qn$. 
        On the other hand, for any $x \in A$ and $1 \leq i \leq p$, there exist at least $n - s$ indices $1 \leq j \leq n$ such that the vertex $(x,i) \in A^*$ is covered by $M_j$, and hence $\deg_{G'}(x) \geq p(n-s)$. 
        Since $p/q > c'$, we can choose $n$ so that $p(n-s) > c'qn$, completing the proof.
    \end{scproof}

    Theorem~\ref{theo:asi_two_conjectures} now easily follows:

    \begin{scproof}[ of Theorem~\ref{theo:asi_two_conjectures}]
        Suppose Conjecture~\ref{conj:asi} holds for unbalanced bipartitions, i.e., the following statement is true:
        \begin{quote}
            \textsl{If a bipartition $(A,B)$ of a locally finite Borel graph $G$ is unbalanced and $\asi(G) < \infty$, then $G$ has a Borel matching covering $A$.}
        \end{quote}
        Now let $G$ be a locally finite Borel graph such that $\asi(G) < \infty$ and let $(A,B)$ be a Borel bipartition of $G$ that has combinatorial expansion. We want to find a Borel matching in $G$ that covers $A$. We may replace $G$ by its underlying simple graph to arrange that $G$ itself is simple. Write $G = (V,E)$. Let $c > 1$ be such that for every finite set $S \subseteq A$, we have $|N_G(S)| \geq c\,|S|$. Fix any $1 < c' < c$. By Lemma~\ref{lemma:make_unbalanced}, there is a locally finite Borel graph $G' = (V, E', \pi)$ such that:
        \begin{enumerate}[label=\ep{\normalfont\arabic*}]
            \item the vertex set of $G'$ is $V$,
            
            \item the underlying simple graph of $G'$ is a subgraph of $G$, and 
            
            \item the bipartition $(A,B)$ in $G'$ is of type $(\high, \low)$ for some $\high$, $\low$ such that $\high > c' \,\low$.
        \end{enumerate}
        Since $c' > 1$, the bipartition $(A,B)$ in $G'$ is unbalanced.  As $\asi(G') \leq \asi(G) < \infty$, $G'$ has a Borel matching $M \subseteq E'$ covering $A$ by our assumption. Then $\pi(M)$ is a Borel matching in $G$ that covers $A$, as desired. (The set $\pi(M)$ is Borel because the map $\pi$ is finite-to-one.)
    \end{scproof}

    A similar argument allows us to derive Theorem~\ref{theo:hyperfinite} from Theorem~\ref{theo:measurable}:

    \begin{scproof}[ of Theorem~\ref{theo:hyperfinite}]
        Let $G = (V,E)$ be a (without loss of generality) simple locally finite Borel graph and let $(A,B)$ be a Borel bipartition of $G$ with combinatorial expansion. Suppose $G$ is hyperfinite and let $\mu$ be a probability measure on $V$. We seek a Borel matching in $G$ that covers $\mu$-almost every vertex in $A$. 
        By Theorem~\ref{theo:asi_ae}, 
        we can remove a $\mu$-null $G$-invariant Borel subset from $V$ to arrange that $\asi(G) \leq 1$. As in the proof of Theorem~\ref{theo:asi_two_conjectures} presented above, Lemma~\ref{lemma:make_unbalanced} then can be used to find a locally finite Borel graph $G' = (V,E',\pi)$ such that
        \begin{enumerate}[label=\ep{\normalfont\arabic*}]
            \item the vertex set of $G'$ is $V$,
            
            \item the underlying simple graph of $G'$ is a subgraph of $G$, and 
            
            \item the bipartition $(A,B)$ in $G'$ is unbalanced.
        \end{enumerate}
        By Theorem~\ref{theo:measurable}, there is a Borel matching $M \subseteq E'$ that covers $\mu$-almost every vertex in $A$, and then the desired matching in $G$ is $\pi(M)$.  
    \end{scproof}

    To conclude this subsection, we establish the $\asi = 1$ case of Theorem~\ref{theo:asi} via a modification of the proof of Lemma~\ref{lemma:make_unbalanced}.

    \begin{scproof}[ of Theorem~\ref{theo:asi} in the $\asi = 1$ case]
        Let $G$ be a locally finite Borel graph with a Borel bipartition $(A,B)$ such that for every finite set $S \subseteq A$, we have $|N_G(S)| \geq 2\,|S|$. 
        Suppose $\asi(G) \leq 1$. We need to find a Borel matching in $G$ covering $A$. 

        We may assume $G = (V,E)$ is a simple graph. Form an auxiliary simple graph $G^*$ by ``duplicating'' every vertex in $A$; that is, $G^* \defeq (V^*, E^*)$, where we let
        \[
            A^* \,\defeq\, A \times \set{0,1}, \qquad V^* \,\defeq\, A^* \sqcup B, \qquad \text{and} \qquad E^* \,\defeq\, \big\{\,\set{(x,i), y} \,:\, i \in \set{0,1}, \, \set{x,y} \in E\,\big\}.
        \]
        We also have a two-to-one projection map $\pi^* \colon E^* \to E$ given by
        \[
            \pi^*\big(\set{(x,i), y} \big) \,\defeq\, \set{x,y}.
        \]
        Fix a partition $A = A_0 \sqcup A_1$ of $A$ into two Borel $G^4$-finite sets, which exists since $\asi(G) \leq 1$. For $i \in \set{0,1}$, let $B_i \defeq N_G(A_i)$, and note that the set $A_i \cup B_i$ is $G$-finite. 
        As in the proof of Lemma~\ref{lemma:make_unbalanced}, we observe that the bipartition $(A^*, B^*)$ of $G^*$ satisfies Hall's condition, so, by Theorem~\ref{theo:Hall}, there exists a (possibly non-Borel) matching in $G^*$ that covers $A^*$. 
        Thus, by Theorem~\ref{theo:smooth}, there exists a Borel matching $M_i$ in the graph $G^*[(A_i \times \set{0,1}) \cup B_i]$ covering every vertex in $A_i \times \set{0,1}$. Define
        \[
            F_i \,\defeq\, \pi^*(M_i) \,\subseteq\, E.
        \]
        The following properties of the set $F_i$ are immediate from its construction:
        \begin{itemize}
            \item Every vertex in $A_i$ is incident to exactly two edges in $F_i$.
            \item Every vertex in $B_i$ is incident to at most one edge in $F_i$.
            \item The vertices in $V \setminus (A_i \cup B_i)$ are not incident to any edges in $F_i$.
        \end{itemize}
        Now we consider the spanning subgraph $G'$ of $G$ with edge set $F_0 \cup F_1$. The above statements imply that the bipartition $(A,B)$ in $G'$ is of type $(2,2)$, so $G'$ has a (not necessarily Borel) matching that covers $A$. 
        To finish the proof, we need one more simple observation:

        \begin{claim}\label{claim:small_paths}
            The graph $G'$ is component-finite.
        \end{claim}
        \begin{claimproof}
            If not, then there is an infinite path in $G'$, which we can represent by an infinite sequence $(x_n,y_n)_{n \in \N}$ of distinct vertices, where $x_n \in A$, $y_n \in B$, and $\set{x_n, y_n}$, $\set{y_n, x_{n+1}} \in F_0 \cup F_1$. In particular, every vertex $y_n$ is incident to $2$ edges in $F_0 \cup F_1$. This can only happen if $y_n$ belongs to $B_0 \cap B_1$, i.e., if it has neighbors in both $A_0$ and $A_1$. Since $\dist_G(y_n,y_{n+1}) = 2$ for all $n \in \N$, the vertices $y_0$, $y_1$, \ldots{} are all adjacent to the same component of the graph $G^4[A_0]$. But that is impossible, since the graph $G^4[A_0]$ is component-finite and $G$ is locally finite.
        \end{claimproof}

        Thanks to Claim~\ref{claim:small_paths}, we can apply Theorem~\ref{theo:smooth} to the graph $G'$ and obtain a Borel matching in $G'$ covering $A$. As $G' \subseteq G$, this completes the proof.
    \end{scproof}

    \subsection{Independent complete sections, the Local Lemma, and 
    Theorem~\ref{theo:asi}}\label{subsec:asi_general}

    The proof of Theorem~\ref{theo:asi} for graphs with asymptotic separation index $1$ given in the previous subsection is short and completely elementary. Unfortunately, it does not seem to generalize to any higher values of $\asi$. The argument we present here, whose only assumption is that $\asi < \infty$, is considerably more involved and relies on fairly sophisticated combinatorial tools. As a reward for this extra work, we in fact obtain a much stronger and more general result concerning independent complete sections. We expect this general result to 
    have applications to various other problems beyond matchings. 
    Note that it is of significant interest for graphs with $\asi = 1$ as well, and we do not see any way to simplify its proof even in that setting.
    
    A \emphd{complete section} for a family $\mathcal{F}$ of pairwise disjoint sets is a set $X$ such that $X \cap F \neq \0$ for all $F \in \mathcal{F}$. The following is a classical theorem in extremal combinatorics: 

    \begin{theorem}[{Alon \cite{LinArbor}/Fellows \cite{Fellows}, Haxell \cite{Haxell}}]\label{theo:Alon}
        Let $G = (V,E,\pi)$ be a graph of maximum degree at most $d \in \N^+$, 
        and let $\mathcal{F}$ be a family of pairwise disjoint subsets of $V$. If $|S| \geq 2d$ for all $S \in \mathcal{F}$, then there exists a complete section $X \subseteq V$ for $\mathcal{F}$ that is independent in $G$.
    \end{theorem}

    A ``qualitative'' form of Theorem~\ref{theo:Alon}, with the bound $|S| \geq 2d$ replaced by $|S| \geq C\,d$ for some constant $C > 0$, is sufficient for most applications. In this form, the theorem was established by Alon \cite{LinArbor} and independently Fellows \cite{Fellows} using a probabilistic argument. A simplified version of the proof, presented in \cite{Alon2}, gives the value $C = 2\mathsf{e}$. This was lowered to $C = 2$ by Haxell \cite{Haxell,Haxell2} via a different, purely combinatorial approach (see also \cite{top1,top2} for alternative proofs based on topological ideas). The value $C = 2$ is generally best possible \cite{BES,SzT}, although it can be further reduced in various special cases; see, e.g., \cite{LohSudakov}.

    \begin{remark}
        A related notion considered in the literature is that of an \emph{independent transversal}---i.e., a set that meets each member of $\mathcal{F}$ in {exactly} one point. For our purposes, working with complete sections will be somewhat more natural, although in most cases of interest independent transversals and independent complete sections are interchangeable.  
    \end{remark}


    We are interested in Borel versions of Theorem~\ref{theo:Alon}. To phrase our results, it will be convenient to assume that the family $\mathcal{F}$ has the form  $\mathcal{F} = \set{C \,:\, \text{$C$ is a component of $G$}}$ for some graph $G$. A~complete section for such a family $\mathcal{F}$, i.e., a set of vertices that meets every component of $G$, will be called a \emphd{complete section for $G$}. The following was proved by the first and third named authors: 

    \begin{theorem}[{AB--FW \cite[Thm.~1.53]{bw_lll}}]\label{theo:Borel_section_log}
        Let $H$ be a locally finite Borel graph with $\asi(H) < \infty$ and let $G_1$, $G_2$ be two spanning Borel subgraphs of $H$. Assume that:
        \begin{itemize}
        \item the maximum degree of $G_1$ is at most $d \in \N$, where $d \geq 2$, and
        \item every component of $G_2$ contains at least $2^{25}\,d\log d$ vertices.
        \end{itemize}
        Then there exists a Borel complete section for $G_2$ that is independent in $G_1$.
    \end{theorem}

    Notice that the lower bound on the size of the components of $G_2$ required by Theorem~\ref{theo:Borel_section_log} is (slightly) superlinear in $d$, in contrast to the linear bound in Theorem~\ref{theo:Alon}. The authors of \cite{bw_lll} conjectured that a linear bound is sufficient \cite[Conj.~1.54]{bw_lll}. We prove this conjecture: 

    \begin{tcolorbox}
    \begin{theorem}[Borel independent complete sections]\label{theo:Borel_section}
        There is a constant $C > 0$ with the following property. Let $H$ be a locally finite Borel graph with $\asi(H) < \infty$ and let $G_1$, $G_2$ be two spanning Borel subgraphs of $H$. Assume that:
        \begin{itemize}
        \item the maximum degree of $G_1$ is at most $d \in \N^+$, and
        \item every component of $G_2$ contains at least $C\,d$ vertices.
        \end{itemize}
        Then there exists a Borel complete section for $G_2$ that is independent in $G_1$.
    \end{theorem}
    \end{tcolorbox}

    Before proving Theorem~\ref{theo:Borel_section}, let us explain how it quickly yields Theorem~\ref{theo:asi}.

    \begin{scproof}[ of Theorem~\ref{theo:asi}]
        Let $C > 0$ be the constant from Theorem~\ref{theo:Borel_section}. We first prove a special case of Theorem~\ref{theo:asi} for graphs with unbalanced bipartitions:
        
        \begin{claim}\label{claim:asi_unbalanced}
            Suppose $G = (V, E, \pi)$ is a locally finite Borel graph with $\asi(G) < \infty$. Let $(A,B)$ be a bipartition of $G$ of type $(\high,\low)$ such that $\high > C\,\low$. Then $G$ has a Borel matching that covers $A$.
        \end{claim}
        \begin{claimproof}
            We may, of course, assume $A \neq \0$, and hence also $E \neq \0$, and thus $\low \geq 1$. Define simple Borel graphs $H$, $G_1$, and $G_2$ with vertex set $E$ and the following edge sets:
        \begin{itemize}
            \item the edge set of $H$ is $\set{\set{e,e'} \in [E]^2 \,:\, \text{$e$ and $e'$ have a common endpoint}}$,
            \item the edge set of $G_1$ is $\set{\set{e,e'} \in [E]^2 \,:\, \text{$e$ and $e'$ have a common endpoint in $B$}}$,
            \item the edge set of $G_2$ is $\set{\set{e,e'} \in [E]^2 \,:\, \text{$e$ and $e'$ have a common endpoint in $A$}}$.
        \end{itemize}
        It is easy to see that $\asi(H) = \asi(G) < \infty$. Note that $G_1$ and $G_2$ are spanning Borel subgraphs of $H$, the maximum degree of $G_1$ is at most $\low$, and every component of $G_2$ has at least $\high > C\,\low$ elements. Therefore, we may apply Theorem~\ref{theo:Borel_section} to obtain a Borel complete section $X \subseteq E$ for $G_2$ that is independent in $G_1$. In other words, $X$ is a Borel set of edges in $G$ such that every vertex $x$ is incident to \emph{at least one} edge in $X$ if $x \in A$ and \emph{at most one} edge in $X$ if $x \in B$. Using Luzin--Novikov Uniformization (Lemma~\ref{lemma:LNU}), we can pass to a further Borel subset $M \subseteq X$ such that every vertex in $A$ is incident to exactly one edge in $M$. Then $M$ is a Borel matching in $G$ covering $A$, as desired. 
        \end{claimproof}

        The rest of the argument proceeds in the same way as the proof of Theorem~\ref{theo:asi_two_conjectures} given in \S\ref{subsec:reduction}. Fix an arbitrary constant $C' > C$. Suppose $G$ is a (without loss of generality) simple locally finite Borel graph with $\asi(G) < \infty$ and let $(A,B)$ be a Borel bipartition of $G$ such that for every finite subset $S \subseteq A$, we have $|N_G(S)| \geq C' \,|S|$. Using Lemma~\ref{lemma:make_unbalanced}, we can find a Borel graph $G' = (V,E', \pi)$ with the same vertex set as $G$ such that the underlying simple graph of $G'$ is a subgraph of $G$, and the bipartition $(A,B)$ in $G'$ is of type $(\high, \low)$, where $\high > C\,\low$. By Claim~\ref{claim:asi_unbalanced}, $G'$ has a Borel matching $M$ that covers $A$, and then $\pi(M)$ is a Borel matching in $G$ covering $A$.
    \end{scproof}


    The remainder of this subsection is dedicated to the proof of Theorem~\ref{theo:Borel_section}.
    
    The main ingredient in the proof of Theorem~\ref{theo:Borel_section} is the \emph{Lov\'asz Local Lemma}, or the \emph{LLL} for short, which is a powerful and versatile tool from probabilistic combinatorics. Roughly, the LLL guarantees the existence of a function satisfying a given set of combinatorial constraints under certain numerical conditions. In \cite{bw_lll}, the first and third named authors showed that under a finite $\asi$ assumption, similar numerical conditions additionally ensure the desired function can be made Borel.  The paper \cite{bw_lll} also describes several applications of this result, including Theorem~\ref{theo:Borel_section_log}. It is well known that the LLL gives a simple proof of Theorem~\ref{theo:Alon} with a larger constant in place of~$2$; see, e.g., \cite[Prop.~5.5.3]{AS}. Unfortunately, this argument does not (at least at present) directly carry over to the Borel setting, for reasons explained below. Nevertheless, we present here an alternative approach that does apply in the Borel context and allows us to establish Theorem~\ref{theo:Borel_section}. 

    The statement of the LLL is somewhat technical and requires a few definitions. First, we introduce a class of combinatorial problems the LLL provides solutions to. 
    
    \begin{definition}[{Constraint satisfaction problems}]\label{defn:CSP}
		Let $q \in \N^+$ be a positive integer, which we identify with the set $\set{0,1,\ldots,q-1}$. 
		Given a set $X$ and a finite subset $D \subseteq X$, an \emphd{$(X,q)$-constraint} \ep{or simply a \emphd{constraint} if $X$ and $q$ are understood} with \emphd{domain} $D$ is a set $B \subseteq q^D$ of mappings $D \to q$. If $B$ is a constraint with domain $D$, we write $\dom(B) \defeq D$. 
        A function $f \colon X \to q$ \emphd{violates} a constraint $B$ 
        if $\rest{f}{\dom(B)} \in B$, and \emphd{satisfies} $B$ otherwise.
			
		A \emphd{constraint satisfaction problem} \ep{a \emphd{CSP} for short} $\mathcal{B}$ on a set $X$ with range $q$, in symbols \[\mathcal{B} \colon X \to^? q,\] is a set of $(X,q)$-constraints. A \emphd{solution} to a CSP $\mathcal{B} \colon X \to^? q$ is a function $f \colon X \to q$ that satisfies every constraint $B \in \mathcal{B}$. We say that $\mathcal{B}$ is \emphd{satisfiable} if it has a solution.
	\end{definition}
	
	In other words, each constraint in a CSP $\mathcal{B} \colon X \to^? q$ specifies a set of finite ``forbidden'' (or ``bad'') patterns that are not allowed to appear in a solution $f \colon X \to q$. This is a very flexible framework, and many problems of combinatorial interest can naturally be formalized as CSPs. Consequently, it is desirable to have general sufficient conditions showing that a given CSP is satisfiable. The LLL provides one such condition, phrased in probabilistic terms.

    \begin{definition}[$\mathsf{p}(\mathcal{B})$ and $\mathsf{d}(\mathcal{B})$]
        Let $\mathcal{B} \colon X \to^? q$. The \emphd{probability} of a constraint $B \in \mathcal{B}$ is
        \[
            \P[B] \,\defeq\, \frac{|B|}{q^{|\dom(B)|}}.
        \]
        In other words, $\P[B]$ is the probability that a uniformly random map $\dom(B) \to q$ violates $B$. Let
        \[
            \mathsf{p}(\mathcal{B}) \,\defeq\, \sup_{B \in \mathcal{B}} \P[B].
        \]
        The \emphd{maximum dependency degree} of the CSP $\mathcal{B}$ is the quantity
        \[
            \mathsf{d}(\mathcal{B}) \,\defeq\, \sup_{B \in \mathcal{B}} \big|\big\{B' \in \mathcal{B} \,:\, B' \neq B \text{ and } \dom(B') \cap \dom(B) \neq \0\big\}\big|.
        \]
    \end{definition}

    \begin{theorem}[{Lov\'asz Local Lemma \cite{EL,SpencerRamsey}; \cite[Cor.~5.1.2]{AS}}]\label{theo:LLL}
        If $\mathcal{B}$ is a CSP such that 
        \[
            \mathsf{e} \cdot \mathsf{p}(\mathcal{B}) \cdot (\mathsf{d}(\mathcal{B}) + 1) \,\leq\, 1,
        \]
        where $\mathsf{e} = 2.718\ldots$ is the base of the natural logarithm, then $\mathcal{B}$ is satisfiable.
    \end{theorem}

    The above statement is really a special case of the LLL; in its full generality, the LLL can be applied to abstract probability spaces and not just spaces of functions $X \to q$. That being said, the above framework encompasses virtually all standard applications and is often viewed as the ``right one'' for constructive considerations; see, e.g., \cite{Beck,MT,decompose7,KolipakaSzegedy,RSh,subexp1,subexp2,BernshteynDistributed,BerShift,BGRDeterministicLLL,CPS}. We should also note that in the combinatorics literature, the CSP $\mathcal{B}$ is usually assumed to be finite. However, Theorem~\ref{theo:LLL} holds without any restrictions on the cardinality of $\mathcal{B}$ by compactness; see, e.g., \cite[proof of Thm.~5.2.2]{AS}. 

    To prove Theorem~\ref{theo:Borel_section}, we need a version of the LLL that yields a \emph{Borel} function satisfying every constraint. In this context, it is natural to assume that the CSP $\mathcal{B}$ is itself Borel in an appropriate sense. In brief, if $X$ is a standard Borel space and $q \in \N^+$, then each $(X,q)$-constraint is a finite set. As a result, the space $\mathsf{Const}(X,q)$ of all $(X,q)$-constraints carries a Borel $\sigma$-algebra, allowing us to call a CSP \emphd{Borel} if it is a Borel subset of $\mathsf{Const}(X,q)$. See \cite[\S2.1]{bw_lll} for details.

    Unfortunately, a Borel CSP $\mathcal{B}$ as in Theorem~\ref{theo:LLL} may fail to have a Borel solution unless the condition $\mathsf{e} \,\mathsf{p}(\mathcal{B}) \, (\mathsf{d}(\mathcal{B}) + 1) \leq 1$ is replaced by a much stronger one, namely $\mathsf{p}(\mathcal{B})\, 2^{\mathsf{d}(\mathcal{B})} < 1$, which is too restrictive to be widely applicable; see \cite[Thm.~1.26]{bw_lll} and the following paragraphs for a discussion of this issue. 
    That being said, for most applications the inequality $\mathsf{e} \,\mathsf{p}(\mathcal{B}) \, (\mathsf{d}(\mathcal{B}) + 1) \leq 1$ can be replaced by any 
    \emph{polynomial criterion}, i.e., a bound of the form $\mathsf{p}(\mathcal{B}) \, \leq\, c \,\mathsf{d}(\mathcal{B})^{-k}$ for some constants $c > 0$ and $k \in \N^+$ \cite{Beck, BernshteynDistributed, decompose7, CP, GHK}. A Borel version of the LLL with a polynomial criterion does exist, provided that 
    a certain graph defined using $\mathcal{B}$ has finite asymptotic separation index.

    \begin{definition}[Graphs associated to CSPs]
        Let $\mathcal{B} \colon X \to^? q$ be a CSP. Then $G_\mathcal{B}$ is the simple graph with vertex set $X$ in which distinct vertices $x$, $y \in X$ are adjacent if and only if there is a constraint $B \in \mathcal{B}$ such that $\set{x,y} \subseteq \dom(B)$.
    \end{definition}

    \begin{theorem}[{AB--FW \cite[Thm.~1.31]{bw_lll}}]\label{theo:BLLL}
        Let $\mathcal{B} \colon X \to^? q$ be a Borel CSP on a standard Borel space $X$. Suppose that $\asi(G_\mathcal{B}) < \infty$. If
        \[
            2^{15} \cdot \mathsf{p}(\mathcal{B}) \cdot (\mathsf{d}(\mathcal{B}) + 1)^8 \,\leq\, 1,
        \]
        then $\mathcal{B}$ has a Borel solution $f \colon X \to q$.
    \end{theorem}

    A major open question is whether the standard LLL condition $\mathsf{e} \,\mathsf{p}(\mathcal{B}) \, (\mathsf{d}(\mathcal{B}) + 1) \leq 1$ is already sufficient for the conclusion of Theorem~\ref{theo:BLLL}. Although in most combinatorial arguments, the classical LLL can be replaced by Theorem~\ref{theo:BLLL} without much difficulty, there are a few combinatorial problems for which the exact form of the LLL bound is crucial. Finding independent complete sections happens to be one such problem (arguably the most important one). This is the reason Alon's simple LLL-based proof cannot be easily adapted to derive Theorem~\ref{theo:Borel_section}. Instead, we prove Theorem~\ref{theo:Borel_section} via an inductive approach, with Theorem~\ref{theo:Borel_section_log} serving as the base case.

    \begin{scproof}[ of Theorem~\ref{theo:Borel_section}]
        The case $d=1$ in Theorem~\ref{theo:Borel_section} follows from the case $d=2$ (at the cost of replacing the constant $C$ by $2C$), so we shall restrict our attention to $d \geq 2$ from now on, which will allow us to apply Theorem~\ref{theo:Borel_section_log}. 
        
        To facilitate an inductive argument, we establish a slightly more technical statement. Let $C > 0$ be a sufficiently large constant (it will become clear in due course how large $C$ needs to be). Let $H$ be a locally finite Borel graph with $\asi(H) < \infty$ and let $G_1$, $G_2$ be two spanning Borel subgraphs of $H$. Suppose that, for some integer $d \geq 2$:
        \begin{itemize}
            \item the maximum degree of $G_1$ is at most $d$, and
            \item the number of vertices in every component of $G_2$ is at least
            \begin{equation}\label{eq:factor}
                C\,d \,-\, d^{2/3}. 
            \end{equation}
        \end{itemize}
        We will prove that, under these assumptions, there exists a Borel complete section for $G_2$ that is independent in $G_1$. This will, of course,  imply the statement of Theorem~\ref{theo:Borel_section}. 
        
        We proceed by induction on $d$, i.e., we assume that the statement holds when $d$ is replaced by any strictly smaller value. Let the common vertex set of $H$, $G_1$, $G_2$ be $V$. Without loss of generality, we may assume the graphs 
        $H$, $G_1$, $G_2$ are simple. 
        Let $k \in \N^+$ be the ceiling of the quantity in \eqref{eq:factor}, so that every component of $G_2$ has at least $k$ vertices. As explained in the first paragraph of \cite[proof of Thm.~1.53]{bw_lll}, it is enough to consider the case when every component of $G_2$ has \emph{exactly} $k$ vertices. Furthermore, we may replace the graphs $H$ and $G_2$ by their $k$-th powers and assume that every component of $G_2$ is a $k$-clique. 

        If $k \geq 2^{25} d \log d$, then we are done by Theorem~\ref{theo:Borel_section_log}, so we may assume $k < 2^{25} d \log d$. This implies that $C < 2^{25} \log d + 1$. 
        In particular, by making $C$ large enough, we can arrange that $d$ exceeds any given constant. 
        Define
        \[
            d' \,\defeq\, \left\lfloor \frac{d}{2} \,+\, 4\sqrt{d \log k} \right\rfloor \qquad \text{and} \qquad k' \,\defeq\, \left\lceil \frac{k}{2} \,-\, 4 \sqrt{k \log k} \right\rceil.
        \]
        The key idea is to use the following observation, which is an easy consequence of Theorem~\ref{theo:BLLL}:

        \begin{claim}\label{claim:reduce}
            There exists a Borel set $Z \subseteq V$ such that:
            \begin{itemize}
                \item in the graph $G_1$, every vertex $x \in V$ has at most $d'$ neighbors in $Z$, and
                \item every component of $G_2$ meets $Z$ in at least $k'$ vertices.
            \end{itemize}
        \end{claim}
        \begin{claimproof}
            For brevity, given a vertex $x \in V$, we write $N_1(x) \defeq N_{G_1}(x)$.  We will show that there exists a Borel function $f \colon V \to 2$ with the following properties:
            \begin{itemize}
                \item for all $x \in V$, $|f^{-1}(1) \cap N_1(x)| \leq d'$, and
                \item for every component $C$ of $G_2$, $|f^{-1}(1) \cap C| \geq k'$.
            \end{itemize}
            Setting $Z \defeq f^{-1}(1)$ will then complete the proof of the claim.
            
            The problem of finding the desired function $f$ can naturally be represented as a CSP $\mathcal{B} \colon V \to^? 2$. Namely, for each vertex $x \in V$, we define a constraint $B_x$ with domain $N_1(x)$ by
            \[
                B_x \,\defeq\, \big\{\,\phi \colon N_1(x) \to 2 \ :\  |\phi^{-1}(1)| > d'\,\big\}.
            \]
            Similarly, for a component $C$ of $G_2$, we let $B_C$ be the constraint with domain $C$ given by
            \[
                B_C \,\defeq\, \big\{\,\psi \colon C \to 2 \ :\  |\psi^{-1}(1)| < k'\,\big\}.
            \]
            Our goal is to find a Borel solution $f \colon V \to 2$ to the CSP
            \[
                \mathcal{B} \,\defeq\, \set{B_x \,:\, x \in V} \,\cup\, \set{B_C \,:\, \text{$C$ is a component of $G_2$}}.
            \]
            It is routine to check that the CSP $\mathcal{B}$ is Borel. Furthermore, $G_\mathcal{B}$ is a subgraph of $H^2$, so $\asi(G_\mathcal{B}) < \infty$. Now we need to bound $\mathsf{p}(\mathcal{B})$ and $\mathsf{d}(\mathcal{B})$.

            We start with $\mathsf{d}(\mathcal{B})$. Consider a constraint of the form $B_x$ for $x \in V$. Since $G_1$ has maximum degree at most $d$, there are at most $d^2$ vertices $y \in V$ such that $\dom(B_y) \cap \dom(B_x) = N_1(y) \cap N_1(x) \neq \0$. Also, the set $\dom(B_x) = N_1(x)$ has size at most $d$, so it can intersect at most $d$ components of $G_2$. Therefore, in total there are at most $d^2 + d$ constraints $B' \in \mathcal{B}$ with $\dom(B') \cap \dom(B_x) \neq \0$. Now take any constraint of the form $B_C$ for a component $C$ of $G_2$. If $\dom(B') \cap \dom(B_C) \neq \0$ for some $B' \in \mathcal{B} \setminus \set{B_C}$, then $B' = B_y$ for a vertex $y \in V$ that has a $G_1$-neighbor in $C$. Since $|C| = k$, there are at most $d k$ such vertices. To summarize, we have
            \begin{equation}\label{eq:d}
                \mathsf{d}(\mathcal{B}) \,\leq\, \max \set{d^2 + d, \, d k} \,<\, k^2.
            \end{equation}

            Now we turn our attention to $\mathsf{p}(\mathcal{B})$. Define 
            \[
                \lambda \,\defeq\, 4 \sqrt{\log k}.
            \]
            Consider a constraint $B_x$ for some $x \in V$. Note that if $\phi \colon N_1(x) \to 2$ is chosen uniformly at random, then $|\phi^{-1}(1)|$ is a binomial random variable with distribution $\mathsf{Bin}(n,1/2)$, where $n \defeq |N_1(x)| \leq d$. Thus, the standard Chernoff bound \cite[Thm.~A.1.1]{AS} yields
            \[
                \P[B_x] \,=\, \P\left[\,|\phi^{-1}(1)| > d'\,\right] \,=\, \P\left[\,|\phi^{-1}(1)| > d/2 + \lambda \sqrt{d}\,\right] \,\leq\, \mathsf{e}^{-2\lambda^2}.
            \]
            Similarly, if $C$ is a component of $G_2$ and $\psi \colon C \to 2$ is chosen uniformly at random, then the random variable $|\psi^{-1}(1)|$ has distribution $\mathsf{Bin}(k,1/2)$, and hence
            \[
                \P[B_C] \,=\, \P\left[\,|\psi^{-1}(1)| < k'\,\right] \,=\, \P\left[\,|\psi^{-1}(1)| < k/2 - \lambda \sqrt{k}\,\right] \,\leq\, \mathsf{e}^{-2\lambda^2}.
            \]
            In conclusion, we have the overall bound
            \begin{equation}\label{eq:p}
                \mathsf{p}(\mathcal{B}) \,\leq\, \mathsf{e}^{-2\lambda^2} \,=\, k^{-32}.
            \end{equation}

            Finally, we put \eqref{eq:d} and \eqref{eq:p} together and observe that
            \[
                2^{15} \cdot \mathsf{p}(\mathcal{B}) \cdot (\mathsf{d}(\mathcal{B}) + 1)^8 \,\leq\, 2^{15} \cdot k^{-32} \cdot k^{16} \,=\, 2^{15} \cdot k^{-16} \,<\, 1.
            \]
            Hence, Theorem~\ref{theo:BLLL} provides a Borel solution to $\mathcal{B}$, as desired.
        \end{claimproof}

        Let $Z \subseteq V$ be a Borel set given by Claim~\ref{claim:reduce}. Notice that, since $k < 2^{25} d \log d$, we have
        \[
            d' \,\leq\, \frac{d}{2} \,+\, 4\sqrt{d \log k} \,=\, \frac{d}{2} \,+\, O\big(\sqrt{d \log d}\big) \,<\, d, 
        \]
        assuming $d$ is large enough. 
        It remains to show that
        \begin{equation}\label{eq:k_bound}
            k' \,\geq\, C\,d' \,-\, (d')^{2/3}.
        \end{equation}
        Indeed, if \eqref{eq:k_bound} holds, we can apply the inductive hypothesis to find a Borel complete section $X \subseteq Z$ for the graph $G_2[Z]$ that is independent in $G_1[Z]$. (Note that every connected component of $G_2[Z]$ has at least $k'$ vertices because all components of $G_2$ are cliques.) Clearly, $X$ is also a Borel complete section for $G_2$ that is independent in $G_1$, as desired.

        To prove \eqref{eq:k_bound}, 
        it suffices to argue that
        \[
            \frac{C\,d}{2} \,-\, \frac{d^{2/3}}{2} \,-\, 4\sqrt{k \log k} \,\geq\, \frac{C\,d}{2} \,+\, 4\,C\sqrt{d \log k} \,-\, (d')^{2/3}. 
        \]
        After rearranging the terms, this becomes
        \begin{equation}\label{eq:roots}
            2 \,(d')^{2/3} \,-\, d^{2/3} \,\geq\, 8\sqrt{k \log k} \,+\, 8 \,C \sqrt{d \log k}. 
        \end{equation}
        Since $d' \geq d/2$, the left-hand side of \eqref{eq:roots} can be bounded below as
        \[
            2 \,(d')^{2/3} \,-\, d^{2/3} \,\geq\, (2^{1/3} - 1) \,d^{2/3}
            \,=\, \Theta(d^{2/3}).
        \]
        On the other hand, since $C < 2^{25} \log d + 1$ and $k < 2^{25} d \log d$, the right-hand side of \eqref{eq:roots} is
        \[
            8\sqrt{k \log k} \,+\, 8 \,C \sqrt{d \log k} \,\leq\, O\left(\sqrt{d}\, (\log d)^{3/2}\right). 
        \]
        It follows that \eqref{eq:roots} is satisfied for all sufficiently large $d$, and we are done.
    \end{scproof}

    With minimal modifications, the above proof also yields the measure-theoretic variant of Theorem~\ref{theo:Borel_section}, which we record here \ep{although it will not be used later}. If $G = (V,E,\pi)$ is a Borel graph and $\mu$ is a probability measure on $V$, we call a set $X \subseteq V$ a \emphd{$\mu$-complete section} for $G$ if for $\mu$-almost every vertex $x \in V$, the connected component of $x$ in $G$ meets $X$.

    \begin{tcolorbox}
    \begin{theorem}[Independent complete sections almost everywhere]\label{theo:measure_section}
        There is a constant $C > 0$ with the following property. Let $G_1$, $G_2$ be two locally finite Borel graphs with the same vertex set $V$. Assume that:
        \begin{itemize}
        \item the maximum degree of $G_1$ is at most $d \in \N^+$, and
        \item every component of $G_2$ contains at least $C\,d$ vertices.
        \end{itemize}
        Then, for any probability measure $\mu$ on $V$, there exists a Borel $\mu$-complete section for $G_2$ that is independent in $G_1$.
    \end{theorem}
    \end{tcolorbox}

    The proof of Theorem~\ref{theo:measure_section} verbatim repeats the proof of Theorem~\ref{theo:Borel_section}, except that Theorem~\ref{theo:BLLL}, i.e., the Borel version of the LLL, is replaced by its measurable analogue, \cite[Thm.~2.20(i)]{BernshteynDistributed}. Note that Theorem~\ref{theo:Borel_section_log} is also proved in \cite{bw_lll} via Theorem~\ref{theo:BLLL}, and its proof similarly works in the measure-theoretic setting if \cite[Thm.~2.20(i)]{BernshteynDistributed} is used instead.

    \begin{remark}
        By invoking the distributed algorithm for the LLL due to Rozho\v{n} and Ghaffari \cite[Cor.~3.9]{decompose6}, it is in fact possible to obtain a randomized distributed algorithm for the independent complete section problem on bounded degree graphs with $\operatorname{poly}(\log \log n)$ running time (we introduce distributed algorithms in the next subsection). This can be seen as a common generalization of Theorems~\ref{theo:Borel_section} and \ref{theo:measure_section}, thanks to \cite[Thm.~1.42]{bw_lll} and \cite[Thm.~2.14(i)]{BernshteynDistributed}, respectively.
    \end{remark}

    \subsection{List-colorings and the proof of Theorem~\ref{theo:asi_simple}}\label{subsec:asi_simple}

    In this subsection we establish Theorem~\ref{theo:asi_simple}, which improves Theorem~\ref{theo:asi} in the case of simple and, more generally, $k$-simple graphs $G$ with an unbalanced bipartition by replacing the constant $C$ with a $1 + o(1)$ factor. The main idea is to observe that the matching problem in simple bipartite graphs can be seen as a special case of the so-called \emph{list-coloring problem}, which can be addressed using powerful techniques from graph theory and computer science.

    List-coloring was introduced in the 1970s by Vizing \cite{Viz} and, independently, Erd\H{o}s, Rubin, and Taylor \cite{ERT}; for textbook introductions, see \cites[\S5.4]{Diestel}[\S8.4]{West}. A \emphd{list assignment} for a graph $G = (V,E,\pi)$ is a function $L \colon V \to \finset{\N}$, where $\finset{\N}$ is the set of all finite subsets of $\N$. For a vertex $x \in V$, the set $L(x)$ is called the \emphd{list} of $x$, and the elements of $L(x)$ are \emphd{available} to $x$. An \emphd{$L$-coloring} of $G$ is a map $f \colon V \to \N$ such that $f(x) \in L(x)$ for all $x \in V$. An $L$-coloring $f$ is \emphd{proper at a vertex} $x \in V$ if $f(x) \neq f(y)$ for all neighbors $y$ of $x$. An $L$-coloring is \emphd{proper} if it is proper at every vertex. The graph $G$ is \emphd{$L$-colorable} if it has a proper $L$-coloring.
    
    We call $L$ an \emphd{$(\ell,d)$-list assignment} for some $\ell$, $d \in \N$ if it has the following two properties:
    \begin{itemize}
        \item for all $x \in V$, $|L(x)| \geq \ell$,
        \item for all $x \in V$ and $i \in L(x)$, $|\set{y \in N_G(x) \,:\, i \in L(y)}| \leq d$.
    \end{itemize}
    Reed \cite{Reed} conjectured that if $L$ is an $(\ell,d)$-list assignment with $\ell > d$, then $G$ is $L$-colorable. Somewhat surprisingly, this was disproved by Bohman and Holzman \cite{BohHol}. Nevertheless, we have the following asymptotic result due to Reed and Sudakov \cite{ReedSud}:

    \begin{theorem}[{Reed--Sudakov \cite{ReedSud}}]\label{theo:ReedSudakov}
        For each $\epsilon > 0$, there is $d_0 = d_0(\epsilon) \in \N$ such that if $G$ is a graph and $L$ is an $(\ell,d)$-list assignment for $G$ with $\ell > (1+\epsilon) \, d$ and $d \geq d_0$, then $G$ is $L$-colorable.
    \end{theorem}

    In \cite{BernshteynDistributed}, the first named author proved a version of Theorem~\ref{theo:ReedSudakov} for Borel graphs:

    \begin{theorem}[{AB \cite[Thm.~3.13]{BernshteynDistributed}}]\label{theo:BorelReedSudakov}
        For each $\epsilon > 0$, there is $d_0 = d_0(\epsilon) \in \N$ with the following property. Let $G = (V,E,\pi)$ be a Borel graph and let $L$ be a Borel $(\ell,d)$-list assignment for $G$ with $\ell > (1+\epsilon) \,d$ and $d \geq d_0$. If $\mu$ is a probability measure and $\tau$ is a compatible topology on $V$, then:
        \begin{enumerate}[label=\ep{\normalfont\arabic*}]
            \item $G$ has a Borel $L$-coloring that is proper at $\mu$-almost every vertex,
            \item $G$ has a Borel $L$-coloring that is proper at all but a $\tau$-meager set of vertices. 
        \end{enumerate}
    \end{theorem}

    The proof of Theorem~\ref{theo:BorelReedSudakov} makes crucial use of \emph{distributed algorithms}, specifically the \emph{\LOCAL model} of distributed computation. This model was formally introduced by Linial in \cite{Linial} \ep{although there is also some earlier related work, e.g., \cite{ABI,GPSh,Luby}}. For comprehensive overviews, we recommend the book \cite{BE} by Barenboim and Elkin and the recent survey article \cite{Rozhon_survey} by Rozho\v{n}. The \LOCAL model is intended to quantify the difficulty of transforming local data into a global solution to a problem in a large decentralized communication network. Informally, the model views each vertex of a finite graph $G$ as an agent that may pass messages to its neighbors. The length of the messages is unrestricted, and the computational power available to each agent is unlimited; the only resource that is constrained in this model is the number of communication steps. Eventually, every agent must decide on its own part of the output, and the time complexity of a \LOCAL algorithm is the required number of communication rounds. See \cite{BE,Rozhon_survey} for a rigorous definition of this model and various examples of \LOCAL algorithms.

    A line of research initiated by the first named author in \cite{BernshteynDistributed} establishes close links between distributed computing and descriptive combinatorics. It turns out that efficient distributed algorithms for various combinatorial problems on finite graphs can be used, in a completely black-box manner, to obtain well-behaved solutions to these problems on Borel graphs. (Moreover, sometimes this implication can be reversed \cite{grebik.rozhon,btoast,Ber_cont}.) Theorem~\ref{theo:BorelReedSudakov} is a typical instance of this general phenomenon. Chung, Pettie, and Su \cite[\S4.4]{CPS} gave a randomized \LOCAL algorithm that, for an $n$-vertex graph $G$ and a list assignment $L$ as in Theorem~\ref{theo:ReedSudakov}, finds a proper $L$-coloring of $G$ in time $O(d^{-\delta} \log n)$, where $\delta = \delta(\epsilon) > 0$ depends only on $\epsilon$. As explained in \cite[\S3.3]{BernshteynDistributed}, this implies Theorem~\ref{theo:BorelReedSudakov} via a general ``distributed-to-descriptive'' result, namely \cite[Thm.~2.14]{BernshteynDistributed}.
    
    As shown in \cite{bw_lll}, Theorem~\ref{theo:BLLL} implies a version of \cite[Thm.~2.14]{BernshteynDistributed} avoiding the loss of a null/meager set under a finite 
    asymptotic separation index
    assumption. Substituting this result---more specifically \cite[Thm.~3.9]{bw_lll}---in place of \cite[Thm.~2.14]{BernshteynDistributed} in the proof of Theorem~\ref{theo:BorelReedSudakov} yields the following:

    \begin{tcolorbox}
    \begin{theorem}[Borel list-coloring for graphs with $\asi < \infty$]\label{theo:BorelReedSudakov_asi}
        For every $\epsilon > 0$, there is $d_0 = d_0(\epsilon) \in \N$ with the following property. Let $G = (V,E,\pi)$ be a locally finite Borel graph and let $L$ be a Borel $(\ell,d)$-list assignment for $G$ with $\ell > (1+\epsilon) \, d$ and $d \geq d_0$. If $\asi(G) < \infty$, then $G$ has a Borel proper $L$-coloring.
    \end{theorem}
    \end{tcolorbox}

    Modulo replacing \cite[Thm.~2.14]{BernshteynDistributed} by \cite[Thm.~3.9]{bw_lll}, the proof of Theorem~\ref{theo:BorelReedSudakov_asi} is exactly the same as that of Theorem~\ref{theo:BorelReedSudakov}, so there is no need to reproduce it here; the details can be found in \cite[\S3.3]{BernshteynDistributed}. With Theorem~\ref{theo:BorelReedSudakov_asi} in hand, we are ready to verify Theorem~\ref{theo:asi_simple}:

    \begin{scproof}[ of Theorem~\ref{theo:asi_simple}]
        Fix $\epsilon > 0$ and $k \in \N^+$ and let $G = (V,E,\pi)$ be a $k$-simple locally finite Borel graph with a bipartition $(A,B)$ of type $(\high, \low)$. Assume that $\high > (1+\epsilon) \,\low$ and $\low \geq b_0(\epsilon,k)$, where $b_0$ is a sufficiently large function of $\epsilon$ and $k$ (it will become clear during the proof how large $b_0$ needs to be). Supposing that $\asi(G) < \infty$, we seek a Borel matching in $G$ that covers $A$.

        We first consider the case when $k = 1$, i.e., $G$ is a simple graph. As usual, we may then assume that $E \subseteq [V]^2$ and $\pi$ is the inclusion map $E \hookrightarrow [V]^2$. Let $d_0 = d_0(\epsilon) > 0$  be the quantity given by Theorem~\ref{theo:BorelReedSudakov_asi}. We will show that taking $b_0(\epsilon,1) \defeq d_0$ works. 

        Using Corollary~\ref{corl:spaced}, we partition the vertex set of $G$ into countably many Borel $4$-independent sets as $V = \bigsqcup_{i \in \N} U_i$. Let $H \defeq G^2[A]$ and let $L \colon A \to \finset{\N}$ be the list assignment for $H$ given by
        \[
            L(x) \,\defeq\, \set{i \in \N \,:\, \text{$x$ has a $G$-neighbor in $U_i$}}.
        \]
        Each $x \in A$ has at most one $G$-neighbor in every $U_i$. Since $G$ is simple, we conclude that
        \[
            |L(x)| \,=\, |N_G(x)| \,=\, \deg_G(x) \,\geq\, \high.
        \]
        Now let $i \in L(x)$ and let $y$ be the unique vertex in $N_G(x) \cap U_i$. If $z \in A$ is an $H$-neighbor of $x$ with $i \in L(z)$, then, since $U_i$ is $4$-independent in $G$, the unique vertex in $N_G(z) \cap U_i$ must be $y$. In other words, $z \in N_G(y)$, and thus there are at most $\deg_G(y) \leq \low$ options for $z$.
        
        We have shown that $L$ is an $(\high,\low)$-list assignment for $H$. Since $\asi(H) = \asi(G) < \infty$,  Theorem~\ref{theo:BorelReedSudakov_asi} yields a Borel proper $L$-coloring $f \colon A \to \N$ of $H$. Define
        \[
            M \,\defeq\, \big\{\set{x,y} \,:\, \text{$x \in A$, $y$ is the unique vertex in $N_G(x) \cap U_{f(x)}$}\big \}.
        \]
        The fact that $f$ is proper implies that $M$ is a matching in $G$. Since $M$ is Borel and covers $A$ by construction, this completes the proof of the $k = 1$ case.

        Now suppose $k > 1$. Our plan is to reduce this situation to the $k=1$ case by passing to a suitable simple spanning subgraph of $G$.
        
        We may assume that $\epsilon < 1$ and $\high < 2\low$. Set $\delta \defeq \epsilon/100$ and $N \defeq \lceil 2k^3/\delta \rceil$. 
            %
        Define
        \[
            \high' \,\defeq\, \left\lceil (1-2\delta)\frac{\high}{N} \right\rceil \qquad \text{and} \qquad \low' \,\defeq\, \left\lfloor (1+\delta)\frac{\low}{N}\right\rfloor.
        \]

        \begin{claim}\label{claim:simplify}
            Assuming $\low$ is large enough as a function of $\epsilon$ and $k$, there exists a simple spanning Borel subgraph $G' \subseteq G$ such that, in the graph $G'$, the bipartition $(A,B)$ is of type $(\high',\low')$.
        \end{claim}
        \begin{claimproof}
            By Lemma~\ref{lemma:reg_subgraph}, we may replace $G$ by a spanning subgraph to arrange that every vertex in $A$ has degree exactly $\high$. 
            Given a partial function $f \colon E \pto N$, 
            we say that an edge $e \in E$ 
            is \emphd{duplicated} in $f$ if $f(e) = 0$ and there exists some other edge $e' \in f^{-1}(0) \setminus \set{e}$ such that $\pi(e') = \pi(e)$. Suppose $f \colon E \to N$ is a Borel function with the following properties:
            \begin{itemize}
                \item every vertex $x \in A$ is incident to at least $(1-\delta)\high/N$ edges in $f^{-1}(0)$,
                \item every vertex $y \in B$ is incident to at most $(1+\delta)\low/N$ edges in $f^{-1}(0)$,
                \item every vertex $x \in A$ is incident to at most $\delta \high/N$ edges that are duplicated in $f$.
            \end{itemize}
            Then, letting $G'$ be the spanning subgraph of $G$ with edge set $\set{e \in f^{-1}(0)  \,:\, e \text{ is not duplicated in } f}$ completes the proof of the claim. It remains to show that such a function $f$ exists.

            The problem of finding the desired function $f$ can be represented as a CSP $\mathcal{B} \colon E \to^? N$. For each vertex $x \in V$, we let $E_x \subseteq E$ be the set of all edges incident to $x$. For $x \in A$, we define two constraints, $B_x$ and $B_x'$, with domain $E_x$ as follows:
            \begin{align*}
                B_x \,&\defeq\, \big\{\,\phi \colon E_x \to N \ :\  |\phi^{-1}(0)| < (1-\delta)\high/N\,\big\}, \\
                B_x' \,&\defeq\, \big\{\,\phi \colon E_x \to N \ :\  \text{more than $\delta\high/N$ edges in $E_x$ are duplicated in $\phi$}\,\big\}.
            \end{align*}
            Similarly, for $y \in B$, we let $B_y$ be the constraint with domain $E_y$ given by
            \begin{align*}
                B_y \,&\defeq\, \big\{\,\psi \colon E_y \to N \ :\  |\psi^{-1}(0)| > (1+\delta)\low/N\,\big\}.
            \end{align*}
            Our goal is to find a Borel solution $f \colon E \to N$ to the CSP \[\mathcal{B} \,\defeq\, \set{B_x,\,B_x',\,B_y \,:\, x \in A, \,y \in B}.\] It is routine to check that the CSP $\mathcal{B}$ is Borel and that $\asi(G_\mathcal{B}) \leq \asi(G) < \infty$.
            
            To apply Theorem~\ref{theo:BLLL}, we now need to bound $\mathsf{p}(\mathcal{B})$ and $\mathsf{d}(\mathcal{B})$. For $\mathsf{d}(\mathcal{B})$, it is easy to see that
            \[
                \mathsf{d}(\mathcal{B}) \,\leq\, \max \set{1 + \high,\, 2\low} \,\leq\, 2\low,
            \]
            where we use the assumption that $\high < 2 \low$. Now we turn to $\mathsf{p}(\mathcal{B})$. Take $x \in A$ and suppose a mapping $\phi \colon E_x \to N$ is chosen uniformly at random. Then $|\phi^{-1}(0)|$ is a binomial random variable with distribution $\mathsf{Bin}(\high, 1/N)$. Thus, by the Chernoff bound \cite[Thm.~A.1.4]{AS},
            \begin{equation}\label{eq:prob_simplify_1}
                \P[B_x] \,=\, \P\left[\,|\phi^{-1}(0)| < (1-\delta)\high/N\,\right] \,\leq\, \exp\left(-\frac{2\delta^2}{N^2}\,\high\right).
            \end{equation}
            Similarly, for all $y \in B$, we have
            \begin{equation}\label{eq:prob_simplify_2}
                \P[B_y] \,\leq\, \exp\left(-\frac{2\delta^2}{N^2}\,\low\right).
            \end{equation}
            Next we consider the constraints $B_x'$ for $x \in A$. Suppose $x$ has $n$ distinct neighbors $y_1$, \ldots, $y_n \in B$. Pick $\phi \colon E_x \to N$ uniformly at random and, for each $1 \leq i \leq n$, let $D_i$ be the random event that at least two edges in the set $\phi^{-1}(0)$ join $x$ to $y_i$. There are at most $k$ edges joining $x$ to $y_i$, so we have
            \[
                \P[D_i] \,\leq\, \frac{1}{N^2} \,{k \choose 2} \,\leq\, \frac{k^2}{N^2} \,\leq\, \frac{\delta}{2kN},
            \]
            where the last step holds by the choice of $N$. 
            Since the events $D_i$ are mutually independent and $n \leq \high$, the number of indices $1 \leq i \leq n$ such that $D_i$ holds is stochastically dominated by $\mathsf{Bin}(\high,\delta/(2kN))$. Therefore, using the $k$-simplicity of $G$ and the Chernoff bound, we obtain
            \begin{equation}\label{eq:prob_simplify_3}
                \P[B_x'] \,\leq\, \P[\text{$D_i$ holds for at least $\delta\high/(kN)$ indices $i$}] \,\leq\, \exp\left(-\frac{\delta^2}{2k^2N^2} \,\high\right).
            \end{equation}
            Comparing the bounds \eqref{eq:prob_simplify_1}, \eqref{eq:prob_simplify_2}, and \eqref{eq:prob_simplify_3}, we conclude that
            \[
                \mathsf{p}(\mathcal{B}) \,\leq\, \exp\left(-\frac{\delta^2}{2k^2N^2} \,\low\right).
            \]

            Finally, we combine our bounds on $\mathsf{d}(\mathcal{B})$ and $\mathsf{p}(\mathcal{B})$ and observe that
            \[
                2^{15} \cdot \mathsf{p}(\mathcal{B}) \cdot (\mathsf{d}(\mathcal{B}) + 1)^8 \,\leq\, 2^{15} \cdot \exp\left(-\frac{\delta^2}{2k^2N^2} \,\low\right) \cdot (2\low + 1)^{8} \,<\, 1,
            \]
            provided $\low$ is large enough. Theorem~\ref{theo:BLLL} then yields a Borel solution to $\mathcal{B}$, as desired.
        \end{claimproof}

        Let $G'$ be the graph given by Claim~\ref{claim:simplify}. Note that, assuming $\low$ is large enough as a function of $\epsilon$ and $k$, we have $\low' \geq b_0(\epsilon/2,1)$. Furthermore, our choice of $\delta$ implies that 
        \[
            \frac{\high'}{\low'} \,\geq\, \frac{(1-2\delta)\high}{(1+\delta) \low} \,\geq\, \frac{(1-2\delta) (1+\epsilon)}{1+\delta} \,\geq\, 1 + \frac{\epsilon}{2}. 
        \]
        Hence, by the simple graph case of Theorem~\ref{theo:asi_simple}, the graph $G'$ has a Borel matching $M$ covering $A$. Since $G' \subseteq G$, $M$ is also a matching in $G$, and the proof is complete.
    \end{scproof}

    \subsection{Bounded distance matchings: Theorem~\ref{theo:asi_bounded_distance} and related results}\label{subsec:bd}

    Here we apply Theorem~\ref{theo:asi_simple} to derive Theorem~\ref{theo:asi_bounded_distance}. We will need the following lemma:

    \begin{lemma}\label{lemma:lumping}
        Let $G = (V,E,\pi)$ be a locally finite Borel graph all whose components are infinite. Then, for all $n \in \N$, there exists a Borel bijection $F \colon V \to V$ such that for all $x \in V$:
        \begin{itemize}
            \item $\dist_G(x,F(x)) \leq 4n$, and
            \item the points $x$, $F(x)$, \ldots, $F^{n-1}(x)$ are all distinct.
        \end{itemize}
    \end{lemma}
    \begin{scproof}
        Using Corollary~\ref{corl:spaced}\ref{item:max_spaced}, we pick a Borel inclusion-maximal $2n$-independent set $C \subseteq V$. Then, for every vertex $x \in V$, we have $\dist_G(x,C) \leq 2n$. Let $c_x \in C$ be a vertex that minimizes the distance to $x$; if there are several such vertices, we use Luzin--Novikov Uniformization (Lemma~\ref{lemma:LNU}) to pick one of them so that the map $V \to C \colon x \to c_x$ is Borel. For each $c \in C$, let $V_c \defeq \set{x \in V \,:\, c_x = c}$. By construction, the sets $(V_c \,:\, c \in C)$ form a partition of $V$. Note also that for all $c \in C$,
        \[
            B_G(c,n) \,\subseteq\, V_c \,\subseteq\, B_G(c,2n).
        \]
        The second inclusion holds since $\dist_G(x,c_x) \leq 2n$ for all $x \in V$ by construction. The first inclusion is true because for each $x \in B_G(c,n)$ and $c ' \in C \setminus \set{c}$, we have
        \[
            \dist_G(x,c') \,\geq\, \dist_G(c,c') \,-\, \dist_G(x,c) \,>\, n \,\geq\, \dist_G(x,c),
        \]
        which implies that $c_x = c$. Since every component of $G$ is infinite, every $n$-ball in $G$ contains more than $n$ vertices. Therefore, each $c \in C$ satisfies
        \[
            n \,<\, |B_G(c,n)| \,\leq\, |V_c| \,\leq\, |B_G(c,2n)| \,<\,\infty.
        \]
        Let $F_c \colon V_c \to V_c$ be a permutation of the finite set $V_c$ with a single cycle. By Luzin--Novikov Uniformization again, out of finitely many such permutations, we can choose one in a Borel way, so that the map $F \colon V \to V$ given by $F(x) \defeq F_{c_x}(x)$ is Borel. By construction, $F$ is a bijection and \[\dist_G(x,F(x)) \,\leq\, \dist_G(x,c_x) \,+\, \dist_G(c_x, F(x)) \,\leq\, 4n.\] 
        Finally, since $|V_c| > n$ for all $c \in C$, the points $x$, $F(x)$, \ldots, $F^{n-1}(x)$ are distinct for each $x \in V$. 
    \end{scproof}

    Now we are ready for the proof of Theorem~\ref{theo:asi_bounded_distance}, restated here for ease of reference:

    \begin{theocopy}{theo:asi_bounded_distance}
        Suppose that a Borel bipartition $(A,B)$ of a locally finite Borel graph $G = (V,E,\pi)$ has combinatorial expansion. If $\asi(G) < \infty$, then there exist $d \in \N$ and a Borel injection $f \colon A \to B$ such that, for all $x \in A$, $\dist_G(x,f(x)) \leq d$.
    \end{theocopy}
    \begin{scproof}
        As in the proof of Theorem~\ref{theo:asi_two_conjectures} given in \S\ref{subsec:reduction}, Lemma~\ref{lemma:make_unbalanced} shows that it is enough to consider the case when the bipartition $(A,B)$ is unbalanced.  
        Hence, we will assume from now on that the bipartition $(A,B)$ is of type $(\high,\low)$, where $\high > \low$.
        
        Let $V_{<\infty} \subseteq V$ be the union of all finite components of $G$, and let $V_\infty \defeq V \setminus V_{<\infty}$. Since the sets $V_{<\infty}$ and $V_\infty$ are Borel and $G$-invariant, we may work in each of them separately. The graph $G[V_{<\infty}]$ is component-finite, so it has a Borel matching $M$ that covers $A \cap V_{<\infty}$ by Theorem~\ref{theo:smooth}. Then, for each $x \in V_{<\infty}$, we can let $f(x)$ be the neighbor to which $x$ is matched by $M$. Thus, it only remains to prove the theorem for the graph $G[V_\infty]$, so we will assume from now on that $V = V_\infty$, i.e., every component of $G$ is infinite.

        Pick $\epsilon > 0$ so that $\high > (1+\epsilon)\, \low$. Let $b_0(\epsilon,k)$ be as in Theorem~\ref{theo:asi_simple} and define $n \defeq b_0\left(\epsilon, \low\right)$.

        \begin{claim}\label{claim:F_defined}
            There exists a Borel bijection $F \colon A \to A$ such that for all $x \in A$,
        \begin{itemize}
            \item $\dist_G(x,F(x)) \leq 8n$, and
            \item the points $x$, $F(x)$, \ldots, $F^{n-1}(x)$ are all distinct.
        \end{itemize}
        \end{claim}
        \begin{claimproof}
            Let $H$ be the square of the underlying simple graph of $G$. Note that every component of the induced subgraph $H[A]$ is infinite, because if $C$ is a component of $H[A]$, then the vertices in $C$ or adjacent to $C$ in $G$ form a component of $G$, which must be infinite by assumption. Therefore, we can obtain the desired function $F$ by applying Lemma~\ref{lemma:lumping} to the graph $H[A]$. 
        \end{claimproof}

        Fix a function $F$ given by Claim~\ref{claim:F_defined} and define a graph $G' = (V, E', \pi')$ as follows: Let
        \[
            E' \,\defeq\, E \times \set{0,1,\ldots, n-1}.
        \]
        For $0 \leq i < n$ and an edge $e \in E$ with $\pi(e) = \set{x,y}$, where $x \in A$ and $y \in B$, we let
        \[
            \pi'(e,i) \,\defeq\, \set{F^i(x), y}.
        \]
        In other words, we put in $G'$ edges joining each of the vertices $x$, $F(x)$, \ldots, $F^{n-1}(x)$ to $y$. In particular, if $x \in A$ and $y \in B$ are adjacent in $G'$, then
        \[
            \dist_G(x,y) \,\leq\, 8n(n-1) + 1.
        \]
        It follows that $\asi(G') \leq \asi(G) < \infty$, since the underlying simple graph of $G'$ is contained in a finite power of the underlying simple graph of $G$. 
        
        It is clear from the construction that for all $y \in B$,
        \[
            \deg_{G'}(y) \,=\, n \deg_G(y) \,\leq\, n \low.
        \]
        On the other hand, for all $x \in A$,
        \[
            \deg_{G'}(x) \,=\, \sum_{i = 0}^{n-1} \deg_G(F^{-i}(x)) \,\geq\, n \high.
        \]
        Hence, $(A,B)$ is a bipartition of $G'$ of type $(n\high, n\low)$. Note also that the graph $G'$ is $\low$-simple, because for any pair of vertices $x \in A$, $y \in B$, the number of edges joining $y$ to $x$ in $G'$ is equal to the number of edges joining $y$ to the set $\set{F^{-i}(x) \,:\, 0 \leq i < n}$ in $G$, which cannot exceed $\deg_G(y) \leq \low$. Since
        \[
            n\high \,>\, (1+\epsilon)\, n\low \qquad \text{and} \qquad n \low \,\geq\, n \,=\, b_0(\epsilon, b),
        \]
        we can use Theorem~\ref{theo:asi_simple} to find a Borel matching $M$ in $G'$ that covers $A$. Now, for each $x \in A$, we let $f(x)$ be the $G'$-neighbor to which $x$ is matched by $M$, and observe that the resulting function $f$ satisfies the conclusion of Theorem~\ref{theo:asi_bounded_distance} with $d \defeq 8n(n-1) + 1$.
    \end{scproof}

    As mentioned in \S\ref{subsec:intro_Borel}, the conclusion of Theorem~\ref{theo:asi_bounded_distance} no longer holds if we drop the assumption that $\asi(G) < \infty$. Indeed, we have the following general fact:

    \begin{prop}\label{prop:no_dist_matching}
        Let $G = (V,E,\pi)$ be a locally finite Borel graph with a Borel bipartition $(A,B)$ such that $\deg_G(y) \leq 2$ for all $y \in B$. The following statements are equivalent:
        \begin{enumerate}[label=\ep{\normalfont\roman*}]
            \item\label{item:actual_matching} $G$ has a Borel matching that covers $A$,
            \item\label{item:distance_matching} there exist $d \in \N$ and a Borel injection $f \colon A \to B$ such that $\dist_G(x,f(x)) \leq d$ for all $x \in A$.
        \end{enumerate}
    \end{prop}

    By Theorem~\ref{theo:Marks}, there exist graphs $G$ for which condition \ref{item:actual_matching} in Proposition~\ref{prop:no_dist_matching} fails, and hence so does condition~\ref{item:distance_matching}. 
        
    \begin{scproof}
        Implication \ref{item:actual_matching} $\Longrightarrow$ \ref{item:distance_matching} is clear (if \ref{item:actual_matching} holds, then \ref{item:distance_matching} is satisfied with $d = 1$). Now suppose \ref{item:distance_matching} holds, i.e., for some $d \in \N$, there is a Borel injection $f \colon A \to B$ such that $\dist_G(x,f(x)) \leq d$ for all $x \in A$. Our goal is to find a Borel matching in $G$ that covers $A$. The approach we follow is similar to the one used in \cite[proof of Thm.~1.7]{HypBorCom}.

        For each $x \in A$, let $d(x) \defeq \dist_G(x,f(x)) \leq d$, and let
        \[
            \mathsf{path}(x) \,=\, \left(e_1(x), \, e_2(x), \,\ldots,\, e_{d(x)}(x)\right) \,\in\, E^{d(x)}
        \]
        be a sequence of edges on some shortest path in $G$ from $x$ to $f(x)$. Using Luzin--Novikov Uniformization, we can arrange that the assignment $x \mapsto \mathsf{path}(x)$ is Borel.
        
        Define a function $w \colon E \to \Z$ as follows:
        \[
            w(e) \,\defeq\, |\set{x \in A \,:\, e = e_i(x) \text{ for some odd } i}| \,-\, |\set{x \in A \,:\, e = e_i(x) \text{ for some even } i}|.
        \]
        Note that $w(e)$ is indeed a well-defined finite integer, since the only vertices $x \in A$ that may affect $w(e)$ are the ones within distance $d$ from the endpoints of $e$, and there are only finitely many of them. A crucial observation now is that for all $x \in A$,
        \begin{equation}\label{eq:sum_at_x}
            \sum_{\text{$e \in E$ is incident to $x$}} w(e) \,=\, 1,
        \end{equation}
        because for each vertex $x' \in A \setminus \set{x}$, $\mathsf{path}(x')$ contributes $0$ to this sum, while $\mathsf{path}(x)$ contributes $1$. Similar reasoning shows that for all $y \in B$,
        \begin{equation}\label{eq:sum_at_y}
            \sum_{\text{$e \in E$ is incident to $y$}} w(e) \,\in\, \set{0,1};
        \end{equation}
        the actual value is either $1$ or $0$ depending on whether $y$ is in $f(A)$.

        Let $P \defeq \set{e \in E \,:\, w(e) > 0}$. By \eqref{eq:sum_at_x}, every vertex in $A$ is incident to at least one edge in $P$. On the other hand, a vertex $y \in B$ cannot be incident to more than one edge in $P$, because $\deg_G(y) \leq 2$, and if $y$ is incident to $2$ edges in $P$, then the sum in \eqref{eq:sum_at_y} would be at least $2$. By Luzin--Novikov Uniformization, we can pass to a further Borel subset $M \subseteq P$ such that every vertex in $A$ is incident to exactly one edge in $M$. Then $M$ is a Borel matching in $G$ covering $A$, as desired. 
    \end{scproof}

    \begin{remark}
        The function $w$ constructed in the proof of Proposition~\ref{prop:no_dist_matching} above is essentially a \emph{flow} in $G$ that removes one unit of ``charge'' from each vertex in $A$ and adds at most one unit to each vertex in $B$. Such objects are frequently used in the study of matchings in Borel graphs and definable equidecompositions; see, e.g., \cite{marks-unger,ciesla.sabok,bowen2025uniform,bowen.kun.sabok,bowen2024measurable,mathe2022circle,kun2021measurable}.
    \end{remark}


    
    Using Proposition~\ref{prop:no_dist_matching} and Theorem~\ref{theo:asi_bounded_distance}, it is easy to derive Theorem~\ref{theo:asi_32}:

    \begin{theocopy}{theo:asi_32}
        Let $(A,B)$ be a Borel bipartition of a locally finite Borel graph $G$ with combinatorial expansion. If $\asi(G) < \infty$ and $\deg_G(y) \leq 2$ for all $y \in B$, then $G$ has a Borel matching covering $A$.
    \end{theocopy}
    \begin{scproof}
        By Theorem~\ref{theo:asi_bounded_distance}, we can find $d \in \N$ and a Borel injection $f \colon A \to B$ with $\dist_G(x,f(x)) \leq d$ for all $x \in A$. By Proposition~\ref{prop:no_dist_matching}, this implies that $G$ has a Borel matching covering $A$. 
    \end{scproof}



    We can also give a different, somewhat more direct proof of Theorem~\ref{theo:asi_32} in the special case when every vertex in $A$ has degree at least $3$.\footnote{It is not hard to reduce the full statement of Theorem~\ref{theo:asi_32} to this case by observing that every vertex in $A$ must have degree at least $2$, and paths of degree-$2$ vertices in $G$ must have bounded length. We omit the details, since we have a different proof of Theorem~\ref{theo:asi_32} already.} It uses the following Borel version of Brooks' theorem \cites{Brooks}[Thm.~5.2.4]{Diestel} (see \S\ref{subsec:chi} for the definition of $\chi_\mathsf{B}(G)$):

    \begin{theorem}[{AB--FW \cite[Thm.~1.17]{bw_lll}}]\label{theo:Brooks_asi}
        Let $G$ be a Borel graph with finite maximum degree $\Delta \geq 3$ and without a $(\Delta + 1)$-clique. If $\asi(G) < \infty$, then $\chi_\mathsf{B}(G) \leq \Delta$.
    \end{theorem}

    \begin{scproof}[ of Theorem~\ref{theo:asi_32} for bipartitions of type $(3,2)$ \ep{alternative}]
        Let $G = (V,E,\pi)$ be a locally finite Borel graph with $\asi(G) < \infty$ and a bipartition $(A,B)$ of type $(3,2)$. 
        Using Lemma~\ref{lemma:reg_subgraph}, we may arrange that every vertex in $A$ has degree exactly $3$. 
        Let $H$ be the simple graph with vertex set $E$ such that $\set{e,e'} \in [E]^2$ is an edge of $H$ if and only if $e$ and $e'$ have a common endpoint. Clearly, the maximum degree of $H$ is at most $2 + 1 = 3$. Furthermore, it is easy to see that $H$ has no $4$-cliques. Therefore, by Theorem~\ref{theo:Brooks_asi}, there exists a partition $E = M_0 \sqcup M_1 \sqcup M_2$ into $3$ Borel sets that are independent in $H$; in other words, each $M_i$ is a Borel matching in $G$. Since every vertex in $A$ is incident to $3$ edges, it must be incident to one edge in each of $M_0$, $M_1$, $M_2$. Hence, every one of these matchings covers $A$, as desired.
    \end{scproof}

    Although the above argument is quite simple, it is worth noting that the proof of Theorem~\ref{theo:Brooks_asi} itself uses distributed algorithms and ultimately relies on Theorem~\ref{theo:BLLL}. 

    \section{Applications}\label{sec:applications}

    \subsection{Remarks on quasi-invariance and related issues}\label{subsubsec:qi}



    A useful simplifying assumption in measurable combinatorics is quasi-invariance of the measure:

    \begin{definition}[Quasi-invariant measures]
        Let $G = (V,E,\pi)$ be a locally finite Borel graph. For a set $X \subseteq V$, we let $[X]_G$ be the \emphd{$G$-saturation} of $X$, i.e., the set of all vertices reachable from $X$ by a finite path in $G$. Equivalently, $[X]_G$ is the smallest $G$-invariant set containing $X$. A probability measure $\mu$ on $V$ is \emphd{$G$-quasi-invariant} if for every $\mu$-null set $X \subseteq V$, the set $[X]_G$ is also $\mu$-null. 
    \end{definition}

    Most of the time, assuming that the measure under consideration is quasi-invariant results in no loss of generality, thanks to the following easy and well-known fact:

    \begin{prop}[{\cite[Prop.~4.16]{kechrisCBER}}]\label{prop:quasi-invariant}
        Let $G$ be a locally finite Borel graph with a probability measure $\mu$ on its vertex set. Then there exists a $G$-quasi-invariant probability measure $\mu^*$ such that $\mu$ is absolutely continuous with respect to $\mu^*$ \ep{meaning that every $\mu^*$-null set is also $\mu$-null}.
    \end{prop}


    Combining Theorem~\ref{theo:measurable} with Proposition~\ref{prop:quasi-invariant} yields the following:

    \begin{cor}\label{corl:invariant}
        Let $(A,B)$ be an unbalanced bipartition of a 
        locally finite Borel graph $G = (V,E,\pi)$. Then, for any probability measure $\mu$ on $V$, there exist a $\mu$-null $G$-invariant Borel set $X \subseteq V$ and a Borel matching $M$ in $G$ that covers every vertex in $A \setminus X$.
    \end{cor}
    \begin{scproof}
        Let $\mu^*$ be a $G$-quasi-invariant probability measure such that $\mu$ is absolutely continuous with respect to $\mu^*$. Applying Theorem~\ref{theo:measurable} with $\mu^*$ in place of $\mu$ gives a Borel matching $M$ that covers $\mu^*$-almost every vertex in $A$. Now we let $Y \subseteq A$ be the set of all vertices in $A$ that are not covered by $M$ and take $X \defeq [Y]_G$.
    \end{scproof}


    An analogous strengthening of Theorem~\ref{theo:BM} was proved by Marks and Unger in \cite{marks2016baire} (modulo the small changes described in \S\ref{subsec:Baire}):

    \begin{theorem}[{ess.~Marks--Unger \cite[Thm.~1.3]{marks2016baire}}]\label{theo:BM_invariant}
        Suppose that a Borel bipartition $(A,B)$ of a locally finite Borel graph $G = (V,E,\pi)$ has combinatorial expansion. Then, for any compatible topology $\tau$ on $V$, there exist a $\tau$-meager $G$-invariant Borel set $X \subseteq V$ and a Borel matching $M$ in $G$ that covers every vertex in $A \setminus X$.
    \end{theorem}

    We remark that in a forthcoming paper \cite{ChenWeilacher}, Chen and the third named author 
    establish a version of Proposition~\ref{prop:quasi-invariant} for Baire category, which shows that Theorem~\ref{theo:BM_invariant} is really a consequence of Theorem~\ref{theo:BM}, in the same way as Corollary~\ref{corl:invariant} is a consequence of Theorem~\ref{theo:measurable}.

    Of the applications presented in this paper, Corollary~\ref{corl:invariant} and Theorem~\ref{theo:BM_invariant} will only be explicitly invoked in the proofs of Theorem~\ref{theo:balance} (almost balanced orientations) and Theorem~\ref{theo:unperf} (almost unperforation of type semigroups). 
    However, they can 
    be of use in other situations as well, and in
    general they often make applying Theorems~\ref{theo:measurable} and \ref{theo:BM}  
    more straightforward (see, e.g., Remark~\ref{remk:invariance}).

    \subsection{Shannon's bound for Borel multigraphs}\label{sec:shannon}

    \subsubsection{Proofs of Theorems~\ref{theo:Shannon_meas} and \ref{theo:Shannon_Borel}}\label{subsubsec:Shannon1}


    Here we prove Theorems~\ref{theo:Shannon_meas} and \ref{theo:Shannon_Borel}, restated in a combined form for the reader's convenience:
    
    \begin{theorem}[Theorems~\ref{theo:Shannon_meas} and \ref{theo:Shannon_Borel}]\label{theo:shannon_combined}
        Let $G = (V,E,\pi)$ be a Borel graph of maximum degree $\Delta \in \N$ and let $q \defeq \lfloor 3\Delta/2\rfloor$. 
        \begin{enumerate}[label=\ep{\normalfont\arabic*}]
            \item\label{item:shannon_measure} If $\mu$ is a probability measure on $V$, then $\chi'_\mu(G) \leq q$.

            \item If $\tau$ is a compatible topology on $V$, then $\chi'_\tau(G) \leq q$.

            \item If $G$ is of subexponential growth or everywhere $2$-ended, then $\chi'_\mathsf{B}(G) \leq q$.
        \end{enumerate}
    \end{theorem}

    In addition to the results concerning Borel matchings in unbalanced bipartite graphs, we will also rely on descriptive set-theoretic versions of Brooks' theorem from graph theory \cites{Brooks}[Thm.~5.2.4]{Diestel}, summarized in the following statement:

    \begin{theorem}[{Conley--Marks--Tucker-Drob \cite{MeasurableBrooks}, AB \cite[\S3.1]{BernshteynDistributed}, Greb\'ik--Vidny\'anszky \cite[Thm.~8.4]{GV_from}, AB--FW \cite[Thm.~1.17]{bw_lll}}]\label{theo:Brooks_meas}
        Let $\Delta \geq 3$ be an integer and let $G = (V,E,\pi)$ be a Borel graph of maximum degree at most $\Delta$ and without a $(\Delta + 1)$-clique. Then there exists a Borel $G$-invariant set $X \subseteq V$ such that $\chi_\mathsf{B}(G - X) \leq \Delta$, and furthermore:
        \begin{enumerate}[label=\ep{\normalfont\arabic*}]
            \item\label{item:Brooks_meas} if $\mu$ is a probability measure on $V$, we can arrange that $X$ is $\mu$-null,

            \item\label{item:Brooks_cat} if $\tau$ is a compatible topology on $V$, we can arrange that $X$ is $\tau$-meager,

            \item\label{item:Brooks_subexp} if $G$ is of subexponential growth, we can arrange that $X = \0$,

            \item\label{item:Brooks_asi} if $\asi(G) < \infty$, we can arrange that $X = \0$.
        \end{enumerate}
    \end{theorem}

    Parts~\ref{item:Brooks_meas} and \ref{item:Brooks_cat} of Theorem~\ref{theo:Brooks_meas} were proved by Conley, Marks, and Tucker-Drob in \cite{MeasurableBrooks}. Part~\ref{item:Brooks_subexp} also follows from their general results, although it is not stated in \cite{MeasurableBrooks} explicitly. An explicit statement, with a different proof, appears in \cite[\S3.1]{BernshteynDistributed}, and yet another proof was given by Greb\'ik and Vidny\'anszky \cite[Thm.~8.4]{GV_from}. Part~\ref{item:Brooks_asi} is a restatement of Theorem~\ref{theo:Brooks_asi}, which was established in \cite{bw_lll}. As far as Theorem~\ref{theo:shannon_combined} is concerned, Theorem~\ref{theo:Brooks_meas}\ref{item:Brooks_asi} is relevant because it includes the case when $G$ is a subgraph of an everywhere $2$-ended locally finite Borel graph. (See \S\ref{subsec:asi_defns} for a discussion of asymptotic separation index.)
    
    To prove Theorem~\ref{theo:shannon_combined}, we need the following consequence of Theorem~\ref{theo:Brooks_meas}:

    \begin{cor}\label{corl:3}
        Let $G = (V,E,\pi)$ be a Borel graph with maximum degree at most $3$. Suppose no two vertices of degree $3$ in $G$ are adjacent.
        \begin{enumerate}[label=\ep{\normalfont\arabic*}]
            \item\label{item:3_meas} If $\mu$ is a probability measure on $V$, then $\chi_\mu'(G) \leq 3$.

            \item\label{item:3_cat} If $\tau$ is a compatible topology on $V$, then $\chi_\tau'(G) \leq 3$.

            \item\label{item:3_subexp} If $G$ is of subexponential growth, then $\chi_\mathsf{B}'(G) \leq 3$.

            \item\label{item:3_asi} If $\asi(G) < \infty$, then $\chi_\mathsf{B}'(G) \leq 3$.
        \end{enumerate}
    \end{cor}

    A special case of Corollary~\ref{corl:3}\ref{item:3_asi} was already used in the alternative proof of Theorem~\ref{theo:asi_32} for bipartitions of type $(3,2)$  in \S\ref{subsec:bd}.
    
    \begin{scproof}
        We may assume that every vertex of $G$ is incident to at least one edge. Let $H$ be the simple graph with vertex set $E$ such that $\set{e,e'} \in [E]^2$ is an edge in $H$ if and only if $e$ and $e'$ have a common endpoint. It is easy to check that the maximum degree of $H$ is at most $3$ and there are no $4$-cliques in $H$. The desired result now follows by applying Theorem~\ref{theo:Brooks_meas} to $H$. For example, to prove \ref{item:3_meas}, we use Luzin--Novikov Uniformization (Lemma~\ref{lemma:LNU}) to find a Borel map $\tau \colon V \to E$ such that $\tau(x)$ is incident to $x$ for all $x \in V$. Let $\nu \defeq \tau_\ast(\mu)$ be the pushforward of $\mu$ by $\tau$. By Theorem~\ref{theo:Brooks_meas}\ref{item:Brooks_meas} applied to the graph $H$ and the measure $\nu$, there exists a $\nu$-null Borel $H$-invariant set $F \subseteq E$ such that $\chi_\mathsf{B}(H - F) \leq 3$. Then the set $X \subseteq V$ of all vertices incident to an edge in $F$ is $G$-invariant and $\mu$-null, and the graph $G - X$ has a proper Borel $3$-edge-coloring $f$. Extending $f$ to the edges in $G[X]$ arbitrarily \ep{say, by assigning $0$ to all of them} yields a Borel $3$-edge-coloring of $G$ that is proper at $\mu$-almost every vertex, as desired. The other parts of Corollary~\ref{corl:3} are proved in the same way, \emph{mutatis mutandis}.
    \end{scproof}

    We will also require the following simple observation: 

    \begin{lemma}\label{lemma:extend}
        Let $G = (V,E,\pi)$ be a Borel graph of maximum degree $\Delta \in \N$ and fix an integer $q \geq \Delta$. Suppose $X \subseteq V$ is a Borel set and $f'$ is a proper Borel $q$-edge-coloring of the graph $G - X$. Then $f'$ can be extended to a Borel $q$-edge-coloring $f \colon E \to q$ that is proper at every vertex $x \notin X$.
    \end{lemma}

    \begin{remark}\label{remk:invariance}
        By using somewhat stronger versions of Corollary~\ref{corl:invariant} and Theorem~\ref{theo:BM_invariant}, it is possible to ensure that in the proof of Theorem~\ref{theo:shannon_combined}, we only need to apply Lemma~\ref{lemma:extend} in the trivial case when the set $X$ is $G$-invariant (and thus \emph{any} extension of $f'$ has the desired property). The needed strengthening of Corollary~\ref{corl:invariant} can be derived from Proposition~\ref{prop:quasi-invariant}, while for Theorem~\ref{theo:BM_invariant} it follows from the analogous Baire-categorical result in \cite{ChenWeilacher}. Since the paper \cite{ChenWeilacher} has not yet appeared, we give here a proof of Lemma~\ref{lemma:extend} for arbitrary $X$, sidestepping the issue of $G$-invariance. 
    \end{remark}

    \begin{scproof}[ of Lemma~\ref{lemma:extend}]
        The edges in $G[X]$ may be colored arbitrarily (say, we can give all of them the color $0$), so we only need to be careful when extending $f'$ to the edges joining $V \setminus X$ to $X$. For each $x \in V \setminus X$, let $E_0(x)$ be the set of all edges that join $x$ to $X$, and let $E_1(x)$ be the set of all edges that join $x$ to another vertex in $V \setminus X$. Let $C(x) \defeq \set{f'(e) \,:\, e \in E_1(x)}$ and note that
        \begin{equation}\label{eq:enough_colors}
            |E_0(x)| + |C(x)| \,=\, |E_0(x)| + |E_1(x)| \,=\, \deg_G(x) \,\leq\, \Delta \,\leq\, q.
        \end{equation}
        Fix an arbitrary Borel linear ordering $\leq$ on $E$ and list the elements of $E_0(x)$ in the increasing order as $e(1,x) < e(2,x) < \cdots$. Thanks to \eqref{eq:enough_colors}, we can then let $f(e(i,x))$ be the $i$-th smallest element in the set $\set{0,1,\ldots,q-1} \setminus C(x)$. The Luzin--Novikov theorem (see \S\ref{subsec:LN}) ensures that the resulting extension of $f'$ is Borel, and, by construction, it is proper at all vertices in $V \setminus X$, as desired.
    \end{scproof}

    Now we are ready to start the proof of Theorem~\ref{theo:shannon_combined}. For concreteness, we will only present the proof of Theorem~\ref{theo:shannon_combined}\ref{item:shannon_measure}, i.e., the measure-theoretic part of the result. The other statements in Theorem~\ref{theo:shannon_combined} follow via exactly the same argument, except that Theorem~\ref{theo:measurable} should be replaced by Theorems~\ref{theo:BM},~\ref{theo:growth}, or~\ref{theo:ends1}. 
    We employ an inductive approach derived from the work of Ghaffari, Kuhn, Maus, and Uitto  \cite{ghaffari2018deterministic}. Specifically, we rely on the following key lemma:

    \begin{lemma}\label{lemma:3match}
        Let $G = (V,E,\pi)$ be a Borel graph of finite maximum degree $\Delta \geq 2$ and let $\mu$ be a probability measure on $V$. Then there exist a $\mu$-null Borel set $X \subseteq V$ and pairwise disjoint Borel matchings $M_0$, $M_1$, $M_2$ in $G$ such that $\Delta(G - X - M_0 - M_1 - M_2) \leq \Delta -2$. 
    \end{lemma}

    Recall that $G - X - M_0 - M_1 - M_{2}$ here denotes the subgraph obtained from $G$ by deleting the vertices in $X$ together with their incident edges as well as the edges in the matchings $M_0$, $M_1$, $M_{2}$, and $\Delta(G - X - M_0 - M_1 - M_{2})$ is the maximum degree of this graph. 
    Before proving Lemma~\ref{lemma:3match}, let us explain how it implies Theorem~\ref{theo:shannon_combined}\ref{item:shannon_measure}:

    \begin{scproof}[ of Theorem~\ref{theo:shannon_combined}\ref{item:shannon_measure}]
        Thanks to Lemma~\ref{lemma:extend}, it is enough to find a $\mu$-null Borel set $X \subseteq V$ such that $\chi'_\mathsf{B}(G-X) \leq q$. To this end, we proceed by induction on $\Delta$, with the base cases $\Delta = 0$ and $\Delta = 1$ being trivial (when $\Delta \leq 1$, we clearly have $\chi'_\mathsf{B}(G) \leq \Delta$). Suppose $\Delta \geq 2$ and let $X$, $M_0$, $M_1$, $M_2$ be given by Lemma~\ref{lemma:3match}. Recall that $q = \lfloor 3\Delta/2 \rfloor$ and note that $\lfloor 3(\Delta - 2)/2\rfloor = q - 3$.  
        Hence, applying the inductive hypothesis to the graph \[G' \,\defeq\, G - X - M_0 - M_1 - M_2\] yields a $\mu$-null Borel set $X' \subseteq V\setminus X$ and a proper Borel $(q - 3)$-edge-coloring $f'$ of $G' - X'$. 
        We then extend $f'$ to a proper Borel $q$-edge-coloring of $G - X - X'$ by letting $f'(e) \defeq q-3 + i$ for all $e \in M_i$. This shows that $\chi'_\mathsf{B}(G - X - X') \leq q$, and we are done since the set $X \cup X'$ is $\mu$-null.
    \end{scproof}

    It remains to prove Lemma~\ref{lemma:3match}. The argument below is a Borel adaptation of \cite[\S3]{ghaffari2018deterministic}.

    \begin{scproof}[ of Lemma~\ref{lemma:3match}]
        The core of the proof is in the following claim:

        \begin{claim}\label{claim:3}
            There exist a Borel set $F \subseteq E$ and a $\mu$-null Borel set $U \subseteq V$ such that $\Delta(G - U - F) \leq \Delta-2$ and every vertex of $G$ is incident to at most $3$ edges in $F$.
        \end{claim}
        \begin{claimproof}
            Let
            $
                V_0 \defeq \set{x \in V \,:\, \deg_G(x) = \Delta}$. Using Theorem~\ref{theo:KST}\ref{item:maximal}, we find a Borel inclusion-maximal matching  $F_0$ in the graph $G[V_0]$. 
        Let $G_1 \defeq G - F_0$ and define 
        \[
            V_1 \,\defeq\, \set{x \in V \,:\, \deg_{G_1}(x) = \Delta}\,=\, \set{x \in V_0 \,:\, \text{$x$ is not covered by $F_0$}}.
        \]
        By the maximality of $F_0$, $V_1$ is an independent set in $G$. 
        Therefore, the bipartition $(V_1, V \setminus V_1)$ of the graph $G_1[V_1, V \setminus V_1]$ is unbalanced of type $(\Delta,\Delta-1)$. 
        By Theorem~\ref{theo:measurable}, there exist a Borel matching $F_1$ in $G_1$ and a $\mu$-null Borel set $Y \subseteq V_1$ such that $F_1$ covers $V_1 \setminus Y$. Let $G_2 \defeq G_1 - Y - F_1$. Since every vertex in $V_0 \setminus Y$ is covered by at least one of the matchings $F_0$, $F_1$, we have $\Delta(G_2) \leq \Delta - 1$. Note also that every vertex covered by $F_0$ has degree $\Delta$ in $G$.

        Now we repeat the same construction starting with the graph $G_2$, which yields two disjoint Borel matchings $F_2$, $F_3$ in $G_2$ and a $\mu$-null Borel set $Z \subseteq V \setminus Y$ such that $\Delta(G_2 - Z - F_2 - F_3) \leq \Delta - 2$ and every vertex covered by $F_2$ has degree $\Delta - 1$ in $G_2$.

        We claim that taking $F \defeq F_0 \cup F_1 \cup F_2 \cup F_3$ and $U \defeq Y \cup Z$ works. Since $\Delta(G - U - F) \leq \Delta - 2$ by construction, it only remains to check that every vertex is incident to at most $3$ edges in $F$. That is, we need to rule out the possibility that a vertex $x \in V$ is covered by all $4$ of the matchings $F_0$, $F_1$, $F_2$, $F_3$. But, indeed, if $x$ is covered by both $F_0$ and $F_1$, then $\deg_{G_2}(x) \leq \Delta - 2$, and thus $x$ cannot be covered by $F_2$.
        \end{claimproof}
        
        Let $F$ and $U$ be given by Claim~\ref{claim:3}. Since $U$ is $\mu$-null, we may without loss of generality replace $G$ by $G - U$ and assume $U = \0$. That is, we assume $F \subseteq E$ is a Borel set such that $\Delta(G - F) \leq \Delta - 2$ and every vertex of $G$ is incident to at most $3$ edges in $F$. Let
        \[
            T \,\defeq\, \set{e \in F \,:\, \text{both endpoints of $e$ are incident to $3$ edges in $F$}}.
        \]
        Using Theorem~\ref{theo:KST}\ref{item:maximal} again, we find a Borel inclusion-maximal matching $M \subseteq T$ and let $F' \defeq F \setminus M$.
        
        \begin{claim}\label{claim:still_low_deg}
            $\Delta(G - F') \leq \Delta - 2$.
        \end{claim}
        \begin{claimproof}
            Take any vertex $x \in V$. We are done if $\deg_{G - F'} (x) = \deg_{G-F}(x)$. Otherwise, $x$ is covered by $M$. But then $x$ is incident to $2$ edges in $F'$, and thus $\deg_{G-F'}(x) = \deg_G(x) - 2 \leq \Delta - 2$.
        \end{claimproof}

        Let $H$ be the spanning subgraph of $G$ with edge set $F'$. Then $H$ has maximum degree at most $3$ and, by the maximality of $M$, no two vertices of degree $3$ in $H$ are adjacent. By Corollary~\ref{corl:3}\ref{item:3_meas}, it follows that $\chi'_\mu(H) \leq 3$. In other words, there exists a $\mu$-null Borel set $X \subseteq V$ such that the edge set of $H - X$ can be partitioned into $3$ Borel matchings $M_0$, $M_1$, $M_2$. Since $G - X - M_0 - M_1 - M_2 = G - X - F'$ has maximum degree at most $\Delta - 2$ by Claim~\ref{claim:still_low_deg}, the proof is complete.
    \end{scproof}

    \subsubsection{A partial result using asymptotic separation index: Proof of Theorem~\ref{theo:Shannon_Borel_asi}}

    Now we turn to Theorem~\ref{theo:Shannon_Borel_asi}, restated here for convenience:

    \begin{theocopy}{theo:Shannon_Borel_asi}\label{theo:asi_Shannon_restated}
        Let $G = (V,E,\pi)$ be a Borel graph of maximum degree $\Delta \in \N$ and let $q \defeq \lfloor 3\Delta/2\rfloor$. Suppose that $\asi(G) = s < \infty$. Then $\chi'_\mathsf{B}(G) \leq q + s$.
    \end{theocopy}


    The proof of Theorem~\ref{theo:asi_Shannon_restated} follows the same overall strategy as the proof of Theorem~\ref{theo:shannon_combined} presented in \S\ref{subsubsec:Shannon1}. However, since Conjecture~\ref{conj:asi} about Borel matchings in graphs with finite asymptotic separation index remains open, we cannot implement the argument from \S\ref{subsubsec:Shannon1} directly. Instead we perform an ``approximate'' version of the proof, introducing a small set of ``errors'' in the construction. The error set is responsible for the ``$+s$'' term in the final bound on $\chi'_\mathsf{B}(G)$ in Theorem~\ref{theo:asi_Shannon_restated}. This technique is commonly used in the finite $\asi$ regime, leading to the appearance of similar extra terms in bounds on various other Borel parameters \cite{qian2022descriptive,bowen2023definable,conley2020borel,weilacher.abelian}.



    
    The central role in the proof of Theorem~\ref{theo:asi_Shannon_restated} is played by the machinery of well-spaced shattering systems; the relevant terminology is defined in \S\ref{subsec:reduction}. In place of Lemma~\ref{lemma:3match}, we rely on the following statement:

    \begin{lemma}\label{lemma:3match+error}
        Let $G = (V,E,\pi)$ be a Borel graph of finite maximum degree $\Delta \geq 2$ and with $\asi(G) = s < \infty$. Let $\mathcal{S}_0$ and $\mathcal{S}_1$ be $(s,3)$-shattering systems in $G$ such that $\dist_G(\mathcal{S}_0, \mathcal{S}_1) \geq 3$. Then there exist Borel matchings $M_0$, $M_1$, $M_2$, $N_0$, \ldots, $N_{s-1}$ in $G$ with the following properties:
        \begin{itemize}
            \item $\Delta(G - M_0 - M_1 - M_2 - N_0 - \cdots - N_{s-1}) \leq \Delta -2$, and
            \item for all $0 \leq i < s$, every edge in $N_i$ has at least one endpoint in $\mathcal{S}_0[i] \cup \mathcal{S}_1[i]$.
        \end{itemize}
    \end{lemma}
    \begin{scproof}
        We use the following modification of Claim~\ref{claim:3}:

        \begin{claim}\label{claim:3+error}
            There exist a Borel set $F \subseteq E$ and Borel matchings $N_0$, \ldots, $N_{s-1}$ in $G$ such that:
            \begin{itemize}
                \item $\Delta(G - F - N_0 - \cdots - N_{s-1}) \leq \Delta - 2$,
                \item every vertex of $G$ is incident to at most $3$ edges in $F$, and
                \item for all $0 \leq i < s$, every edge in $N_i$ has at least one endpoint in $\mathcal{S}_0[i] \cup \mathcal{S}_1[i]$.
            \end{itemize}
        \end{claim}
        \begin{claimproof}
        %
        As in the proof of Claim~\ref{claim:3}, we set $V_0 \defeq \set{x \in V \,:\, \deg_G(x) = \Delta}$ and find a Borel maximal matching $F_0$ in $G[V_0]$. Next, we let $G_1 \defeq G - F_0$ and define $V_1 \defeq \set{x \in V \,:\, \deg_{G_1}(x) = \Delta}$. The bipartition $(V_1, V \setminus V_1)$ of the graph $G_1[V_1, V \setminus V_1]$ is unbalanced of type $(\Delta, \Delta-1)$. 
        Since $\mathcal{S}_0$ is an $(s,3)$-shattering system, the set $U \defeq \set{x \in V \,:\, \dist_G(x, \mathcal{S}_0^*) \leq 1}$ 
        is $G$-finite. 
        Hence, by Theorem~\ref{theo:smooth} applied to the graph $G_1[U]$, 
        there exists a Borel matching $F_1$ in $G_1$ that covers $V_1 \cap \mathcal{S}_0^*$. Similarly, there exist Borel matchings 
        %
        $N_{0,0}$, \ldots, $N_{0,s-1}$ in $G_1$ such that 
         $N_{0,i}$ covers 
         $V_1 \cap \mathcal{S}_0[i]$. Additionally, we may 
         choose $N_{0,i}$ so that every edge in $N_{0,i}$  has an endpoint in $\mathcal{S}_0[i]$. 
        
        Let $G_2 \defeq G_1 - F_1 - N_{0,0} - \cdots - N_{0,s-1}$ and note that $\Delta(G_2) \leq \Delta - 1$. Repeating the construction from the previous paragraph with $G_2$ and $\mathcal{S}_1$ in place of $G$ and $\mathcal{S}_0$ produces Borel matchings $F_2$, $F_3$, $N_{1,0}$, \ldots, $N_{1,s-1}$ in $G_2$ such that:
        \begin{itemize}
            \item every vertex covered by $F_2$ has degree $\Delta - 1$ in $G_2$,
            \item $\Delta(G_2 - F_2 - F_3 - N_{1,0} - \cdots - N_{1,s-1}) \leq \Delta - 2$, and
            \item every edge in $N_{1,i}$ has at least one endpoint in $\mathcal{S}_1[i]$.
        \end{itemize}   
        Taking $F \defeq F_0 \cup F_1 \cup F_2 \cup F_3$ and $N_i \defeq N_{0,i} \cup N_{1,i}$ works. Indeed, as in the proof of Claim~\ref{claim:3}, every vertex is incident to at most $3$ edges in $F$ because a vertex covered by both $F_0$ and $F_1$ must be uncovered by $F_2$. It remains to note that each $N_i$ is a matching in $G$ since $\dist_G(\mathcal{S}_0[i], \mathcal{S}_1[i]) \geq 3$.  
        \end{claimproof}

        Let $F$, $N_0$, \ldots, $N_{s-1}$ be given by Claim~\ref{claim:3+error} and let $G' \defeq G - N_0 - \cdots - N_{s-1}$. To finish the proof, we repeat the last part of the proof of Lemma~\ref{lemma:3match} to find Borel matchings $M_0$, $M_1$, $M_2 \subseteq F$ such that $\Delta(G' - M_0 - M_1 - M_2) \leq \Delta - 2$. In brief, we let
        \[
            T\,\defeq\, \set{e \in F \,:\, \text{both endpoints of $e$ are incident to $3$ edges in $F$}}
        \]
        and pick a Borel maximal matching $M \subseteq T$. Let $F' \defeq F \setminus M$. By construction, $\Delta(G' - F') \leq \Delta - 2$ (see Claim~\ref{claim:still_low_deg}). Let $H$ be the spanning subgraph of $G$ with edge set $F'$. Then $H$ has maximum degree at most $3$ and no two vertices of degree $3$ in $H$ are adjacent. By Corollary~\ref{corl:3}\ref{item:3_asi}, $\chi_\mathsf{B}'(H) \leq 3$. In other words, the edge set of $H$ (i.e., the set $F'$) can be covered by $3$ Borel matchings, as desired.
    \end{scproof}


    \begin{scproof}[ of Theorem~\ref{theo:asi_Shannon_restated}]
        Let $t \defeq \lfloor \Delta/2 \rfloor$ and note that $q = \lfloor 3\Delta/2\rfloor$ is equal to $3t$ if $\Delta$ is even and $3t+1$ if $\Delta$ is odd. Using Lemma~\ref{lemma:shatter}, we find a $3$-spaced family  $(\mathcal{S}_{j,0}, \mathcal{S}_{j,1})_{1 \leq j \leq t}$
        of $(s,3)$-shattering systems in $G$. Starting with $G_0 \defeq G$, we iteratively apply Lemma~\ref{lemma:3match+error} to obtain a sequence of Borel graphs $G_0 \supseteq G_1 \supseteq \cdots \supseteq G_t$ and Borel matchings $(M_{j,0}, M_{j,1}, M_{j,2}, N_{j,0}, \ldots, N_{j,s-1})_{1 \leq j \leq t}$, where:
        \begin{itemize}
            \item $G_{j} = G_{j-1} - M_{j,0} - M_{j,1} - M_{j,2} - N_{j,0}- \cdots - N_{j,s-1}$ for all $1 \leq j \leq t$,
            \item $\Delta(G_j) \leq \Delta - 2j$ for all $0 \leq j \leq t$,
            \item every edge in each matching $N_{j,i}$ has at least one endpoint in $\mathcal{S}_{j,0}[i] \cup \mathcal{S}_{j,1}[i]$.
        \end{itemize}
        For each index $0 \leq i < s$, we define $N_i \defeq \bigcup_{j = 1}^t N_{j,i}$ and note that $N_i$ is a matching in $G$ because the family $(\mathcal{S}_{j,0}, \mathcal{S}_{j,1})_{1 \leq j \leq t}$ is $3$-spaced. If $\Delta$ is even, then $\Delta(G_t) = 0$, i.e., $G_t$ has no edges, and hence
        \[
            E \,=\, \bigcup_{j = 1}^t \left(M_{j,0} \cup M_{j,1} \cup M_{j,2}\right)\,\cup\, \bigcup_{i=0}^{s-1} N_i.
        \]
        On the other hand, if $\Delta$ is odd, then $\Delta(G_t) \leq 1$, so the edge set $E'$ of $G_t$ is a matching and
        \[
            E \,=\, \bigcup_{j = 1}^t \left(M_{j,0} \cup M_{j,1} \cup M_{j,2}\right)\,\cup\, \bigcup_{i=0}^{s-1} N_i \,\cup\, E'.
        \]
        In both cases, $E$ is covered by $q + s$ Borel matchings, which means that $\chi'_\mathsf{B}(G) \leq q + s$, as desired.
    \end{scproof}

    \subsection{Almost balanced orientations}\label{subsec:orientations}

    To prove Theorem~\ref{theo:balance}, we need the following simple fact, which in essence is  
    a Borel version of the Schr\"oder--Bernstein theorem. The observation that the Schr\"oder--Bernstein theorem can be viewed as a statement about matchings in bipartite graphs 
    dates back to K\H{o}nig \cite[\S6.3]{konig_book}.

    \begin{prop}[{Borel Schr\"oder--Bernstein}]\label{prop:Bernstein}
        Let $G = (V,E,\pi)$ be a Borel graph with a bipartition $(A,B)$ and let $X \subseteq A$, $Y \subseteq B$ be Borel sets. Suppose $M_X$, $M_Y$ are Borel matchings in $G$ such that $M_X$ covers $X$ and $M_Y$ covers $Y$. Then there is a Borel matching $M \subseteq M_X \cup M_Y$ that covers $X \cup Y$.
    \end{prop}
    \begin{scproof}
        Without loss of generality, we may assume that $E = M_X \cup M_Y$. In particular, $\Delta(G) \leq 2$. Let $F$ be the set of all edges $e \in E$ such that the connected component of $G$ containing $e$ includes a vertex $x \in X$ with $\deg_G(x) = 1$. Taking $M \defeq (M_X \cap F) \cup (M_Y \setminus F)$ works (as in the usual proof of the Schr\"oder--Bernstein Theorem avoiding the Axiom of Choice; see, e.g., \cite{SB}).
    \end{scproof}

    Now we present the proof of Theorem~\ref{theo:balance}\ref{item:balance_measure}, restated here for convenience:

    \begin{theocopy}{theo:balance}
        {\normalfont \ref{item:balance_measure}} Let $G = (V,E,\pi)$ be a $d$-regular Borel graph for some $d \in \N$. If $\mu$ is a probability measure on $V$, then $G$ has a Borel orientation that is almost balanced at $\mu$-almost every vertex.
    \end{theocopy}
    %
    \begin{scproof}
        We introduce the following quantities: $k^- \defeq \lfloor (d-1)/2 \rfloor$ and $k^+ \defeq \lceil (d+1)/2 \rceil$. 
        We seek a Borel orientation $\mathcal{O}$ of $G$ such that $k^- \leq \deg_\mathcal{O}^+(x) \leq k^+$ for $\mu$-almost all $x \in V$. Let
        \[
            V^- \,\defeq\, V \times \set{1, \ldots, k^-} \qquad \text{and} \qquad V^+ \,\defeq\, V \times \set{1, \ldots, k^+}.
        \]
        Without loss of generality, $E \cap V^+ = \0$. Define a simple Borel graph $H$ with a bipartition $(E, V^+)$ by putting an edge between $e \in E$ and $(x,i) \in V^+$ if and only if $x \in \pi(e)$. 
        The vertex set of $H$ is equipped with the pushforward measure $\nu \defeq \iota_\ast(\mu)$, where $\iota$ is the embedding $\iota \colon V \to E \cup V^+ \colon x \mapsto (x,1)$. 
        
        The bipartition $(E, V^+)$ of $H$ is of type $(2k^+, d)$. Since $2k^+ \geq d+1 > d$, Corollary~\ref{corl:invariant} yields a $\nu$-null $H$-invariant Borel set $W \subseteq E \cup V^+$ and a Borel matching in $H$ covering every point in $E \setminus W$. By the construction of $H$ and the definition of $\nu$, the Borel set $X \defeq \set{x \in V \,:\, (x,1) \in W}$ is $\mu$-null and $G$-invariant. Thus, we may without loss of generality replace $G$ by $G - X$ (the edges in $G[X]$ may be oriented arbitrarily) 
        and assume that there is a Borel matching in $H$ that fully covers $E$.

        Next  we observe that the bipartition $(V^-, E)$ of the subgraph $H[V^-, E]$ is of type $(d, 2k^-)$. Since $d > d-1 \geq 2k^-$, Corollary~\ref{corl:invariant} again allows us to remove some $\mu$-null $G$-invariant Borel set from $V$ to arrange that there is a Borel matching in $H$ covering $V^-$. By Proposition~\ref{prop:Bernstein}, we now conclude that $H$ has a Borel matching $M$ that covers both $E$ and $V^-$.

        Finally, we define a Borel orientation $\mathcal{O}$ of $G$ by directing each edge $e \in E$ away from the unique vertex $x \in \pi(e)$ such that $\set{e, (x,i)} \in M$ for some $1 \leq i \leq k^+$. Consider any $x \in V$ and note that
        \[
            \deg^+_\mathcal{O}(x) \,=\, |\set{i \,:\,  \text{$(x,i) \in V^+$ is covered by $M$}}|.
        \]
        This quantity is at most $k^+$ and, since $M$ covers $V^-$, it is also at least $k^-$, as desired.
    \end{scproof}

    The other parts of Theorem~\ref{theo:balance} are proved in exactly the same way, but with Corollary~\ref{corl:invariant} replaced by Theorems~\ref{theo:growth}, \ref{theo:ends1}, and \ref{theo:BM_invariant}.

    \subsection{Equidecomposition theory}\label{subsec:tilings}

    \subsubsection{Preliminary observations}

    In \S\ref{subsec:tilings}, we shall describe applications of our techniques in equidecomposition theory. We start with some basic observations about the structure of type semigroups.

    Let $\G \acts X$ be an action of a group $\G$ on a set $X$ and let $\algebra{A} \subseteq \powerset{X}$ be a $\G$-invariant subalgebra. It will be convenient to relate the type semigroup $\mathsf{S}(X,\G,\algebra{A})$ with the type semigroups of certain ``enlarged actions'' on spaces built out of several disjoint copies of $X$.  
    This construction is standard in the equidecomposition theory literature; see, for example, \cite[\S10]{Wagon}.

    For a nonempty finite set $I$, we let $\Sym{I}$ be the group of all permutations of $I$. The set $X \times I$ is equipped with the action $\G \times \Sym{I} \acts X \times I$ given by 
    \[
        (\gamma,\sigma) \cdot (x,i) \,\defeq\, (\gamma \cdot x, \sigma(i)).
    \]
    That is, we view $X \times I$ as a family of disjoint copies of $X$ indexed by $I$ and let $\G$ act on each copy of $X$ separately, while $\Sym{I}$ permutes the copies. For $A \subseteq X \times I$ and $i \in I$, the \emphd{$i$-th level} of $A$ is 
    \[
        \level_i(A) \,\defeq\, \set{x \in X \,:\, (x,i) \in A}.
    \]
    Define a $(\G \times \Sym{I})$-invariant subalgebra $A^{\oplus I} \subseteq \powerset{X\times I}$ by
    \[
        \algebra{A}^{\oplus I} \,\defeq\, \set{A \subseteq X \times I \,:\, \level_i(A) \in \algebra{A} \text{ for all } i \in I}.
    \]
    The type semigroups $\mathsf{S}(X,\G,\algebra{A})$ and $\mathsf{S}(X\times I,\G\times \Sym{I},\algebra{A}^{\oplus I})$ are naturally isomorphic: 

    \begin{prop}\label{prop:congruence}
        Let $\G \acts X$ be an action of a group $\G$ on a set $X$ and let $\algebra{A} \subseteq \mathcal{P}(X)$ be a $\G$-invariant subalgebra. Fix a nonempty finite set $I$. We have an isomorphism
        \[
            \phi \colon \mathsf{S}(X\times I,\G\times \Sym{I},\algebra{A}^{\oplus I}) \xrightarrow{\sim} \mathsf{S}(X,\G,\algebra{A}),
        \]
        which is defined on the generators $[A]$ of $\mathsf{S}(X\times I,\G\times \Sym{I},\algebra{A}^{\oplus I})$ by the formula
        \[
            \phi \left([A]\right) \,\defeq\, \sum_{i \in I} [\level_i(A)] \,\in\, \mathsf{S}(X,\G,\algebra{A}). 
        \]
        The inverse of $\phi$ applied to a generator $[A]$ of $\mathsf{S}(X,\G,\algebra{A})$ can be computed via
        \[
            \phi^{-1} \left([A]\right) \,=\, [A \times \set{i}] \,\in\, \mathsf{S}(X\times I,\G\times \Sym{I},\algebra{A}^{\oplus I}), 
        \]
        where $i \in I$ is an arbitrary element. 
    \end{prop}
    \begin{scproof}
        Routine direct verification.
    \end{scproof}


    \subsubsection{Proof of Theorem~\ref{theo:Borel_comparison_asi}}

    To prove Theorem~\ref{theo:Borel_comparison_asi}, we follow, with suitable modifications, the strategy used in \cite[proof of Thm.~10.20]{Wagon}, which dates back to the insightful early work of K\H{o}nig \cite{KonigFrench}. For simplicity, we first establish the following special case of Theorem~\ref{theo:Borel_comparison_asi}, from which the full result will follow easily.

    \begin{lemma}\label{lemma:Borel_comparison_asi_for_sets}
        Let $\G \acts X$ be a Borel action of a finitely generated group $\G$ on a standard Borel space $X$ with Borel $\sigma$-algebra $\algebra{B}$. Suppose either that $\G$ is of subexponential growth, or that the asymptotic separation index of the action $\G \acts X$ is finite. Let $A$, $B \subseteq X$ be disjoint Borel sets. If for some $n > m \in \N$, 
        we have $n[A] \leq m[B]$ in $\mathsf{S}(X,\G)$, then $[A] \leq [B]$ in $\mathsf{S}(X,\G,\algebra{B})$.
    \end{lemma}
    \begin{scproof}
        Let $I \defeq \set{1,\ldots,n}$ and $J \defeq \set{1,\ldots, m}$. Note that $J \subset I$. Consider the isomorphism
        \[
            \phi \colon \mathsf{S}(X\times I, \G \times \Sym{I})  \xrightarrow{\sim} \mathsf{S}(X,\G)
        \]
        given by Proposition~\ref{prop:congruence}. Using the formula for $\phi$, we can write
        \[
            \phi([A \times I]) \,=\, n [A] \qquad \text{and} \qquad \phi([B \times J]) \,=\, m[B].
        \]
        Therefore, we have 
        $[A \times I] \leq [B \times J]$ in $\mathsf{S}(X\times I, \G \times \Sym{I})$, i.e., the set $A \times I$ is equidecomposable by $\G \times \Sym{I}$ with a subset of $B \times J$. Thus, we have a partition $A \times I = C_1 \sqcup \ldots \sqcup C_k$ and disjoint subsets $D_1$, \ldots, $D_k$ of $B \times J$ such that for all $1 \leq t \leq k$, there exist $\gamma_t \in \G$ and $\sigma_t \in \Sym{I}$ with
        \[
            D_t = (\gamma_t, \sigma_t) \cdot C_t. 
        \]
        We emphasize that the sets $C_t$, $D_t$ need not be Borel. 
        
        We define two graphs with a bipartition $(A,B)$. The first one, $G$, is simple and has edge set
        \[
            \set{\set{x, \gamma_t \cdot x} \,:\, \text{$x \in A$ and $1 \leq t \leq k$ such that $\gamma_t \cdot x \in B$}}.
        \]
        Note that the graph $G$ is Borel. Furthermore, $G$ is a subgraph of the Schreier graph $\mathsf{Sch}(X, \G,S)$ corresponding to any finite generating set $S$ for $\G$ with $\gamma_1$, \ldots, $\gamma_k \in S$ (see Definition~\ref{defn:schreier}). In particular, $G$ either is of subexponential growth or satisfies $ \asi (G) \leq \asi(\G \acts X) < \infty$.

        The second graph we define, $H$, is not necessarily Borel. It is a multigraph with vertex set $V \defeq A \sqcup B$, edge set $E \defeq A \times I$, and endpoint map $\pi \colon E \to [V]^2$ given by
        \[
            \pi(x,i) \,\defeq\, \set{x, \gamma_t \cdot x}, \qquad \text{where } (x,i) \in C_t.
        \]
        If $(x,i) \in C_t$, then $\gamma_t \cdot x$ is in one of the levels of $D_t$, so $\gamma_t \cdot x \in B$. Thus, $H$ is indeed a graph with a bipartition $(A,B)$. Furthermore, the underlying simple graph of $H$ is a subgraph of $G$.

        Next we observe that the bipartition $(A,B)$ of $H$ is of type $(n,m)$. Indeed, for each $x \in A$, there are exactly $n$ edges in $H$ that are incident to $x$, namely $(x,i)$ with $i \in I$. On the other hand, for $y \in B$, the number of edges in $H$ that are incident to $y$ is equal to the number of indices $j \in J$ such that $(y,j) \in D_1 \sqcup \ldots \sqcup D_k$, which cannot exceed $m$. Since $n > m$, we conclude that the bipartition $(A,B)$ is unbalanced in $H$. In particular, it has combinatorial expansion. The latter property only depends on the underlying simple graph of $H$ and is preserved by adding edges. Therefore, the bipartition $(A,B)$ has combinatorial expansion in $G$ as well.

        Now, using either Theorem~\ref{theo:growth} or Theorem~\ref{theo:asi_bounded_distance} as appropriate, we find $d \in \N$ and a Borel injective function $f \colon A \to B$ such that, for all $x \in A$,
        \[
            \dist_G(x,f(x)) \,\leq\, d.
        \]
        (If $G$ is of subexponential growth, Theorem~\ref{theo:growth} allows us to take $d = 1$.) The definition of $G$ implies that for all $x \in A$, we can write $f(x) = \gamma_x \cdot x$, where $\gamma_x \in \G$ belongs to the finite set
        \[
            F \,\defeq\, \left\{\gamma_1, \ldots, \gamma_k, \gamma_1^{-1}, \ldots, \gamma_k^{-1}, \mathbf{1}\right\}^d. 
        \]
        The assignment $x \mapsto \gamma_x$ can be made Borel (e.g., using Luzin--Novikov Uniformization). 
        Letting \[A_\gamma \,\defeq\, \set{x \in A \,:\, \gamma_x  = \gamma},\] we obtain a partition $A = \bigsqcup_{\gamma \in F} A_\gamma$ of $A$ into finitely many Borel pieces. Since $f$ is injective, the sets \[B_\gamma \,\defeq\, \gamma \cdot A_\gamma \,=\, \set{f(x) \,:\, x \in A_\gamma} \,\subseteq\, B\] are pairwise disjoint. This shows that $A$ is equidecomposable with a Borel subset of $B$---namely with $\bigsqcup_{\gamma \in F} B_\gamma$---using Borel pieces. In other words, $[A] \leq [B]$ in $\mathsf{S}(X,\G,\algebra{B})$, as desired.
    \end{scproof}

    \begin{scproof}[ of Theorem~\ref{theo:Borel_comparison_asi}]
        Recall that the set-up for Theorem~\ref{theo:Borel_comparison_asi} is the same as for Lemma~\ref{lemma:Borel_comparison_asi_for_sets}, but instead of a pair of Borel sets, we are given a pair of elements $\alpha$, $\beta \in \mathsf{S}(X,\G,\algebra{B})$ such that $n\alpha \leq m\beta$ in $\mathsf{S}(X,\G)$ for some $n > m \in \N$, and our goal is to show that $\alpha \leq \beta$ in $\mathsf{S}(X,\G,\algebra{B})$. 

        Write $\alpha = [A_1] + \cdots + [A_k]$ and $\beta = [B_1] + \cdots + [B_\ell]$, where $A_1$, \ldots, $A_k$, $B_1$, \ldots, $B_\ell$ are Borel subsets of $X$, and let $I \defeq \set{1, \ldots, k + \ell}$. Consider the sets
        \[
            A \,\defeq\, \bigsqcup_{i=1}^k (A_i \times \set{i}) \qquad \text{and} \qquad B \,\defeq\, \bigsqcup_{j=1}^\ell (B_j \times \set{k + j}).
        \]
        Note that $\algebra{B}^{\oplus I}$ is just the Borel $\sigma$-algebra on $X \times I$. By Proposition~\ref{prop:congruence}, we have an isomorphism
        \[
            \phi \colon \mathsf{S}(X\times I, \G \times \Sym{I})  \xrightarrow{\sim} \mathsf{S}(X,\G). 
        \]
        Using the formula for $\phi$, we see that
        \[
            \phi([A]) \,=\, [A_1] + \cdots + [A_k] \,=\, \alpha \qquad \text{and} \qquad \phi([B]) \,=\, [B_1] + \cdots + [B_\ell] \,=\, \beta.
        \]
        It follows that $n [A] \leq m[B]$ in $\mathsf{S}(X\times I, \G \times \Sym{I})$. Now we make two simple observations:
        \begin{itemize}
            \item If $\G$ is of subexponential growth, then so is $\G \times \Sym{I}$.

            \item The actions $\G \acts X$ and $\G \times \Sym{I} \acts X \times I$ have the same asymptotic separation index.
        \end{itemize}
        Since $A$ and $B$ are disjoint Borel subsets of $X \times I$, Lemma~\ref{lemma:Borel_comparison_asi_for_sets} shows that $[A] \leq [B]$ in the type semigroup $\mathsf{S}(X \times I, \G \times \Sym{I}, \algebra{B}^{\oplus I})$. Thanks to the isomorphism
        \[
            \mathsf{S}(X\times I, \G \times \Sym{I},\algebra{B}^{\oplus I}) \,\cong\, \mathsf{S}(X,\G,\algebra{B}),
        \]
        this implies that $\alpha \leq \beta$ in $\mathsf{S}(X, \G, \algebra{B})$, and the proof is complete.
    \end{scproof}

    \subsubsection{Almost unperforation: Proof of Theorem~\ref{theo:unperf}}\label{subsubsec:unperf}

    Now we turn our attention to Theorem~\ref{theo:unperf}, i.e., the fact that type semigroups corresponding to certain invariant $\sigma$-algebras are almost unperforated. The following fact is well known (as observed by K\H{o}nig \cite{KonigFrench}, it is a consequence of Hall's Theorem~\ref{theo:Hall}):

    \begin{theorem}[{\cite[Thm.~10.20]{Wagon}}]\label{theo:truly_unperforated}
       If $\G \acts X$ is an action of a group $\G$ on a set $X$, then the 
       semigroup $\mathsf{S}(X,\G)$ is 
       {\upshape\emphd{unperforated}}, i.e., for all $\alpha$, $\beta \in \mathsf{S}(X,\G)$ and $n \in \N^+$, $n \alpha \leq n\beta$ implies $\alpha \leq \beta$. 
   \end{theorem}
    
    Parts \ref{item:unperf_subexp} and \ref{item:unperf_asi} of Theorem~\ref{theo:unperf} are immediate corollaries of Theorem~\ref{theo:Borel_comparison_asi}:

    \begin{theocopy}{theo:unperf}\label{theo:unp_copy_Borel}
        {\normalfont\ref{item:unperf_subexp}, \ref{item:unperf_asi}} Let $\G \acts X$ be an action of a finitely generated group $\G$ on a standard Borel space $X$ and let $\algebra{B}$ be the Borel $\sigma$-algebra of $X$. Suppose either that $\G$ is of subexponential growth, or that $\asi(\G \acts X) < \infty$. Then the type semigroup $\mathsf{S}(X,\G,\algebra{B})$ is almost unperforated.
    \end{theocopy}
    \begin{scproof}
        Let $\alpha$, $\beta \in \mathsf{S}(X,\G,\algebra{B})$ be such that for some $n > m \in \N$, we have $n \alpha \leq m \beta$ in $\mathsf{S}(X,\G,\algebra{B})$. Then we also have $n \alpha \leq m \beta$ in $\mathsf{S}(X,\G)$, and thus $\alpha \leq \beta$ in  $\mathsf{S}(X,\G,\algebra{B})$ by Theorem~\ref{theo:Borel_comparison_asi}.
    \end{scproof}

    Part~\ref{item:unperf_BM} of Theorem~\ref{theo:unperf} also follows from Theorem~\ref{theo:Borel_comparison_asi}. We need a lemma that allows us to assume that Baire-measurable sets appearing in the proof are in fact Borel: 

    \begin{lemma}\label{lemma:inv_make_Borel}
        Let $\G \acts X$ be an action of a countable group $\G$ on a standard Borel space $X$, and let $\tau$ be a compatible topology on $X$. Suppose that the $\sigma$-algebra $\algebra{BM}$ of all $\tau$-Baire-measurable sets is $\G$-invariant. Let $\mathcal{U} \subseteq \algebra{BM}$ be a countable family of $\tau$-Baire-measurable sets. Then there exists a partition $X = Y \sqcup Z$ into $\G$-invariant Borel sets such that:
            \begin{itemize}
                \item the set $Y$ is $\tau$-meager,
                \item the action $\G \acts Z$ is Borel, and
                \item for all $U \in \mathcal{U}$, the set $U \cap Z$ is Borel.
            \end{itemize}
    \end{lemma}
    \begin{scproof}
        Let $I \subseteq X$ be the set of all $\tau$-isolated points of $X$ and let $I^* \defeq \G \cdot I$. Note that the sets $I$ and $I^*$ are countable (hence Borel) and $I^*$ is $\G$-invariant.
        We claim that if $M \subseteq X \setminus I^*$ is $\tau$-meager, then for all $\gamma \in \G$, $\gamma \cdot M$ is also $\tau$-meager. Indeed, every subset of $M$ is $\tau$-meager, so, by the $\G$-invariance of $\algebra{BM}$, every subset of $\gamma \cdot M$ must be $\tau$-Baire-measurable. Since $\gamma \cdot M$ contains no isolated points of $X$, this is only possible if $\gamma \cdot M$ is $\tau$-meager, as claimed. 

        Without loss of generality, we may assume that the family $\mathcal{U}$ is $\G$-invariant and includes a generating set for the Borel $\sigma$-algebra on $X$. Since each set $U \in \mathcal{U}$ is $\tau$-Baire-measurable, there is a $\tau$-meager Borel set $M_U \subseteq X$ such that $U \setminus M_U$ is Borel. Define
        \[
            Y \,\defeq\, \bigcup_{U \in \mathcal{U}, \, \gamma \in \G} \left(\gamma \cdot (M_U \setminus I^*)\right) \qquad \text{and} \qquad Z \,\defeq\, X \setminus Y.
        \]
        Then $Y$ is a $\G$-invariant Borel set. Furthermore, it is a countable union of $\tau$-meager sets, so it is itself $\tau$-meager. For each $U \in \mathcal{U}$, we have
        \[
            U \cap Z \,=\, ((U \setminus M_U) \cup (U \cap I^*)) \cap Z.
        \]
        Since the sets $U \setminus M_U$, $U \cap I^*$ (which is countable), and $Z$ are Borel, it follows that $U \cap Z$ is Borel as well. Finally, since $\mathcal{U}$ is $\G$-invariant and includes a generating set for the Borel $\sigma$-algebra on $X$, we conclude that the action $\G \acts Z$ is Borel, completing the proof.
    \end{scproof}

    \begin{theocopy}{theo:unperf}\label{theo:unperf_BM_copy}
        {\normalfont\ref{item:unperf_BM}} Let $\G \acts X$ be an action of a finitely generated group $\G$ on a standard Borel space $X$ and let $\algebra{BM}$ be the $\sigma$-algebra of $\tau$-Baire-measurable sets for some compatible topology $\tau$ on $X$. Suppose $\algebra{BM}$ is $\G$-invariant. Then the type semigroup $\mathsf{S}(X,\G,\algebra{BM})$ is almost unperforated.
    \end{theocopy}
    \begin{scproof}
        Suppose $\alpha = \sum_{i=1}^k [A_i]$ and $\beta = \sum_{j=1}^\ell [B_j]$ are elements of $\mathsf{S}(X,\G,\algebra{BM})$ such that $n\alpha \leq m\beta$ for some $n > m \in \N$. Here the sets $A_i$ and $B_j$ belong to $\algebra{BM}$. By Lemma~\ref{lemma:inv_make_Borel}, we can find a partition $X = Y \sqcup Z$ into $\G$-invariant Borel sets such that:
        \begin{itemize}
                \item the set $Y$ is $\tau$-meager,
                \item the action $\G \acts Z$ is Borel, and
                \item for all $1 \leq i \leq k$ and $1 \leq j \leq \ell$, the sets $A_i \cap Z$ and $B_j \cap Z$ are Borel.
            \end{itemize}
        Furthermore, thanks to Theorem~\ref{theo:asi_ae}, we may assume that $\asi(\G \acts Z) \leq 1$. Let
        \[
            \alpha \cap Y \,\defeq\, \sum_{i=1}^k [A_i \cap Y], \qquad \alpha \cap Z \,\defeq\, \sum_{i=1}^k [A_i \cap Z],
        \]
        and define $\beta \cap Y$ and $\beta \cap Z$ analogously. Since the sets $Y$ and $Z$ are $\G$-invariant, we have
        \[
            n(\alpha \cap Y) \,\leq\, m (\beta \cap Y) \qquad \text{and} \qquad n(\alpha \cap Z) \,\leq\, m(\beta \cap Z).
        \]
        By Theorem~\ref{theo:truly_unperforated}, the first of these inequalities implies that $\alpha \cap Y \leq \beta \cap Y$ in $\mathsf{S}(Y,\G)$. Since all subsets of $Y$ are $\tau$-Baire-measurable, we also have $\alpha \cap Y \leq \beta \cap Y$ in $\mathsf{S}(X,\G,\algebra{BM})$. On the other hand, since $\asi(\G \acts Z) \leq 1$, Theorem~\ref{theo:Borel_comparison_asi} shows that $\alpha \cap Z \leq \beta \cap Z$ in $\mathsf{S}(Z,\G,\algebra{B})$, where $\algebra{B}$ is the Borel $\sigma$-algebra on $Z$. In particular, $\alpha \cap Z \leq \beta \cap Z$ in $\mathsf{S}(X,\G,\algebra{BM})$. To complete the proof, we observe that, in $\mathsf{S}(X,\G,\algebra{BM})$,
            $\alpha = \alpha \cap Y + \alpha \cap Z \leq \beta \cap Y + \beta \cap Z = \beta$, 
        as desired.
    \end{scproof}

    It remains to verify part \ref{item:unperf_meas} of Theorem~\ref{theo:unperf}, which concerns the $\sigma$-algebra of $\mu$-measurable sets. We will follow a strategy similar to the one employed in the proof of Theorem~\ref{theo:Borel_comparison_asi}, except that Theorems~\ref{theo:growth} and \ref{theo:asi_bounded_distance} will be replaced by Theorem~\ref{theo:measurable}. We will require a version of Lemma~\ref{lemma:inv_make_Borel} for $\mu$-measurable sets:

    \begin{lemma}\label{lemma:inv_make_meas_Borel}
        Let $\G \acts X$ be an action of a countable group $\G$ on a standard Borel space $X$, and let $\mu$ be a probability measure on $X$. Suppose that the $\sigma$-algebra $\algebra{M}$ of all $\mu$-measurable sets is $\G$-invariant. Let $\mathcal{U} \subseteq \algebra{M}$ be a countable family of $\mu$-measurable sets. Then there exists a partition $X = Y \sqcup Z$ into $\G$-invariant Borel sets such that:
            \begin{itemize}
                \item the set $Y$ is $\mu$-null,
                \item the action $\G \acts Z$ is Borel, and
                \item for all $U \in \mathcal{U}$, the set $U \cap Z$ is Borel.
            \end{itemize}
    \end{lemma}
    \begin{scproof}
        The proof is the same as that of Lemma~\ref{lemma:inv_make_Borel}, modulo replacing $\tau$-Baire-measurable sets by $\mu$-measurable ones, $\tau$-meager sets by $\mu$-null ones, and $\tau$-isolated points by atoms of $\mu$.
    \end{scproof}

    As in the proof of Theorem~\ref{theo:Borel_comparison_asi}, we will reduce Theorem~\ref{theo:unperf}\ref{item:unperf_meas} to the following special case:

    \begin{lemma}\label{lemma:unp_meas_sets}
        Let $\G \acts X$ be an action of a finitely generated group $\G$ on a standard Borel space $X$ and let $\algebra{M}$ be the $\sigma$-algebra of $\mu$-measurable sets for some probability measure $\mu$ on $X$. Suppose $\algebra{M}$ is $\G$-invariant. Let $A$, $B \subseteq X$ be disjoint $\mu$-measurable sets. If for some natural numbers $n >m$, we have $n[A] \leq m[B]$ in $\mathsf{S}(X,\G, \algebra{M})$, then $[A] \leq [B]$ in $\mathsf{S}(X,\G, \algebra{M})$ as well.
    \end{lemma}
    \begin{scproof}
        As in the proof of Lemma~\ref{lemma:Borel_comparison_asi_for_sets}, we let $I \defeq \set{1,\ldots,n}$ and $J \defeq \set{1,\ldots, m}$ and apply Proposition~\ref{prop:congruence} to obtain the isomorphism
        \[
            \phi \colon \mathsf{S}(X\times I, \G \times \Sym{I}, \algebra{M}^{\oplus I})  \xrightarrow{\sim} \mathsf{S}(X,\G, \algebra{M}).
        \]
        Using the formula for $\phi$, we have $\phi([A \times I]) = n[A]$ and $\phi([B \times J]) = m[B]$, 
        so 
        $[A \times I] \leq [B \times J]$ in $\mathsf{S}(X\times I, \G \times \Sym{I}, \algebra{M}^{\oplus I})$. 
        Notice that $\algebra{M}^{\oplus I}$ is precisely the $\sigma$-algebra of $\mu^{\oplus I}$-measurable subsets of $X \times I$, where $\mu^{\oplus I}$ is the product of $\mu$ with the uniform probability measure on $I$. 
        Therefore, we have a partition $A \times I = C_1 \sqcup \ldots \sqcup C_k$ and disjoint subsets $D_1$, \ldots, $D_k$ of $B \times J$ such that for all $1 \leq t \leq k$, the sets $C_t$, $D_t$ are $\mu^{\oplus I}$-measurable, and there are $\gamma_t \in \G$, $\sigma_t \in \Sym{I}$ with
        \[
            D_t = (\gamma_t, \sigma_t) \cdot C_t. 
        \]

        Every $(\G \times \Sym{I})$-invariant subset of $X \times I$ has the form $Y \times I$ for some $\G$-invariant set $Y \subseteq X$. Hence, by Lemma~\ref{lemma:inv_make_meas_Borel}, there is a partition $X = Y \sqcup Z$ into $\G$-invariant Borel sets such that:
        \begin{itemize}
                \item the set $Y$ is $\mu$-null,
                \item the action $\G \acts Z$ is Borel, and
                \item for all $1 \leq t \leq k$, the sets $C_t \cap (Z \times I)$ and $D_t \cap (Z \times I)$ are Borel.
            \end{itemize}
            As in the proof of Theorem~\ref{theo:unperf_BM_copy}, we have $n[A \cap Y] \leq m[B \cap Y]$, which implies that $[A \cap Y] \leq [B \cap Y]$ in $\mathsf{S}(Y,\G)$ by Theorem~\ref{theo:truly_unperforated}. Since all subsets of $Y$ are $\mu$-measurable, we also have $[A \cap Y] \leq [B \cap Y]$ in $\mathsf{S}(X,\G,\algebra{M})$, and thus it only remains to show that $[A \cap Z] \leq [B \cap Z]$ in $\mathsf{S}(X,\G,\algebra{M})$. Therefore, from now on we shall without loss of generality replace $X$ by $Z$ are assume that:
            \begin{itemize}
                \item the action $\G \acts X$ is Borel, and
                \item for all $1 \leq t \leq k$, the sets $C_t$ and $D_t$ are Borel.
            \end{itemize}
            In particular, the sets $A$ and $B$ are Borel. 

            Again, we follow the proof of Lemma~\ref{lemma:Borel_comparison_asi_for_sets} and define a multigraph $H$ with vertex set $V \defeq A \sqcup B$, edge set $E \defeq A \times I$, and endpoint map $\pi \colon E \to [V]^2$ given by
        \[
            \pi(x,i) \,\defeq\, \set{x, \gamma_t \cdot x}, \qquad \text{where } (x,i) \in C_t.
        \]
           Exactly as in the proof of Lemma~\ref{lemma:Borel_comparison_asi_for_sets}, the bipartition $(A,B)$ in $H$ is unbalanced of type $(n,m)$.  However, in contrast to the situation in the proof of Lemma~\ref{lemma:Borel_comparison_asi_for_sets}, the graph $H$ is Borel, because the sets $C_t$, $D_t$ are Borel. Therefore, we may apply Theorem~\ref{theo:measurable}, or, more precisely, its strengthening given by Corollary~\ref{corl:invariant}, to find a $\mu$-null $H$-invariant Borel set $U \subseteq V$ and a Borel matching $M$ in $H - U$ that covers every vertex in $A \setminus U$. By Hall's Theorem~\ref{theo:Hall}, we also have a (not necessarily Borel) matching $M'$ in $H[U]$ covering every vertex in $A \cap U$. Since the set $U$ is $H$-invariant, $M \sqcup M'$ is a matching in $H$ that covers $A$.
           
           For each $x \in A$, let $\gamma_x \in \set{\gamma_1, \ldots, \gamma_k}$ be a group element such that $x$ is matched to $\gamma_x \cdot x$ by $M \sqcup M'$, 
           where we make the assignment $x \mapsto \gamma_x$ Borel on $A \setminus U$. For each $1 \leq t \leq k$, the set
         \[A_t \,\defeq\, \set{x \in A \,:\, \gamma_x  = \gamma_t}\]
         is $\mu$-measurable, because $A_t \setminus U$ is Borel and $A_t \cap U$ is $\mu$-null. Hence, we have a partition $A = \bigsqcup_{t = 1}^k A_t$ of $A$ into finitely many $\mu$-measurable pieces. The sets $B_t \defeq \gamma_t \cdot A_t$ are pairwise disjoint subsets of $B$, so $A$ is equidecomposable with a subset of $B$ using $\mu$-measurable pieces. In other words, $[A] \leq [B]$ in the type semigroup $\mathsf{S}(X,\G,\algebra{M})$, as desired.
    \end{scproof}

    \begin{theocopy}{theo:unperf}\label{theo:unperf_meas_copy}
        {\normalfont\ref{item:unperf_meas}} Let $\G \acts X$ be an action of a finitely generated group $\G$ on a standard Borel space $X$ and let $\algebra{M}$ be the $\sigma$-algebra of $\mu$-measurable sets for some probability measure $\mu$ on $X$. Suppose $\algebra{M}$ is $\G$-invariant. Then the type semigroup $\mathsf{S}(X,\G,\algebra{M})$ is almost unperforated.
    \end{theocopy}
    \begin{scproof}
        The argument is exactly the same as the reduction of Theorem~\ref{theo:Borel_comparison_asi} to Lemma~\ref{lemma:Borel_comparison_asi_for_sets}. Let $\alpha$, $\beta \in \mathsf{S}(X,\G,\algebra{M})$ be such that $n\alpha \leq m\beta$ in $\mathsf{S}(X,\G, \algebra{M})$ for some $n > m \in \N$. Our goal is to show that $\alpha \leq \beta$ in $\mathsf{S}(X,\G,\algebra{M})$. To this end, we write $\alpha = [A_1] + \cdots + [A_k]$ and $\beta = [B_1] + \cdots + [B_\ell]$, where $A_1$, \ldots, $A_k$, $B_1$, \ldots, $B_\ell$ are $\mu$-measurable subsets of $X$. Let $I \defeq \set{1, \ldots, k + \ell}$ and define
        \[
            A \,\defeq\, \bigsqcup_{i=1}^k (A_i \times \set{i}) \qquad \text{and} \qquad B \,\defeq\, \bigsqcup_{j=1}^\ell (B_j \times \set{k + j}).
        \]
        By Proposition~\ref{prop:congruence}, we have an isomorphism
        \[
            \phi \colon \mathsf{S}(X\times I, \G \times \Sym{I},\algebra{M}^{\oplus I})  \xrightarrow{\sim} \mathsf{S}(X,\G, \algebra{M}) 
        \]
        such that $\phi([A])  = \alpha$ and $\phi([B]) = \beta$. Then $n [A] \leq m[B]$ in $\mathsf{S}(X\times I, \G \times \Sym{I}, \algebra{M}^{\oplus I})$, and, by Lemma~\ref{lemma:unp_meas_sets}, $[A] \leq [B]$ in $\mathsf{S}(X \times I, \G \times \Sym{I}, \algebra{M}^{\oplus I})$. Thus, $\alpha \leq \beta$ in $\mathsf{S}(X,\G,\algebra{M})$, as desired.
    \end{scproof}

    We also observe that if $\algebra{B}$ is the Borel $\sigma$-algebra of $X$, then the type semigroup $\mathsf{S}(X,\G,\algebra{B})$ may in general fail to be almost unperforated; in particular, the additional assumptions in parts~\ref{item:unperf_subexp} and \ref{item:unperf_asi} of Theorem~\ref{theo:unp_copy_Borel} cannot be removed.

    \begin{prop}\label{prop:no_borel_unperf}
        There exists a free Borel action $\G \acts X$ of a finitely generated group $\G$ on a standard Borel space $X$ with Borel $\sigma$-algebra $\algebra{B}$ such that the type semigroup $\mathsf{S}(X,\G,\algebra{B})$ is not almost unperforated.
    \end{prop}
    \begin{scproof}
        Fix $n > 2$ and let $\G \defeq (\Z/2\Z)^{* n}$ be the group that is freely generated by $n$ involutions $\gamma_1$, \ldots, $\gamma_n$. By Marks's Theorem~\ref{theo:grab}, there exists a simple $n$-regular Borel graph $G$ such that:
        \begin{itemize}
            \item every component of $G$ is an $n$-regular tree,
            \item the Borel chromatic index of $G$ is $n$, and
            \item the edge grabbing problem has no Borel solution on $G$.
        \end{itemize}
        (The edge grabbing problem is defined and discussed in \S\ref{subsec:Marks}.) Write $G = (X,E)$ and let $\algebra{B}$ be the Borel $\sigma$-algebra of $X$. Fix a proper Borel coloring $c \colon E \to \set{1,\ldots, n}$ of $G$. 
        For all $x \in X$ and $1 \leq i \leq n$, there exists a unique edge that is incident to $x$ and colored $i$; we let $\gamma_i \cdot x$ be its other endpoint. Since $\gamma_i \cdot (\gamma_i \cdot x) = x$, this definition gives rise to a Borel action $\G \acts X$, with $G$ being its Schreier graph corresponding to the generating set $\set{\gamma_1, \ldots, \gamma_n}$. Since every component of $G$ is a tree, the action $\G \acts X$ is free. We claim that $\mathsf{S}(X,\G,\algebra{B})$ is not almost unperforated.

        By Lemma~\ref{lemma:LNU}, there exists a Borel map $\mathsf{root} \colon E \to X$ that sends every edge $e \in E$ to one of its endpoints, which we refer to as the \emphd{root} of $e$. For each $1 \leq i \leq n$, let
        \[
            B_i \,\defeq\, \set{x \in X \,:\, \mathsf{root}(\set{x, \gamma_i \cdot x}) = x}.
        \]
        For every point $x \in X$, precisely one of $x$, $\gamma_i \cdot x$ belongs to $B_i$. In other words, $X = B_i \sqcup (\gamma_i \cdot B_i)$ is a Borel partition of $X$. It follows that, in $\mathsf{S}(X,\G,\algebra{B})$,
        \begin{equation}\label{eq:Bi}
            2[B_i] \,=\, [B_i] + [\gamma_i \cdot B_i] \,=\, [X].
        \end{equation}
        Now we consider the following two elements of $\mathsf{S}(X,\G,\algebra{B})$: 
        \[
            \alpha \,\defeq\, [X] \qquad \text{and} \qquad \beta\,\defeq\, \sum_{i=1}^n [B_i].
        \]
        By \eqref{eq:Bi}, we have $n\alpha = 2\beta$ in $\mathsf{S}(X,\G,\algebra{B})$. However, we will show that $\alpha \not \leq \beta$ in $\mathsf{S}(X,\G,\algebra{B})$, providing the desired counterexample to the almost unperforation of $\mathsf{S}(X,\G,\algebra{B})$.

        It will be convenient to let $I \defeq \set{1,\ldots, n}$ and use the isomorphism
        \[
            \mathsf{S}(X,\G,\algebra{B}) \,\cong\, \mathsf{S}(X\times I, \G \times \Sym{I},\algebra{B}^{\oplus I})
        \]
        given by Proposition~\ref{prop:congruence}. (Note that $\algebra{B}^{\oplus I}$ is just the Borel $\sigma$-algebra on $X \times I$.) The images of $\alpha$ and $\beta$ in $\mathsf{S}(X\times I, \G \times \Sym{I},\algebra{B}^{\oplus I})$ under this isomorphism are $[A]$ and $[B]$ respectively, where
        \begin{align*}
            A \,&\defeq\,  X \times \set{1},\\
            B \,&\defeq\, \set{(x,i) \,:\, \mathsf{root}(\set{x, \gamma_i \cdot x}) = x}.
        \end{align*}
        Suppose for contradiction that $[A] \leq [B]$ in $\mathsf{S}(X\times I, \G \times \Sym{I},\algebra{B}^{\oplus I})$, i.e., the set $A$ is equidecomposable by $\G \times \Sym{I}$ with a subset of $B$ using Borel pieces. This would give rise to a Borel injection $h \colon A \to B$ sending each point in $A$ to the point in $B$ it is moved to via the equidecomposition. We can then define a Borel injection $f \colon X \to E$ by composing the maps
        \[
            X \xrightarrow{\mbox{$x \mapsto (x,1)$}} A \xrightarrow{\mbox{\quad$h$\quad}} B \xrightarrow{\mbox{$(x,i) \mapsto \set{x, \gamma_i \cdot x}$}} E.
        \]
        Since the pieces of the equidecomposition are only moved by a finite set of group elements, there is some $d \in \N$ such that $\dist_G(x,f(x)) \leq d$ for all $x \in X$. Applying Proposition~\ref{prop:no_dist_matching} to the incidence graph between $X$ and $E$, we conclude that there exists a Borel injective function $g \colon X \to E$ such that every vertex $x \in X$ is incident to the edge $g(x)$. But this precisely means that $g$ is a Borel solution to the edge grabbing problem on $G$, which does not exist by assumption.
    \end{scproof}

    \subsubsection{F\o{}lner tilings and the comparison property}\label{subsubsec:tiling}

    In the last part of \S\ref{subsec:tilings}, we sketch the relationship between our work on type semigroups and the structure of amenable group actions. This discussion is purely expository and is included here solely to help place our results into context. In particular, we will describe certain previously known facts about amenable group actions that are strengthened and generalized by Theorem~\ref{theo:unperf}.

    For the remainder of \S\ref{subsec:tilings}, we fix a finitely generated group $\G$. Given a finite set $K \subseteq \G$ and $\epsilon > 0$, a \emphd{$(K,\epsilon)$-F\o{}lner set} is a nonempty finite subset $F \subseteq \G$ such that 
    \[
        |(KF) \symdif F| \,\leq\, \epsilon |F|.
    \]
    Similarly, if $\G \acts X$ is a free action on a set $X$, then a nonempty finite set $F \subseteq X$ is \emphd{$(K,\epsilon)$-F\o{}lner} 
    if 
    \[
        |(K \cdot F) \symdif F| \,\leq\, \epsilon |F|.
    \]
    The group $\G$ is \emphd{amenable} if it admits a $(K,\epsilon)$-F\o{}lner set $F \subseteq \G$ for every finite $K \subseteq \G$ and $\epsilon > 0$. A \emphd{F\o{}lner sequence} $(F_n)_{n \in \N}$ for $\G$ is a sequence of nonempty finite sets $F_n \subseteq \G$ such that, for all finite $K \subseteq \G$ and $\epsilon > 0$, all but finitely many of the sets $F_n$ are $(K,\epsilon)$-F\o{}lner. 
    Amenability is one of the most important concepts in  group theory; for general overviews, see \cite{amenable1,amenable2}. 
    
    From now on, we make the following standing assumption:

        \begin{assumption*}[for \S\ref{subsubsec:tiling}]
            {$\G$ is amenable and $(F_n)_{n \in \N}$ is a F\o{}lner sequence for $\G$.}
        \end{assumption*}


    \begin{definition}[Quasi-tilings and tilings]\label{defn:tiling}
        Let $\G \acts X$ be a free action of $\G$ on a set $X$. Fix a finite set $K \subseteq \G$ and $\epsilon > 0$. A \emphd{$(K,\epsilon)$-F\o{}lner quasi-tiling} of $X$ is a family $\mathcal{T}$ of pairwise disjoint $(K,\epsilon)$-F\o{}lner subsets of $X$ that satisfy the \emphd{finitely many shapes condition}: there exists a finite list $T_1$, \ldots, $T_n$ of subsets of $\G$, called the \emphd{tiles}, such that each $T \in \mathcal{T}$ is of the form
        \[
            T \,=\, T_i \cdot x \quad \text{for some } 1 \leq i \leq n \text{ and } x \in X.
        \]
        If the sets in $\mathcal{T}$ additionally cover $X$, then $\mathcal{T}$ is called a \emphd{$(K,\epsilon)$-F\o{}lner tiling} of $X$.
    \end{definition}

    Note that the finitely many shapes condition in Definition~\ref{defn:tiling} is equivalent to the statement that there exists a finite set $S \subseteq \G$ such that for all $T \in \mathcal{T}$ and $x \in T$, we have $T \subseteq S \cdot x$.

    We are interested in Borel actions $\G \acts X$ that admit Borel F\o{}lner tilings 
    (here we say that $\mathcal{T}$ is Borel if it is a Borel subset of the standard Borel space $\finset{X}$ of all finite subsets of $X$).
    
    \begin{definition}[Borel tiling property]
        A free Borel action $\G \acts X$ on a standard Borel space $X$ has the \emphd{Borel tiling property} if it has a Borel $(K,\epsilon)$-F\o{}lner tiling for all finite $K \subseteq \G$ and $\epsilon > 0$. 
    \end{definition}

    
    A closely related notion is the \emph{comparison property}. While it was originally defined and studied in the topological context, it is straightforward to adapt the definition for the Borel setting. 
    Our presentation mostly follows the paper \cite{DZ_comparison} by Downarowicz and Zhang, although essentially equivalent ideas also appear in the work of Kerr \cite{Kerr}, Kerr and Szab\'o \cite{KerrSz}, and 
    Conley, Jackson, Kerr, Marks, Seward, and Tucker-Drob~\cite{folner.tilings}.


    \begin{definition}[Banach density]
        Let $\G \acts X$ be a free action of $\G$. For a set $A \subseteq X$, we define the \emphd{lower} and \emphd{upper Banach densities} of $A$ \ep{with respect to the action $\G \acts X$} as, respectively,
        \[
            \underline{D}(A) \,\defeq\, \lim_{n \to \infty} \inf_{x \in X} \frac{|A \cap (F_n \cdot x) |}{|F_n|} \qquad \text{and} \qquad \overline{D}(A) \,\defeq\, \lim_{n \to \infty} \sup_{x \in X} \frac{|A \cap (F_n \cdot x) |}{|F_n|} .
        \]
        For sets $A$, $B \subseteq X$, the \emphd{Banach density advantage} of $B$ over $A$  is
        \[
            \underline{D}(B,A) \,\defeq\, \lim_{n \to \infty} \inf_{x \in X} \frac{|B \cap (F_n \cdot x)| \,-\, |A \cap (F_n \cdot x)|}{|F_n|}. 
        \]
        All these limits exist and do not depend on the F\o{}lner sequence $(F_n)_{n \in \N}$ \cite[Lem.~6.9]{DZ_comparison}.
    \end{definition}



    \begin{definition}[Borel comparison property]
        Let $\G \acts X$ be a free Borel action on a standard Borel space $X$ with Borel $\sigma$-algebra $\algebra{B}$. 
        If for all Borel sets  $A$, $B \subseteq X$ such that $\underline{D}(B,A) > 0$, we have $[A] \leq [B]$ in $\mathsf{S}(X,\G,\algebra{B})$, then the action $\G \acts X$ has the \emphd{Borel comparison property}. 


        
    \end{definition}


    \begin{theorem}[{ess.~Downarowicz--Zhang \cite[Cor.~7.6]{DZ_comparison}, Kerr--Szab\'o \cite[Thm.~6.1]{KerrSz}; see also Conley--Jackson--Kerr--Marks--Seward--Tucker-Drob \cite[\S3]{folner.tilings}}]\label{theo:DZ}
        Let $\G \acts X$ be a free Borel action of $\G$ on a standard Borel space $X$. 
        The following statements are equivalent:
        \begin{enumerate}[label=\ep{\normalfont\roman*}]
            \item\label{item:tiling} the action $\G \acts X$ has the Borel tiling property,\item\label{item:comp} the action $\G \acts X$ has the Borel comparison property.


        \end{enumerate}
    \end{theorem}

    The papers \cite{DZ_comparison,KerrSz} operate in a topological framework, but their arguments can be easily adjusted to obtain Theorem~\ref{theo:DZ} in the Borel setting. The paper \cite{folner.tilings} does work with Borel actions, but it does not explicitly use the term ``comparison property.'' For the reader's convenience, we sketch the proof of Theorem~\ref{theo:DZ} as it is stated here at the end of this subsection. 

    The following simple observation links the properties discussed above to our results:

    \begin{prop}\label{prop:unp_to_comp}
        Let $\G \acts X$ be a free Borel action of $\G$ on a standard Borel space $X$ with Borel $\sigma$-algebra $\algebra{B}$. If the type semigroup $\mathsf{S}(X,\G,\algebra{B})$ is almost unperforated, then the action $\G \acts X$ has the Borel comparison property \ep{and hence also the Borel tiling property}.
    \end{prop}
    \begin{scproof}
        Let $A$, $B \subseteq X$ be Borel sets such that $\underline{D}(B,A) > 0$. This means that, for any large enough $n \in \N$, we have $|A \cap (F_n \cdot x)| < |B \cap (F_n \cdot x)|$ for all $x \in X$. It follows that each point $x \in X$ is covered by the sets $(\gamma^{-1} \cdot B)_{\gamma \,\in\, F_n}$ strictly more times than by the sets $(\gamma^{-1} \cdot A)_{\gamma \,\in\, F_n}$. Thus, letting $m \defeq |F_n|$, we can perform the following calculation in the type semigroup $\mathsf{S}(X,\G,\algebra{B})$:
        \[
            m[B] \,=\, \sum_{\gamma \,\in\, F_n} [\gamma^{-1} \cdot B] \,\geq\, [X] \,+\, \sum_{\gamma \,\in\, F_n} [\gamma^{-1} \cdot A] \,=\, [X] + m[A] \,\geq\, (m+1)[A].
        \]
        Since $\mathsf{S}(X,\G,\algebra{B})$ is almost unperforated, we conclude that $[A] \leq [B]$ in  $\mathsf{S}(X,\G,\algebra{B})$, as desired.
    \end{scproof}

    Thanks to Proposition~\ref{prop:unp_to_comp}, the following combination of previously known results is an immediate consequence of parts \ref{item:unperf_subexp} and \ref{item:unperf_asi} of Theorem~\ref{theo:unperf}\footnote{Parts~\ref{item:unperf_meas} and \ref{item:unperf_BM} of Theorem~\ref{theo:unperf} can similarly be used to prove appropriate measurable/Baire-measurable versions of the tiling property and the comparison property. However, they already follow from Theorem~\ref{theo:tiling}\ref{item:asi_tile}, thanks to Theorem~\ref{theo:asi_ae} and the fact that actions of amenable groups are almost everywhere hyperfinite.}: 

    \begin{theorem}[{Downarowicz--Zhang \cite[Thms.~6.33, 7.5]{DZ_comparison}, Conley--Jackson--Marks--Seward--Tucker-Drob \cite[Thm.~9.1]{conley2020borel}}]\label{theo:tiling}
        Let $\G \acts X$ be a free Borel action of $\G$ on a standard Borel space $X$. Suppose that at least one of the following conditions holds:
        \begin{enumerate}[label=\ep{\normalfont\arabic*}]
            \item\label{item:subexp_tile} $\G$ is of subexponential growth, or
            \item\label{item:asi_tile} the action $\G \acts X$ has finite asymptotic separation index.
        \end{enumerate}
        Then the action $\G \acts X$ has the Borel comparison property, and hence also the Borel tiling property. 
    \end{theorem}

    Part~\ref{item:subexp_tile} of Theorem~\ref{theo:tiling} was proved (in a topological context) by Downarowicz and Zhang~\cite{DZ_comparison}, while part~\ref{item:asi_tile} was established by Conley, Jackson, Marks, Seward, and Tucker-Drob~\cite{conley2020borel}. We remark that our proof of Theorem~\ref{theo:tiling}\ref{item:asi_tile}, which ultimately relies on Borel matchings, appears fundamentally different from the argument in \cite{conley2020borel}. In particular, the approach taken in  \cite{conley2020borel} is to directly construct a Borel F\o{}lner tiling, which only makes sense for amenable groups. On the other hand, Theorem~\ref{theo:unperf}\ref{item:unperf_asi} is a meaningful  statement regardless of the group's amenability.

    To complete this section, we present the proof of Theorem~\ref{theo:DZ} (or, rather, a reduction from Theorem~\ref{theo:DZ} to facts that are explicitly stated in the literature). 

    \begin{scproof}[ of Theorem~\ref{theo:DZ}]
        Let $\algebra{B}$ be the Borel $\sigma$-algebra of $X$.
        To show \ref{item:tiling} $\Longrightarrow$ \ref{item:comp}, suppose that $A$, $B \subseteq X$ are Borel sets with $\underline{D}(B,A) > 0$. By \cite[Lem.~6.9]{DZ_comparison}, there exist a finite set $K \subseteq \G$ and $\epsilon > 0$ such that every $(K,\epsilon)$-F\o{}lner set $F \subseteq X$ satisfies $|B \cap F| > |A \cap F|$. Assuming the action has the Borel tiling property, we can find a Borel $(K,\epsilon)$-F\o{}lner tiling $\mathcal{T}$ of $X$. Then $|B \cap T| > |A \cap T|$ for all $T \in \mathcal{T}$, and thus we can assign to each set $T \in \mathcal{T}$ an injective map
        $
            f_T \colon (A \cap T) \to (B \cap T)
        $.
        Moreover, using Luzin--Novikov Uniformization, we can ensure that this assignment is Borel in an appropriate sense. Letting $f(x) \defeq f_T(x)$ for all $T \in \mathcal{T}$ and $x \in A \cap T$, we obtain a Borel injection $f \colon A \to B$. Thanks to the finitely many shapes condition satisfied by $\mathcal{T}$, there exists some finite set $S \subseteq \G$ such that $f(x) \in S \cdot x$ for all $x \in A$. Therefore, letting $A_\gamma \defeq \set{x \in A \,:\, f(x) = \gamma \cdot x}$ for all $\gamma \in S$, we obtain a partition of $A$ into finitely many Borel sets such that their $\G$-translates $\gamma \cdot A_\gamma$ are disjoint Borel subsets of $B$. It follows that $[A] \leq [B]$ in $\mathsf{S}(X,\G,\algebra{B})$, as desired. 

        Now we turn to the implication \ref{item:comp} $\Longrightarrow$ \ref{item:tiling}. Fix a nonempty finite set $K \subseteq \G$ and $\epsilon \in (0,1/2)$. We may assume $\G$ is infinite. Pick a finite subset $L \supset K$ of $\G$ with $|L| \geq |K|/\epsilon$ and let $\delta \defeq \epsilon/|K|$. 
        By \cite[Lem.~3.4]{folner.tilings}, there is a Borel $(L,\epsilon)$-F\o{}lner quasi-tiling $\mathcal{T}$ of $X$ such that the set $B \subseteq X$ of all points covered by  $\mathcal{T}$ satisfies
        $
            \underline{D}(B) \geq 1 - \delta
        $. 
        Let $A \defeq X \setminus B$, so $\overline{D}(A) \leq \delta$. 
        
        By Luzin--Novikov Uniformization, for each $T \in \mathcal{T}$, we can choose in a Borel way a subset $T' \subseteq T$ of cardinality $\lceil 2\delta |T| \rceil$. Let $B' \subseteq B$ be the union of the sets $(T')_{T \,\in\, \mathcal{T}}$. 
        Then, by \cite[Lem.~4.15]{DZ_comparison}, 
    \[ \underline{D}(B') \,\geq\, 2\delta \cdot \underline{D}(B) \,\geq\, 2\delta(1-\delta) \,>\, \delta \,\geq\, \overline{D}(A). \]
    Therefore, $\underline{D}(B', A) \geq \underline{D}(B') - \overline{D}(A) > 0$, and the Borel comparison property of the action $\G \acts X$ implies that $[A] \leq [B']$ in $\mathsf{S}(X,\G,\algebra{B})$.
    
    Fix a Borel equidecomposition of $A$ with a Borel subset of $B'$ and let $f \colon A \to B'$ be the Borel injection sending each point in $A$ to the point in $B'$ it is moved to by that equidecomposition. Note that there is a finite set $S \subset \G$ such that $f(x) \in S \cdot x$ for all $x \in A$. For each $T \in \mathcal{T}$, define
    \[
        T^* \,\defeq\, T \cup f\inv(T) \,=\, T \cup f\inv(T').
    \]
    In other words, we form $T^*$ by adding to $T$ the points in $A$ that are mapped to $T$ by $f$. Let \[
        \mathcal{T}^* \,\defeq\, \set{T^* \,:\, T \in \mathcal{T}}.
    \]
    By construction, $\mathcal{T}^*$ is a Borel family of finite sets covering every point in $X$ exactly once. Moreover, since for each $T \in \mathcal{T}$, we have $T^* \subseteq S^{-1} \cdot T$, and since $\mathcal{T}$ satisfies the finitely many shapes condition, we see that $\mathcal{T}^*$ satisfies the finitely many shapes condition as well.

    Finally, we observe that for each $T \in \mathcal{T}$, $|f^{-1}(T')| \leq |T'| < 2\delta|T| + 1$, and thus 
    \begin{align*}
        \frac{|(K\cdot T^*) \symdif T^*|}{|T^*|} \,&\leq\, \frac{|(K \cdot T) \symdif T| \,+\, |K \cdot f\inv(T')|}{|T|} \\
        &\leq\, \eps \,+\, \frac{|K| (2\delta|T| + 1)}{|T|} \,=\, 3\epsilon \,+\, \frac{|K|}{|T|} \,\leq\, 
        5\epsilon,
    \end{align*}
    where in the last step we use that for all $T \in \mathcal{T}$,
    \[
        |T| \,\geq\, (1-\eps)|L| \,\geq\, \frac{1-\epsilon}{\epsilon}\, |K| \,\geq\, \frac{1}{2\epsilon}\, |K|.
    \]
        Here the first inequality holds because $T$ is $(L,\epsilon)$-F\o{}lner, and the last uses that $\epsilon < 1/2$.
    In conclusion, $\mathcal{T}^*$ is a Borel $(K,5\epsilon)$-F\o{}lner tiling of $X$, and we are done since $K$ and $\epsilon$ are arbitrary.
    \end{scproof}

\section{Problems}\label{sec:problems}

    We end the paper by collecting some natural questions that so far remain unresolved. The following conjecture is restated here for ease of reference:

    \begin{conjcopy}{conj:asi}\label{conj:asi_copy}
        Suppose that a Borel bipartition $(A,B)$ of a locally finite Borel graph $G$ has combinatorial expansion. If $\asi(G) < \infty$, then $G$ has a Borel matching covering $A$. 
    \end{conjcopy}

    The following special case of Conjecture~\ref{conj:asi_copy} is already open:

    \begin{conj}
        Suppose $(A,B)$ is a bipartition of a locally finite Borel graph $G$ of type $(4,3)$. If $\asi(G) < \infty$, then $G$ has a Borel matching covering $A$.
    \end{conj}

    As discussed in \S\ref{subsec:asi_simple}, a powerful general approach to problems in descriptive combinatorics relies on distributed algorithms. 
    Distributed complexity of the matching problem in unbalanced graphs is not yet settled. In particular, the following is open:

    \begin{question}\label{ques:LOCAL_matching}
        Does there exist a randomized \LOCAL algorithm that, given an $n$-vertex graph $G$ of constant maximum degree $\Delta$ with an unbalanced bipartition $(A,B)$, in time $o(\log n)$ finds a matching in $G$ that covers $A$? (The implied constants in the $o(\cdot)$ notation may depend on $\Delta$.)
    \end{question}

    Answering Question~\ref{ques:LOCAL_matching} in either direction will likely require significant new ideas. A positive solution to Question~\ref{ques:LOCAL_matching} would be an extremely strong result; in particular, it would simultaneously include as special cases Theorem~\ref{theo:measurable} and Conjecture~\ref{conj:asi_copy} (which itself implies Theorems~\ref{theo:BM} and~\ref{theo:ends}), thanks to \cite[Thm.~2.14]{BernshteynDistributed}, \cite[Thm.~1.42]{bw_lll}, and Theorem \ref{theo:asi_two_conjectures}. 

    Note that if we only assume that $(A,B)$ has combinatorial expansion, then a \LOCAL algorithm as in Question~\ref{ques:LOCAL_matching} does not exist. Indeed, if it did, then, by \cite[Thm.~2.14(i)]{BernshteynDistributed},
        Theorem~\ref{theo:measurable} would also hold for Borel graphs $G$ of finite maximum degree and their Borel bipartitions $(A,B)$ with combinatorial expansion, which we know is false \cite[Prop.~16.1]{kechris.marks}. 
        On the other hand, 
        the proofs of Theorems~\ref{theo:asi} and~\ref{theo:asi_simple} given in this paper can be adapted to show that the answer to Question~\ref{ques:LOCAL_matching} is ``yes'' in the following special cases:
        \begin{itemize}
            \item when $(A,B)$ is of type $(\high,\low)$ with $\high > C\,\low$, where $C >0$ is a certain universal constant,
            \item when $G$ is $k$-simple and $(A,B)$ is of type $(\high, \low)$, where $\high > (1+\epsilon)\,\low$ for some fixed $k \in \N^+$ and $\epsilon > 0$ and sufficiently large $\low$. 
        \end{itemize}
        The main subroutine required for this is the sublogarithmic-time randomized \LOCAL algorithm for the LLL due to Fischer and Ghaffari \cite{decompose7} (or the even faster, $\operatorname{poly}(\log \log n)$-time  algorithm of Rozho\v{n} and Ghaffari \cite[Cor.~3.9]{decompose6}). The second bullet point also needs the list-coloring algorithm of Chung, Pettie, and Su \cite[\S4.4]{CPS}, which can be sped up to run in $\operatorname{poly}(\log \log n)$ time using the randomized time hierarchy gap theorem of Chang and Pettie \cite[Thm.~6]{CP}. 
        

    There are also some interesting open questions concerning edge-colorings. One of them we have already mentioned:

    \begin{conjcopy}{conj:asi_Shannon}
        If $G$ is a Borel graph with $\Delta(G) < \infty$ and $\asi(G) < \infty$, then
        $
         \displaystyle   \chi'_\mathsf{B}(G) \leq \left\lfloor \frac{3}{2}\,\Delta(G) \right\rfloor
        $.
    \end{conjcopy}

    The version of this conjecture  for \LOCAL algorithms is also open:

    \begin{question}\label{ques:LOCAL_shannon}
        Does there exist a randomized \LOCAL algorithm that, given an $n$-vertex multigraph $G$ of constant maximum degree $\Delta$, in time $o(\log n)$ finds a proper $\lfloor 3\Delta/2 \rfloor$-edge-coloring of $G$? 
    \end{question}

    The fastest currently available randomized \LOCAL algorithm for this problem, due to Dhawan \cite{dhawan2024edge}, 
    runs in time $O((\log n)^2)$. 

    Measurable edge-colorings of multigraphs warrant further investigation. For example, note that if $G = (V,E,\pi)$ is a Borel graph of maximum degree $\Delta \in \N$ with a probability measure $\mu$ on $V$, then we have $\chi'(G) \geq \Delta$ and, by Theorem~\ref{theo:Shannon_meas}, $\chi'_\mu(G) \leq \lfloor 3\Delta/2\rfloor$. In particular, we can write 
    \[\chi'_\mu(G) \,\leq\, \left\lfloor \frac{3}{2}\,\chi'(G) \right\rfloor.\]
    How close is this to being tight? In other words, how large can the gap between $\chi'_\mu(G)$ and $\chi'(G)$ be? Greb\'ik \cite{MCP2} showed that if $G$ is simple, then $\chi'_\mu(G) \leq \Delta + 1$ and hence $\chi'_\mu(G) \leq \chi'(G) + 1$. The former inequality generally fails for multigraphs, but the latter one is open:



    \begin{question}\label{ques:meas_index_diff}
        Does there exist a Borel multigraph $G$ of finite maximum degree with a probability measure $\mu$ such that $\chi'_\mu(G) > \chi'(G) + 1$? Can the difference $\chi'_\mu(G) - \chi'(G)$ be arbitrarily large?
    \end{question} 

    We can also use edge-colorings to formulate potential stronger versions of our results about matchings. For example, we make the following conjecture that extends Theorem~\ref{theo:measurable}:

    \begin{conj}\label{conj:coloring_bipartite}
        Suppose $(A,B)$ is an unbalanced bipartition of a 
        Borel multigraph $G = (V,E,\pi)$ of maximum degree $\Delta \in \N$. Then $\chi'_\mu(G) \leq \Delta$ for any probability measure $\mu$ on $V$.
    \end{conj}

    To see that Theorem~\ref{theo:measurable} is a consequence of Conjecture~\ref{conj:coloring_bipartite}, let $G = (V,E,\pi)$ be a locally finite Borel graph with a bipartition $(A,B)$ of type $(\high,\low)$, where $\high > \low$. Using Lemma~\ref{lemma:reg_subgraph}, we pass to a spanning subgraph of $G$ to arrange that every vertex in $A$ has degree exactly $\high$. Then $\Delta(G) = \high$, so, given any probability measure $\mu$ on $V$, Conjecture~\ref{conj:coloring_bipartite} yields disjoint Borel matchings $M_1$, \ldots, $M_\high$ in $G$ such that for $\mu$-almost all $x \in V$, the set $M_1 \cup \ldots \cup M_\high$ includes all edges incident to $x$. This implies that each $M_i$ must cover $\mu$-almost every vertex in $A$, as desired.

    Note that a bipartite multigraph $G$ of maximum degree $\Delta \in \N$ satisfies $\chi'(G) \leq \Delta$ by K\H{o}nig's theorem \cites[Prop.~5.3.1]{Diestel}{konig_paper}, so there are no combinatorial obstructions to Conjecture~\ref{conj:coloring_bipartite}. For simple graphs, the following even more general conjecture might hold:

    \begin{conj}\label{conj:F}
        Let $G = (V,E,\pi)$ be a simple Borel graph of maximum degree $\Delta \in \N$. Suppose no two vertices of degree $\Delta$ in $G$ are adjacent. Then $\chi'_\mu(G) \leq \Delta$ for any probability measure $\mu$ on $V$.
    \end{conj}

    Fournier \cite{fournier} showed that under the assumptions of Conjecture~\ref{conj:F}, we have $\chi'(G) \leq \Delta$ \ep{this is a special case of Vizing's Adjacency Lemma \cites[Lem.~3.11]{cranston}{VizLemma}}. Conjecture~\ref{conj:F} implies Conjecture~\ref{conj:coloring_bipartite} for simple graphs, because if $(A,B)$ is an unbalanced bipartition of $G$, all vertices of maximum degree in $G$ must be in $A$, and hence they cannot be adjacent to each other. Note that, by Corollary~\ref{corl:3}\ref{item:3_meas}, Conjecture~\ref{conj:F} is true for $\Delta = 3$.

    In forthcoming work \cite{Sequel}, we make partial progress toward Conjectures~\ref{conj:coloring_bipartite} and \ref{conj:F} using probabilistic techniques similar to the ones in the proof of Theorem~\ref{theo:measurable}. Specifically, we prove both conjectures for multigraphs whose underlying simple graphs are acyclic (it turns out that in this special case, it is not necessary to assume that $G$ is simple in Conjecture~\ref{conj:F}). We also show that bipartite Borel multigraphs $G$ of maximum degree $\Delta \in \N$ satisfy $\chi'_\mu(G) \leq \Delta + 1$. 


    Our results in equidecomposition theory also raise some interesting open questions. For example, we find the following problem quite intriguing:

    \begin{question}\label{ques:paradox}
        Let $\G \acts X$ be a Borel action of a group $\G$ on a standard Borel space $X$ with Borel $\sigma$-algebra $\algebra{B}$. Using the terminology of \cite{TypeSemigroupGroupoids}, does the type semigroup $\mathsf{S}(X,\G,\algebra{B})$ always \emphd{have plain paradoxes}? That is, if an element $\alpha \in \mathsf{S}(X,\G,\algebra{B})$ satisfies $n\alpha \leq m\alpha$ in $\mathsf{S}(X,\G,\algebra{B})$ for some $n > m \in \N$, must we have $2\alpha \leq \alpha$ in $\mathsf{S}(X,\G,\algebra{B})$ as well?
    \end{question}

    If a preordered commutative monoid is almost unperforated, then it does have plain paradoxes \cite[Rmk.~2.20]{TypeSemigroupGroupoids}. In particular, because of parts \ref{item:unperf_meas} and \ref{item:unperf_BM} of Theorem~\ref{theo:unperf}, we cannot use measure or Baire category to obtain a counterexample to Question~\ref{ques:paradox}. The only other available tool seems to be Marks's Borel determinacy technique. Indeed, this technique underlies the proof of Proposition~\ref{prop:no_borel_unperf}, which shows that $\mathsf{S}(X,\G,\algebra{B})$ need not be almost unperforated. However, the example constructed there is not helpful for Question~\ref{ques:paradox}, as, essentially by design, the structures arising in Marks's Theorem~\ref{theo:grab} avoid ``paradoxical'' relations of the form $n \alpha \leq m \alpha$ with $n > m$ (see \cite[\S3.3]{HypBorCom} for an exploration of this phenomenon). Thus, Question~\ref{ques:paradox} appears challenging and may require some significant new insight.

    Lastly, let us consider Borel actions of finitely generated amenable groups. It is not known whether all such actions have the Borel tiling property, or, equivalently, the Borel comparison property \cite[Prob.~3.2]{marks2023measurable}. In view of Proposition~\ref{prop:unp_to_comp}, the positive answer to the following question would be even stronger: 

    \begin{question}\label{ques:amenable_unp}
        Let $\G \acts X$ be a Borel action of a finitely generated amenable group $\G$ on a standard Borel space $X$ with Borel $\sigma$-algebra $\algebra{B}$. Must the type semigroup $\mathsf{S}(X,\G,\algebra{B})$ be almost unperforated?
    \end{question}

    It might even be that all Borel actions of finitely generated amenable groups have finite asymptotic separation index \cite[Q.~3.1.18]{WeilacherThesis}, in which case the answer to Question~\ref{ques:amenable_unp} would be ``yes'' by Theorem~\ref{theo:unperf}.
    



    \subsection*{Acknowledgments}
    FW thanks Forte Shinko for asking interesting questions about measurable edge-chromatic numbers.

    \subsection*{Funding}
    
    This material is based upon work partially supported by the Alfred P. Sloan Foundation, the Simons Foundation under grant 376201, and the National Science Foundation under grants DMS-2528522 and DMS-2402064. Any opinions, findings, and conclusions or recommendations expressed in this material are those of the authors and do not necessarily reflect the views of the funding agencies.

\phantomsection
\addcontentsline{toc}{section}{\hspace*{-8pt}References}

    \printbibliography

\end{document}